\documentclass[12pt]{amsart}
\usepackage[margin=1.0in]{geometry}
\usepackage{amsmath,amssymb,amsfonts,graphicx,color,fancyhdr,psfrag,comment,enumerate}
\usepackage[latin1]{inputenc}
\usepackage[hyperpageref]{backref}
\usepackage[colorlinks=true, pdfstartview=FitV, linkcolor=blue, citecolor=blue, urlcolor=blue]{hyperref}
\allowdisplaybreaks

\usepackage{epstopdf}
\usepackage{mathrsfs}
\usepackage{epsfig}
\usepackage{hyperref}
\usepackage{ifthen}
\usepackage{float}
\usepackage{enumerate}
\usepackage{mathtools}
\usepackage[bottom]{footmisc}
\usepackage{tikz}
\usepackage{soul}
\usetikzlibrary{matrix,shapes,arrows,positioning,chains}
\usepackage{amssymb,amsmath,amsfonts}
\usepackage{bbm}

\numberwithin{equation}{section}

\newtheorem{theorem}{Theorem}[section]
\newtheorem{proposition}[theorem]{Proposition}
\newtheorem{lemma}[theorem]{Lemma}

\newtheorem{remark}[theorem]{Remark}
\newtheorem{corollary}[theorem]{Corollary}

\begin{document}
\title[Bilinear maximal functions on the Heisenberg group]{Bilinear maximal functions associated with some degenerate surfaces on the Heisenberg group}

\author[S. Bagchi, R. Garg and S. K. Rana]
{Sayan Bagchi \and Rahul Garg \and Surya Kanta Rana}

\address[S. Bagchi]{Department of Mathematics and Statistics, Indian Institute of Science Education and Research Kolkata, Mohanpur--741246, West Bengal, India.}
\email{sayansamrat@gmail.com}

\address[R. Garg]{Department of Mathematics, Indian Institute of Science Education and Research Bhopal, Bhopal--462066, Madhya Pradesh, India.}
\email{rahulgarg@iiserb.ac.in}

\address[S. K. Rana]{Department of Mathematics, Indian Institute of Science Education and Research Bhopal, Bhopal--462066, Madhya Pradesh, India.}
\email{surya22@iiserb.ac.in}



\begin{abstract}
Motivated by the works of Lee--Shuin and Borges--Foster on the Euclidean space, in this article we consider bilinear averages on the Heisenberg group $\mathbb{H}^n$ corresponding to some degenerate surfaces associated with the Kor\'anyi norm as well as Nevo--Thangavelu type averaging structures, and study the boundedness of the their full and lacunary maximal functions on various $L^{p}(\mathbb{H}^n) \times L^{q}(\mathbb{H}^n)$ spaces, for $n \geq 2$. For the full maximal functions, we obtain optimal range except possibly at the end points.  
\end{abstract}
\maketitle

\tableofcontents

\section{Introduction} \label{sec:intro}
The study of maximal functions is of fundamental importance in Harmonic Analysis. In his seminal work, Stein \cite{Stein-Elias-Maximal-functions-Spherical-means} investigated the Euclidean spherical maximal function $M^{\mathrm{Euc}}$ defined by
\[
M^{\mathrm{Euc}}f(x)
= \sup_{r>0} \left| A^{\mathrm{Euc}}_{r}f(x) \right|
= \sup_{r>0}\left| \int_{|y|=1} f(x-ry)\, d\sigma_{n-1}(y)\right|, 
\]
and showed that $M^{\mathrm{Euc}}$ is bounded on $L^p(\mathbb{R}^n)$, $n\ge 3$, if and only if $
p>\frac{n}{n-1}.$
Here $\sigma_{n-1}$ denotes the standard surface measure on the unit sphere $ S_1=\{x\in\mathbb{R}^n : |x|=1\}$. 
The case $n=2$ was settled later by Bourgain \cite{Bourgain-Averages-in-the-plane-over-convex-curves-and-maximal-operators}. 

\medskip
Another very important class of maximal operators corresponds to taking supremum over dyadic dilations of suitable subsets $E$ of $[1,2]$. In particular, when $E = \{1\}$, this leads to what is well known as the lacunary spherical maximal operator $M^\mathrm{Euc}_{lac}$, defined by 
\[
M^{\mathrm{Euc}}_{{lac}} f(x)
=
\sup_{j\in \mathbb{Z}}
\left|
\int_{|y|=1} f(x-2^jy)\, d\sigma_{n-1}(y)
\right|.
\]

Clearly, $M^{\mathrm{Euc}}_{lac}$ is smaller than $M^{\mathrm{Euc}}$ and one would expect it to have nice mapping properties on a wider class of functions. That this is indeed the case, was first shown by Calder\'on \cite{Calderon-lacunary-spherical-means}. He showed that $M^{\mathrm{Euc}}_{lac}$ is bounded on $L^p(\mathbb{R}^n)$ for all $1<p\leq \infty$ for $n\ge 2$. More recently, Lacey \cite{Lacey-Sparse-bounds-for-spherical-maximal-functions} revisited the theory of spherical maximal functions, introducing a novel approach to derive sparse domination results.

\medskip
We now turn our attention to the bilinear spherical maximal function, defined by
\[
\mathcal{M}^{\mathrm{Euc}}(f,g)(x)
=
\sup_{t>0}
\left|\int_{S^{2n-1}}
f(x-ty_1)\,g(x-ty_2)\, d\sigma_{2n-1}(y_1,y_2)\right|,
\]
where $\sigma_{2n-1}$ denotes the standard surface measure on the unit sphere $S^{2n-1}$ in $\mathbb{R}^{2n}$. 

\medskip 
Geba et al. \cite{Geba-Greenleaf-multilinear-operators} were the first to study the operator $\mathcal{M}^{\mathrm{Euc}}$, and their results were further improved in \cite{Loukas-Grafakos-bilinear-spherical-maximal-function, Heo-Hong-Improved-bounds-for-the-bilinear-spherical-maximal-function}. In \cite{Maximal-estimates-for-the-bilinear-spherical-averages}, Jeong--Lee proved the sharp range for any $n\ge 2$, showing that with $1\le p,q\le\infty$, $0<r\le\infty$ with $\frac{1}{p}+\frac{1}{q}=\frac{1}{r}$, the operator $\mathcal{M}^{\mathrm{Euc}}$ is $L^p(\mathbb{R}^n)\times L^q(\mathbb{R}^n) \to L^r(\mathbb{R}^n)
$ bounded if and only if $r>\frac{2n}{2n-1}$, except for the case $(p,q) = (1,\infty)$ or $(p,q) = (\infty,1)$, where they obtained the weak-type bounds. Recently, Borges--Foster \cite{Borges-Tainara-Foster-bounds-for-lacunary-bilinear-spherical-maximal-functions} considered the lacunary case for $n\ge 2$ and showed that the operator $\mathcal{M}^{\mathrm{Euc}}_{lac}$ is bounded from $L^p(\mathbb{R}^n)\times L^q(\mathbb{R}^n) \to L^r(\mathbb{R}^n)$ whenever $p,q\in (1,\infty]$, $r\in (0,\infty)$ and $ \frac{1}{p}+\frac{1}{q}=\frac{1}{r}$.

\medskip 
Moving on, the boundedness properties of  bilinear (full) maximal functions associated with some degenerate surfaces were first considered in \cite{Lee_Shuin}. To be more precise, for $n \geq 2$, given $\mathbf{a}=(\mathbf{a}_1,\mathbf{a}_2)\in [1,\infty)^2$, they considered biparametric spherical averages $\mathcal{A}_{\mathbf{s}}^{\mathrm{Euc},\,\mathbf{a}}$, for $\mathbf{s}=(s_1,s_2)\in (0,\infty)^2$, defined over the degenerate surface $S_2^\mathbf{a} :=\{(y,z)\in \mathbb{R}^{2n}: |y|^{\mathbf{a}_1}+|z|^{\mathbf{a}_2}=1\}$, by 
$$ \mathcal{A}_{\mathbf{s}}^{\mathrm{Euc},\,\mathbf{a}}(f,g) (x) := \int_{{S}_2^\mathbf{a}} f(x-s_1y) \, g(x-s_2z)\, d\mu_\mathrm{Euc}^\mathbf{a}(y,z),$$ 
where $d\mu_\mathrm{Euc}^\mathbf{a}$ is the natural surface measure on $S_2^\mathbf{a}$. They proved the sharp strong-type boundedness result for the associated full maximal function $\mathcal{M}^\mathbf{a}(f,g)(x):=\sup\limits_{\mathbf{s}\in (0,\infty)^2}|\mathcal{A}_{{\mathbf{s}}}^{\mathrm{Euc},\,\mathbf{a}}(f,g)(x)|$ except at the boundary. The lacunary counterpart $\mathcal{M}_{lac}^{\mathbf a}$, defined by $\mathcal{M}_{lac}^{\mathbf a}(f,g)(x) : = \sup\limits_{ j \in \mathbb Z} |\mathcal{A}_{{2^j}}^{\mathrm{Euc},\,\mathbf{a}}(f,g)(x)|$, was studied in \cite{Borges-Tainara-Foster-bounds-for-lacunary-bilinear-spherical-maximal-functions}. Let $\mathcal R_{lac}^{\mathbf a}$ denote the H\"older boundedness region for the operator $\mathcal{M}_{lac}^{\mathbf a}$. It was shown in \cite{Borges-Tainara-Foster-bounds-for-lacunary-bilinear-spherical-maximal-functions} that if $\mathcal R^{\mathbf a}_{1}$ denotes the H\"older boundedness region for the single-scale averages $\mathcal{A}_{1}^{\mathrm{Euc},\,\mathbf{a}}$, then $\operatorname{int}(\mathcal R^{\mathbf a}_{1}) \subseteq \mathcal R_{lac}^{\mathbf a}$ whenever $1\leq a_1,a_2\leq 2n$. 

\medskip
Our objective in this article is to study analogous results on the Heisenberg group $\mathbb{H}^n$ corresponding to similar degenerate surfaces associated with both, the Kor\'anyi norm as well as Nevo--Thangavelu type averaging structures. Let $\mathbb{H}^n = \mathbb{C}^n \times \mathbb{R}$ be the $(2n+1)$-dimensional
Heisenberg group equipped with the group law
\[
(z,t) \cdot (w,s) = (z+w,\, t+s+\tfrac{1}{2}\Im(z\cdot \overline{w})).
\]
The Haar measure on $\mathbb{H}^n$ is simply the Lebesgue measure $dz \, dt$, and we shall usually denote it by $dg$. 

\medskip 
Consider the Kor\'anyi norm $| \cdot |_K$ on $\mathbb{H}^n$ which is given by as follows. Given $g = (z,t) \in \mathbb{H}^n$, we define $|g|_K = |(z,t)|_K =(|z|^4+t^2)^{1/4}$. Let $B_K(g, r) = \{g' \in \mathbb{H}^n : |g^{-1} g'|_K < r\}$ denote the Kor\'anyi ball of radius $r$ centered at $g \in \mathbb{H}^n$. Also, let $S_{K,1}$ denote the unit Kor\'anyi sphere $S_{K,1}=\{g \in \mathbb{H}^n: |g|_K=1\}$ in $\mathbb{H}^n$. One can define a spherical maximal function $M_{\mathbb{H}^n}$ with respect to the Kor\'anyi norm $| \cdot |_K$ on $\mathbb{H}^n$ as follows. There exists unique Radon measure $\sigma^{2n}_K$ supported on $S_{K,1}$ with which we have the polar decomposition of the Haar measure $dg$ on $\mathbb{H}^n$: 
\begin{equation} \label{polar-decompostion-withrespectto-the-Koranyi-norm}
\int_{\mathbb{H}^n} f(g) \, dg = \int_0^\infty \int_{S_{K,1}} f(\tilde{\delta_r}g) \, r^{2n+1} \, d\sigma^{2n}_K(g) \, dr, 
\end{equation} 
where $\tilde{\delta}_r$ denotes the non-isotropic dilation $\tilde{\delta}_r (z,t) = (rz, r^2 t)$ on $\mathbb{H}^n$. We will discuss more about $\sigma^{2n}_K$ in Subsection \ref{exp-surface-measure-on-the-linear-Koranyi-sphere}. One defines the spherical maximal function $M_{\mathbb{H}^n}$ by
\begin{align} 
\label{def:linear-Koranyi-sph-maximal}
M_{\mathbb{H}^n}F(g) = \sup_{r>0} \left| A_r^K F(g) \right|, 
\end{align}
corresponding to the spherical averages $A_r^K$ which are given by 
\begin{align} 
\label{def:linear-Koranyi-averages}
A_r^K F(g) = \int_{|g'|_K=1} F(g\cdot\tilde{\delta}_r g'^{-1}) \, d\sigma^{2n}_K(g'). 
\end{align} 

In \cite{Cowling-Littlewood-paley-stein-theory}, M. Cowling initiated the study of the operator $M_{\mathbb{H}^n}$, proving the sharp result that it is bounded on $L^p(\mathbb{H}^n)$ if and only if $p>\frac{2n+1}{2n}$ for any ${n\geq2}$. More recently, the lacunary spherical maximal function $M_{{\mathbb{H}^n},{lac}}$ associated with the Kor\'anyi sphere $\Sigma_{K,1}$, defined by $M_{{\mathbb{H}^n},{lac}}F(g) = \sup\limits_{k\in \mathbb{Z}} | A_{2^k}^K F(g)| $, has been studied in \cite{On-the-lacubary-spherical-maximal-function-on-the-Heisenberg-group,Rajula-Srivastava-Koranyi_sphere, Aswin-Hickman}. 

\medskip
On the other hand, let $\sigma$ denote the rotation invariant surface measure on the $2$ co-dimensional unit sphere $S^{2n-1}$ in $\mathbb{R}^{2n} \times \{0\}$. In \cite{Nevo-Thangavelu-pointwise-ergodic-theorems-for-the-radial-averages-on-the-Heisenberg-groups}, Nevo--Thangavelu considered the following spherical means on $\mathbb{H}^n$: 
\begin{align} 
\label{def:linear-NT-averages}
A^{NT}_r F(g)
= \int_{|w|=1} F\left(z-rw,\, t-\tfrac{r}{2}\Im(z\cdot \overline{w})\right)\, d\sigma(w),
\end{align}
where $g = (z,t)$, and studied the maximal function 
\begin{align} 
\label{def:linear-NT-sph-maximal} 
M^{NT}_{\mathbb{H}^n} F(g) = \sup_{r>0} |A^{NT}_r F(g)|. 
\end{align}
Later, improving upon \cite{Nevo-Thangavelu-pointwise-ergodic-theorems-for-the-radial-averages-on-the-Heisenberg-groups}, M\"uller--Seeger \cite{Muller-Seeger-spherical-maximal-function-on-twostep-nilpotent-group} and Narayanan--Thangavelu \cite{Naryanan-Thangavelu-optimal-theorem-for-the-spherical-maximal-operator-on-the-Heisenberg-group}, independently, established the sharp result which states that the operator $M^{NT}_{\mathbb{H}^n}$ is bounded on $L^p(\mathbb{H}^n)$, $n\ge 2$, if and only if $p>\frac{2n}{2n-1}.$ In \cite{Sayan-Sourav-Luz-Thangavelu-on-the-maximal-function-associated-to-the-spherical-means-on-the-Heisenberg-groups}, Bagchi et al.\ proved that the lacunary spherical maximal function $M^{NT}_{\mathbb{H}^n,lac}$, defined 
by $M^{NT}_{\mathbb{H}^n,lac} F(g) = \sup_{j\in\mathbb{Z}} \left| A^{NT}_{2^j}F(g) \right|$, is bounded on $L^p(\mathbb{H}^n)$ for any $n\geq 2$ and $1<p \leq \infty$. For further developments, including the $L^p$-improving properties and sparse bounds, we refer to \cite{Roos-seeger-Rajula-lebesgue-space-estimates-for-spherical-maximal-functions,Roos-seeger-Rajula-Restricted-dilation-sets}.

\medskip
Moving on, consider the product Heisenberg group 
\[
\mathbb{H}^n \times \mathbb{H}^n
=
(\mathbb{C}^n \times \mathbb{R}) \times (\mathbb{C}^n \times \mathbb{R})= \mathbb{R}^{4n+2},
\]
equipped with the product group law
\[
(g_1,g_2)\cdot (h_1,h_2)
=
(g_1\cdot h_1,\; g_2\cdot h_2). 
\]
The Haar measure on $\mathbb H^n\times\mathbb H^n$ is nothing but the Lebesgue measure $d(g_1,g_2) = dg_1 \, dg_2$. 

\medskip 
Recently, Ghosh--Singh \cite{ghosh2026bilinearsphericalmaximalfunction} defined and studied Nevo--Thangavelu type bilinear maximal operators (full and lacunary) on $\mathbb{H}^n\times \mathbb{H}^n.$ Let $\sigma$ denote the surface measure on $S^{4n-1}.$ Defining the full maximal function $\mathcal{M}^{NT}_{full}$ by 
\begin{align}
\label{def:bilinear-NT-sph-maximal} 
\mathcal{M}^{NT}_{full}(F,G)(g) = \sup \limits_{t>0} \left| \int_{S^{4n-1}} F(g \cdot \tilde{\delta}_{t}(z_1,0)^{-1}) \, G(g \cdot \tilde{\delta}_{t} (z_2,0)^{-1}) \, d\sigma(z_1,z_2) \right|,
\end{align}
they proved a sharp boundedness result except some boundaries for any $n \geq 2$. More precisely, they showed that $\mathcal{ M}^{NT}_{full}$ boundedly maps $L^p \times L^q \rightarrow L^r$ if $\frac{1}{p}+\frac{1}{q}=\frac{1}{r}$ whenever $\left(\frac{1}{p}, \frac{1}{q}\right)\in \mathcal{R}^{NT}$, where 
\begin{align} 
\label{range_for_the_Bilinear_NT}
\mathcal{R}^{NT} :=\left\{\left(\frac{1}{p},\frac{1}{q}\right)\in [0,1)^2 : \frac{1}{p}+\frac{1}{q}<\frac{4n-1}{2n}\right\}. 
\end{align} 

%%%%%%%%%%%%%%%%%%%%%%%%%%%%%%%%%%%%%%%%%%%%%%%%%%%%%%%%%%%%%%%%%%%%%%%%%%%%%%%%%%%%%%%%%%%%%%%%%%%%%%%%%%%%%%%%%%%%%%%%%%%%%%%%%%
\begin{figure}[ht]
\centering

\begin{tikzpicture}[scale=3]

    %================================================
    % Parameter
    %================================================
    \def\n{2}
     \def\a{0.65}
    %================================================
    % Coordinates
    %================================================
    \coordinate (O) at (0,0);
    \coordinate (R) at (0,1);
    \coordinate (P) at (\a,1);
    \coordinate (Q) at (1,\a);
    \coordinate (S) at (1,0);

    %================================================
    % Shaded region P
    %================================================
    \fill[gray!20]
        (O) -- (R) -- (P) -- (Q) -- (S) -- cycle;

    %================================================
    % Auxiliary horizontal and vertical lines
    %================================================
  %  \draw[gray!45, thin]
      %  (P) -- (\a,0);

   % \draw[gray!45, thin]
       % (0,\a) -- (Q);

    %================================================
    % Boundary of the pentagon
    %================================================
    \draw[black!45!black, thick]
        (O) -- (R);

    \draw[black!45!black, thick]
      (S) -- (O);

    %================================================
    % Dotted triangle P-(1,1)-Q
    %================================================
    \draw[black!45, dotted, thick]
        (P) -- (1,1);

    \draw[black!45, dotted, thick]
        (1,1) -- (Q);

    % Dotted slanted edge PQ
    \draw[black!45!black, thick, dotted]
        (P) -- (Q);
            \draw[black!45!black, thick, dotted]
        (R) -- (P);
            \draw[black!45!black, thick, dotted]
        (S) -- (Q);

    %================================================
    % Coordinate axes
    %================================================
    \draw[->, thick]
        (-0.08,0) -- (1.20,0)
        node[right] {$\frac{1}{p}$};

    \draw[->, thick]
        (0,-0.08) -- (0,1.18)
        node[above] {$\frac{1}{q}$};

    %================================================
    % Points
    %================================================

    % O
    \fill (O) circle (0.4pt);

    % Open circles
    \draw[fill=white, thick] (R) circle (0.4pt);
   % \draw[fill=white, thick] (P) circle (0.4pt);
   % \draw[fill=white, thick] (Q) circle (0.4pt);
    \draw[fill=white, thick] (S) circle (0.4pt);

    %================================================
    % Labels
    %================================================

    % O
    \node[above left=2pt] at (O) {$O$};

    % H
    \node[above right=1pt] at (R) {$R$};
    \node[left=5pt] at (R) {$(0,1)$};

    % P
    \node[above right=1pt and -4pt] at (P)
    {$P=\left(\frac{2n-1}{2n},1\right)$};

    % Q
    \node[right=2pt] at (Q)
        {$Q=\left(1,\frac{2n-1}{2n}\right)$};

    % G
    \node[below=0.02pt] at (S) {$S$};
    \node[below=8pt] at (S) {$(1,0)$};

    %================================================
    % Region label
    %================================================

    \node at (0.45,0.45) {$\mathcal{R}^{NT} $};

\end{tikzpicture}

\caption{Boundedness region of bilinear Nevo--Thangavelu type maximal operator $\mathcal{ M}^{NT}_{full}$: 
the union of the open pentagon (with vertices $O, R, P, Q, S$) together with the line segments $[O,R)$ and  $[O,S)$.}
\label{figure-bilinear-NT-full}

\end{figure}
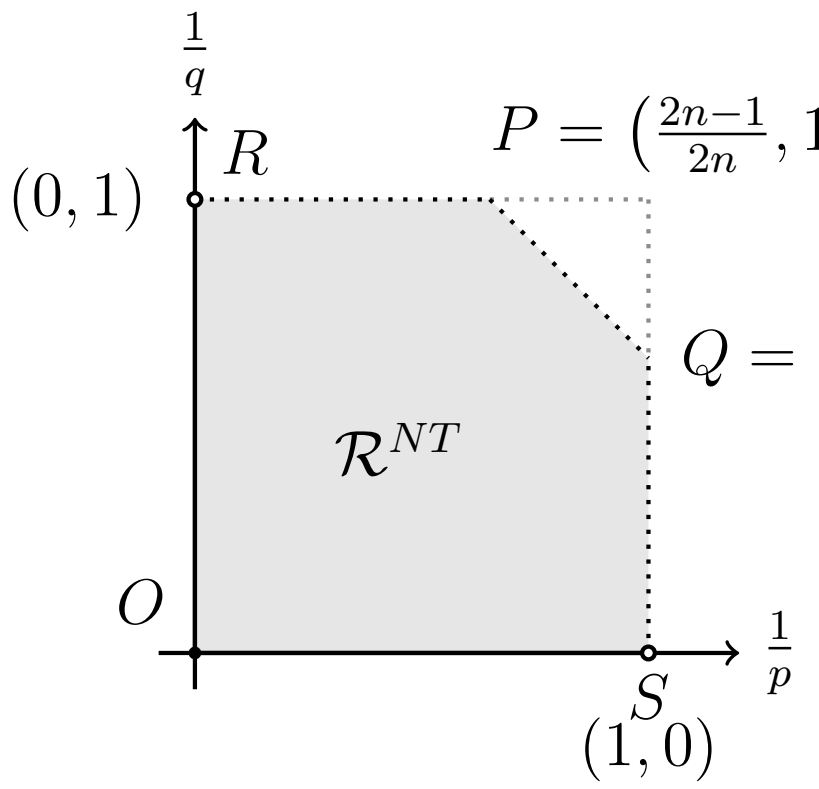

%%%%%%%%%%%%%%%%%%%%%%%%%%%%%%%%%%%%%%%%%%%%%%%%%%%%%%%%%%%%%%%%%%%%%%%%%%
Ghosh--Singh \cite{ghosh2026bilinearsphericalmaximalfunction} also studied the boundedness of the lacunary maximal function
$$ \mathcal{M}^{NT}_{lac}(F,G)(g) = \sup_{k\in \mathbb{Z}} \left| \int _{S^{4n-1}} F(g \cdot \tilde{\delta}_{2^k}(z_1,0)^{-1}) \, G(g \cdot \tilde{\delta}_{2^k} (z_2,0)^{-1}) \, d\sigma (z_1,z_2) \right|.
$$
They established an $L^p\times L^q$-improving region for the corresponding single-scale averaging operator, and showed that $\mathcal{M}^{NT}_{lac}$ is bounded whenever $\left(\frac{1}{p},\frac{1}{q}\right)$ belongs to the interior of that region. 

\medskip
Motivated by the above mentioned results, as mentioned earlier, in this article we study the bilinear maximal functions associated to some degenerate surfaces on $\mathbb{H}^n \times \mathbb{H}^n$, which can be seen as a generalisation of bilinear spherical maximal functions $\mathcal{M}^{NT}_{full}$ and $\mathcal{M}^{NT}_{lac}$ on these surfaces. But, before describing these surfaces and our results in this direction, let us also discuss another setup of our interest, namely, the product Kor\'anyi case. Define the Kor\'anyi norm $| \cdot |_{K\times K}$ on
$\mathbb{H}^{n} \times \mathbb{H}^{n}$ by
\[
|(g_1,g_2)|_{K\times K}
=
\left(|g_1|_K^4 + |g_2|_K^4\right)^{1/4}. 
\]
The product $r-$radius Kor\'anyi sphere here is $\Sigma_{K,r}
= \left\{(g_1,g_2)\in \mathbb{H}^{n} \times \mathbb{H}^{n}
:\ |(g_1,g_2)|_{K\times K} = r \right\}.$
Note also that on $\mathbb{H}^{n} \times \mathbb{H}^{n}$, we have the one parameter family of non-isotropic dilations, $\delta_r : \mathbb{H}^n \times \mathbb{H}^n \rightarrow \mathbb{H}^n \times \mathbb{H}^n$, given by 
\[
\delta_r(g_1,g_2)=(\tilde{\delta_r}(g_1),\tilde{\delta_r}(g_2)).
\]
It follows that $| \cdot |_{K\times K}$ is homogeneous of degree $1$ with respect to $\delta_r $, that is, $|\delta_r(g_1,g_2)|_{K\times K} = r \, |(g_1,g_2)|_{K\times K}$ for any $r>0$. Clearly, $\Sigma_{K,r} = \delta_r \Sigma_{K,1}$. Also, the Jacobian of $\delta_r$ is $r^{4n+4}$, which gives the homogeneous dimension of the group $\mathbb{H}^n \times \mathbb{H}^n$ to be $4n+4$. 

\medskip
It turns out that similar to the linear case, that is, relation \eqref{polar-decompostion-withrespectto-the-Koranyi-norm}, we have an explicit polar decomposition on $\mathbb{H}^n \times \mathbb{H}^n$ as well. More precisely, we have a Radon measure $d\mu^{4n+1}_K$ on the product Kor\'anyi sphere such that 
\begin{equation} 
\label{expres-polar-decoposition-producttt}
\int_{\mathbb{H}^n  \times \mathbb{H}^n} f(g_1,g_2)\,dg_1\,dg_2=\int_0^\infty \left( \int_{\Sigma_{K,1}} f(\delta_r(g_1,g_2)) \, d\mathcal{\mu}^{4n+1}_K(g_1,g_2)\right) r^{4n+3}\,dr. 
\end{equation}
We will discuss the measure $d\mathcal{\mu}^{4n+1}_K$ in details in Subsection \ref{exp-surface-measure-on-the-product-Koranyi-sphere}. We consider the bilinear spherical averages $\mathcal{A}^K_t$ on $\Sigma_{K,1}$ defined by 
\begin{equation}\label{expression-of-spherical-averages-biliear-case}
\mathcal{A}^K_t(F,G)(g) := \int \limits_{\Sigma_{K,1}} F(g \cdot \tilde{\delta_t}{g_1}^{-1})\, G(g \cdot \tilde{\delta_t}{g_2}^{-1}) \, d\mu^{4n+1}_K (g_1,g_2), 
\end{equation} 
and define the corresponding bilinear spherical maximal function 
\begin{equation}\label{definition-of-the-bilinear-spherical-maximal-function}
\mathcal{M}^K_{{full}}(F,G)(g) := \sup_{t>0} \left| \mathcal{A}^K_t(F,G)(g) \right|, 
\end{equation}
and the lacunary one 
\begin{equation}\label{definition-of-lacunary-bilinear-spherical-maximal-function}
\mathcal{M}^K_{lac}(F,G)(g) := \sup_{l\in \mathbb{Z}} \left| \mathcal{A}^K_{2^l}(F,G)(g) \right|.
\end{equation}

More generally, let us consider the following degenerate surfaces associated with the Kor\'anyi norm. Given $(\mathbf{a}_1,\mathbf{a}_2) \in (0,\infty)^2,$ consider the surface 
\begin{align} \label{expression_of_the_non_dengenerate_surface_of_product_Koranyi_type}
\mathcal{S}_K^{\mathbf{a}} := \{(g_1,g_2)\in \mathbb{H}^n\times \mathbb{H}^n : |g_1|_K^{\mathbf{a}_1}\,+ |g_2|_K^{\mathbf{a}_2}=1\}. 
\end{align}
Clearly, the product Kor\'anyi sphere discussed above corresponds to the case when $(\mathbf{a}_1,\mathbf{a}_2) = (4,4)$. Whenever $\mathbf{a}_1, \, \mathbf{a}_2 \geq 2$, we have a natural surface measure $d\mu_K^\mathbf{a}$ on each of $\mathcal{S}_K^{\mathbf{a}}$ whose description is motivated by that of $d\mathcal{\mu}^{4n+1}_K$ (see Subsection \ref{exp-surface-measure-on-the-degenerate-product-Koranyi-sphere}). Given $\mathbf{s}=(s_1,s_2)\in (0,\infty)^2,$ let us consider the biparametric average $\mathcal{A}^{\mathbf{a},K}_\mathbf{s}$ on $\mathcal{S}_K^\mathbf{a}$ given by 
\begin{align} \label{expression_of_averaging_operator_of_Koranyi_type_associated_to_degenerate_surface}
\mathcal{A}^{\mathbf{a},K}_\mathbf{s}(F,G)(g) := \int_{\mathcal{S}_K^\mathbf{a}} F(g\cdot\tilde{\delta}_{s_1}g_1^{-1}) \, G(g\cdot\tilde{\delta}_{s_2}g_2^{-1}) \, d\mu_K^\mathbf{a}(g_1,g_2). 
\end{align}
The corresponding full maximal function is given by 
\begin{align} 
\label{definition_of_full_maximal_degenerate_Koranyi_case} 
\mathcal{M}^{\mathbf{a},K}_{full}(F,G)(g):=\sup\limits_{\mathbf{s}\in (0,\infty)^2}\left| \mathcal{A}^{\mathbf{a},K}_{\mathbf{s}}(F,G)(g) \right|,
\end{align} and the lacunary one by 
\begin{align}\label{definition_of_lacunary_maximal_degenerate_Koranyi_case}
\mathcal{M}^{\mathbf{a},K}_{lac}(F,G)(g):=\sup\limits_{l\in \mathbb{Z}}\left| \mathcal{A}^{\mathbf{a},K}_{2^l}(F,G)(g) \right|.
\end{align}

We are now ready to state our first main result concerning the boundedness of $\mathcal{M}^{\mathbf{a},K}_{{full}}$. For the same, let us begin with illustrating various regions of boundedness. 
\medskip

Writing $l_+=\max\{0,l\}$, let us set
\begin{align} 
\mathcal{R}^K &:= \left\{ \left( \frac{1}{p}, \frac{1}{q} \right) \in [0,1)^2 : \frac{1}{p} + \frac{1}{q} < \frac{4n+1}{2n+1} \right\}, 
\label{description-of-R-range} \\ 
\mathcal{R}^{\mathbf{a},K}_{full} &:= \left\{ \left( \frac{1}{p},\frac{1}{q} \right) \in \mathcal{R}^K : \frac{1}{p}<1-\frac{(\mathbf{a}_2-(2n+2))_+}{\mathbf{a_2}(2n+1)}, \, \frac{1}{q}<1-\frac{(\mathbf{a}_1-(2n+2))_+}{\mathbf{a_1}(2n+1)} \right\}. 
\label{Domain_for_the_Product_Koranyi_type_non_degenerate}
\end{align} 

Depending on whether $\mathbf{a}_1$ or $\mathbf{a}_2$ is less than or greater than the homogeneous dimension $2n+2$ of $\mathbb{H}^n$, three cases arise. We illustrate the same below. 

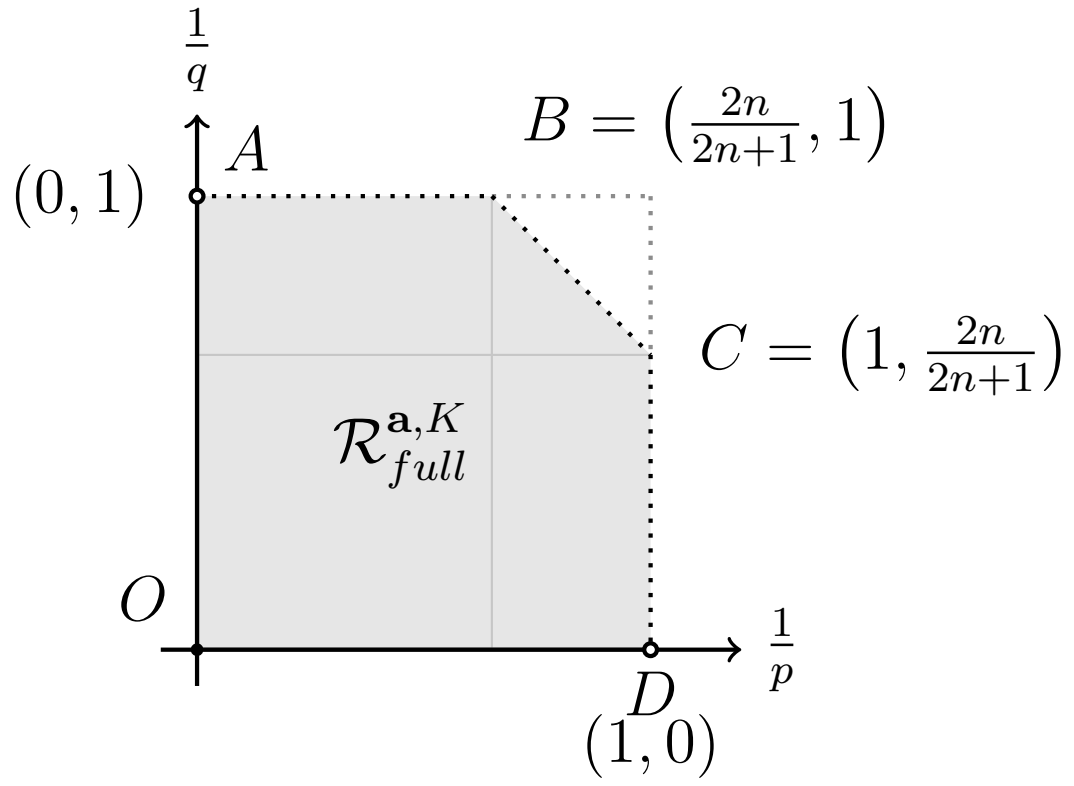
\begin{figure}[ht]
\centering
\begin{tikzpicture}[scale=3]

% -------------------------------------

\def\n{2}

% -------------------------------------------

\def\a{0.65}

%--------------------------------------------

\coordinate (O) at (0,0);
\coordinate (A) at (0,1);
\coordinate (B) at (\a,1);
\coordinate (C) at (1,\a);
\coordinate (D) at (1,0);

%-------------------------------------------
    
\fill[gray!20]
(O) -- (A) -- (B) -- (C) -- (D) -- cycle;

%--------------------------------------------
    
\draw[gray!45, thin]
(B) -- (\a,0);

\draw[gray!45, thin]
(0,\a) -- (C);

%--------------------------------------------

\draw[black!45, thick]
(O) -- (A);

\draw[black!45, thick]
(D) -- (O);

    %------------------------------------------------
    % Horizontal and vertical boundary
    %------------------------------------------------

    % A to B
\draw[black!45!black, dotted, thick]
(A) -- (B);

    % C to D
\draw[black!45!black, dotted, thick]
(C) -- (D);

    % B to (1,1)
\draw[black!45, dotted, thick]
(B) -- (1,1);

    % (1,1) to C
\draw[black!45, dotted, thick]
(1,1) -- (C);

    % B to C
\draw[black!45!black, thick, dotted]
(B) -- (C);

    %------------------------------------------------
    % Coordinate axes
    %------------------------------------------------
\draw[->, thick]
(-0.08,0) -- (1.20,0)
node[right] {$\frac{1}{p}$};

\draw[->, thick]
(0,-0.08) -- (0,1.18)
node[above] {$\frac{1}{q}$};

    %------------------------------------------------
    % Points
    %------------------------------------------------

\fill (O) circle (0.4pt);
\draw[fill=white, thick] (A) circle (0.4pt);
%\draw[fill=white, thick] (B) circle (0.4pt);
% \draw[fill=white, thick] (C) circle (0.4pt);
\draw[fill=white, thick] (D) circle (0.4pt);

    %------------------------------------------------
    % Labels
    %------------------------------------------------
\node[above left=2pt] at (O) {$O$};
\node[above right=1pt] at (A) {$A$};
\node[left=5pt] at (A) {$(0,1)$};
\node[above right=2pt] at (B)
{$B=\left(\frac{2n}{2n+1},1\right)$};
\node[right=5pt] at (C)
{$C=\left(1,\frac{2n}{2n+1}\right)$};
\node[below=0.02pt] at (D) {$D$};
\node[below=8pt] at (D) {$(1,0)$};

    %------------------------------------------------
    % Region label
    %------------------------------------------------
\node at (0.45,0.45) {$\mathcal{R}_{full}^{\mathbf a,K}$};
\end{tikzpicture}

\caption{Region $\mathcal{R}^{\mathbf a,K}_{full}$ for any $\mathbf{a}_1, \mathbf{a}_2 \leq 2n+2$.}
\label{figure-1-bilinear-degenerate-Koranyi-full}
\end{figure}

 \begin{figure}[ht]
\centering
\begin{minipage}[t]{0.48\textwidth}
\centering
\begin{tikzpicture}[scale=3]
 % Visual position of B and C
\def\a{0.65}
    %================================================
    % Coordinates
    %================================================
\coordinate (O) at (0,0);
\coordinate (A) at (0,1);
\coordinate (B) at (\a,1);
\coordinate (C) at (1,\a);
\coordinate (I) at (0,\a);
\coordinate (Z) at (0.85,0.80);
\coordinate (E) at (0.85,0);
\coordinate (J) at (\a,0);

    %================================================
    % Shaded region
    %================================================
    \fill[gray!20]
        (O) -- (A) -- (B) -- (Z) -- (E) -- (J) -- cycle;

    %================================================
    % Auxiliary dotted lines
    %================================================

    % Horizontal line through H and P
    \draw[black!45, dotted, thick]
        (A) -- (1,1);
    % (1,1) to C
    \draw[black!45, dotted, thick]
        (1,1) -- (C);

    % Vertical line through P
    \draw[black!45, dotted, thick]
        (B) -- (B |- O);

    % Horizontal line through Z and Q
    \draw[black!45, dotted, thick]
        (0,\a) -- (C);

    % Vertical line from Z to B
    \draw[black!45!black, dotted, thick]
        (Z) -- (E);

    % Vertical line through Q
    \draw[black!45, dotted, thick]
        (C) -- (C |- O);

    %================================================
    % Boundary of shaded region
    %================================================

    % Left boundary OH
    \draw[black!45, thick]
        (O) -- (A);

    % Bottom boundary OB
    \draw[black!45!black, thick]
        (O) -- (J) -- (E);

    % Top boundary HP
    \draw[black!45!black, dotted, thick]
        (A) -- (B);

    %================================================
    % Straight line P-Z-Q
    %================================================
    \draw[black!45, dotted, thick]
        (B) -- (C);

    %================================================
    % Vertical boundary Z-B
    %================================================
    \draw[black!45!black, dotted, thick]
        (B) -- (Z);

    %================================================
    % Coordinate axes
    %================================================

    \draw[->, thick]
        (-0.08,0) -- (1.20,0)
        node[right] {$\frac{1}{p}$};

    \draw[->, thick]
        (0,-0.08) -- (0,1.18)
        node[above] {$\frac{1}{q}$};

    %================================================
    % Points
       \fill (O) circle (0.4pt);

    \draw[fill=white, thick] (A) circle (0.4pt);
    %\draw[fill=white, thick] (B) circle (0.4pt);
    %\draw[fill=white, thick] (Z) circle (0.4pt);
    \draw[fill=white, thick] (E) circle (0.4pt);%================================================

    %================================================
    % Labels
    %================================================

    \node[above left=2pt] at (O) {$O$};

    \node[ left=2pt] at (I) {$I$};
    
    \node[above right=1pt] at (A) {$A$};
    %\node[left=5pt] at (A) {$(0,1)$};

    \node[above=7pt] at (B) {$B$};

    \node[right=5pt] at (Z) {$Z$};

    \node[below=2pt] at (E) {$E$};

    \node[below=2pt] at (J) {$J$};
     
    \node[right=5pt] at (C) {$C$};

    %================================================
    % Region label
    %================================================

    \node at (0.40,0.43)
        {$\mathcal{R}^{\mathbf a,K}_{full}$};

\end{tikzpicture}

\caption{Region $\mathcal{R}_{full}^{\mathbf a,K}$ for
$\mathbf a_1\leq 2n+2<\mathbf a_2$.}
\label{figure-2-bilinear-degenerate-Koranyi-full}
\end{minipage}
\hfill
\begin{minipage}[t]{0.48\textwidth}
\centering
%%%%%%%%%%%%%%%%%%%%%%%%%%%%%%%%%%%%%%%%%%%%%%%%%%%%%%%%%%%%%%%%%%%%%%%%%%%%%%%%%%%%%%%%%%%%%%%%%%%%%%%%%%%%%%%%%%%%%%%%%%%%%%%%%%
%\begin{figure}[ht]
%\centering

\begin{tikzpicture}[scale=3]
 % Visual position of B and C
\def\a{0.65}
    %================================================
    % Coordinates
    %================================================
    \coordinate (O) at (0,0);
     \coordinate (D) at (1,0);
    \coordinate (A) at (0,1);
    \coordinate (B) at (\a,1);
    \coordinate (C) at (1,\a);
    \coordinate (Z) at (0.85,0.80);
    \coordinate (Y) at (0.80,0.85);
    \coordinate (E) at (0.85,0);
    \coordinate (K) at (0,0.85);
 %   \coordinate (F) at (\a,0);
 %   \coordinate (H) at (0,\a);
    \coordinate (D) at (1,0);

    %================================================
    % Shaded region
    %================================================
    \fill[gray!20]
    (O) -- (K) -- (Y) -- (Z) -- (E) -- cycle;

    %================================================
    % Auxiliary dotted lines
    %================================================

    % Horizontal line through H and P
    \draw[black!45, dotted, thick]
        (A) -- (1,1);
    % (1,1) to D
    \draw[black!45, dotted, thick]
        (1,1) -- (D);
% Vertical line through P
%  \draw[black!45, dotted, thick]
%     (B) -- (B |- O);

% Horizontal line through Z and Q
% \draw[black!45, dotted, thick]
%    (0,\a) -- (C);

    % Vertical line from Z to E
    \draw[black!45, dotted, thick]
        (Z) -- (E);

% Vertical line through Q
%\draw[black!45, dotted, thick]
% (C) -- (C |- O);

    %================================================
    % Boundary of shaded region
    %================================================

    % Left boundary OH
    \draw[black!45, thick]
        (O) -- (K);
 % Vertical line from G to A
    \draw[black!45, dotted, thick]
        (K) -- (A);
    % Bottom boundary OB
    \draw[black!45, thick]
        (O) -- (E);

    % Top boundary HP
    \draw[black!45, dotted, thick]
        (A) -- (B);
  % B to C
    \draw[black!45, thick, dotted]
        (B) -- (C);

         % Y to Z
    \draw[black!45!black, thick, dotted]
        (Y) -- (Z);
    %================================================
    % Straight line P-Z-Q
    %================================================
    \draw[black!45!black, dotted, thick]
        (K) -- (Y);

    %================================================
    % Vertical boundary Z-B
    %================================================
    \draw[black!45!black, dotted, thick]
        (Z) -- (E);

    %================================================
    % Coordinate axes
    %================================================

    \draw[->, thick]
        (-0.08,0) -- (1.20,0)
        node[right] {$\frac{1}{p}$};

    \draw[->, thick]
        (0,-0.08) -- (0,1.18)
        node[above] {$\frac{1}{q}$};

    % Points
       \fill (O) circle (0.4pt);

    \draw[fill=white, thick] (K) circle (0.4pt);
    %\draw[fill=white, thick] (B) circle (0.4pt);
    %\draw[fill=white, thick] (Z) circle (0.4pt);
    \draw[fill=white, thick] (E) circle (0.4pt);%================================================
    %================================================

    \node[above left=2pt] at (O) {$O$};

    \node[above right=1pt] at (A) {$A$};
  %  \node[left=5pt] at (A) {$(0,1)$};

    \node[above=7pt] at (B) {$B$};

    \node[right=5pt] at (Z) {$Z$};

    \node[below=2pt] at (E) {$E$};

    \node[below=2pt] at (D) {$D$};

    \node[left=2pt] at (K) {$K$};

    \node[above=2pt] at (Y) {$Y$};
     
    \node[right=5pt] at (C) {$C$};

    %================================================
    % Region label
    %================================================

    \node at (0.40,0.43)
        {$\mathcal{R}_{full}^{\mathbf a,K}$};

\end{tikzpicture}

\caption{Region $\mathcal{R}_{full}^{\mathbf a,K}$ for
$\mathbf a_1, \mathbf{a}_2>2n+2$ and $\frac{1}{\mathbf{a}_1}+\frac{1}{\mathbf{a}_2}>\frac{1}{2n}.$}
\label{figure-3-bilinear-degenerate-Koranyi-full}
\end{minipage}
\end{figure}

%%%%%%%%%%%%%%%%%%%%
%%%%%%%%%%%%%%%%%%%%
%%%%%%%%%%%%%%%%%%%%

\medskip
\textbf{Case $1$:} If $\mathbf{a_1,a_2}\leq 2n+2.$ Observe that in this case, $\mathcal{R}_{full}^ {\mathbf{a}, K} = \mathcal{R}^K$. 

%%%%%%%%%%%%%%%%%%%%
%%%%%%%%%%%%%%%%%%%%
%%%%%%%%%%%%%%%%%%%%

\medskip
\textbf{Case $2$:} If $\min\{\mathbf{a}_1, \mathbf{a}_2\} \leq 2n+2 < \max \{ \mathbf{a}_1, \mathbf{a}_2 \}$. In this case, two symmetric possibilities arise. We consider only one of those, namely, $\mathbf{a}_1\leq 2n+2 <\mathbf{a}_2$ and denote 
\[
Z = \left( 1 - \frac{1}{2n+1} + \frac{2n+2}{(2n+1) \mathbf{a}_2}, \, 1-\frac{2n+2}{(2n+1) \mathbf{a}_2} \right), \quad E = \left( 1 - \frac{1}{2n+1} + \frac{2n+2}{(2n+1) \mathbf{a}_2}, \, 0 \right).
\]

%%%%%%%%%%%%%%%%%%%
%%%%%%%%%%%%%%%%%%%
%%%%%%%%%%%%%%%%%%%
%%%%%%%%%%%%%%%%%%%

\textbf{Case $3$:} 
If $2n+2<\min \{\mathbf{a}_1, \mathbf{a}_2\}$. In this case, depending on whether $\frac{1}{\mathbf{a}_1}+\frac{1}{\mathbf{a}_2}>\frac{1}{2n+2}$ or $\frac{1}{\mathbf{a}_1}+\frac{1}{\mathbf{a}_2}\leq\frac{1}{2n+2}$, we have two subcases. If $\frac{1}{\mathbf{a}_1}+\frac{1}{\mathbf{a}_2}>\frac{1}{2n+2}$, then the corresponding full region $\mathcal{R}_{full}^ {\mathbf{a},K}$ will intersect the line segment $BC$, and for the other case the region will not intersect the line segment  $BC$. We have provided the diagram only for the former case. Here, 
 \[
Y = \left( 1 - \frac{2n+2}{(2n+1) \mathbf{a}_1}, \, 1 - \frac{1}{2n+1} + \frac{2n+2}{(2n+1) \mathbf{a}_1} \right), \quad K = \left( 0, \, 1 - \frac{1}{2n+1} + \frac{2n+2}{(2n+1) \mathbf{a}_1} \right).
\]

%%%%%%%%%%%%%%%%
%%%%%%%%%%%%%%%%
%%%%%%%%%%%%%%%%
%%%%%%%%%%%%%%%%

With that, here is our main result on the boundedness of $\mathcal{M}^{\mathbf{a},K}_{{full}}$. The boundedness region for the estimate \eqref{estimate_for_the_Product_Koranyi_full_degenrate} is sharp except possibly at the endpoint and the same shall be discussed in Section \ref{sec:sharpness-examples}. 

\begin{theorem} 
\label{Main-theorem-Koranyi-type-for-degenrate}
Let $n\geq 2$. Assume that $\mathbf{a} = (\mathbf{a}_1, \mathbf{a}_2) \in [2,\infty)^2$. Let $1 \leq p, q \leq \infty$ and $r>0$ be such that $\frac{1}{r}=\frac{1}{p}+\frac{1}{q}$. Then, we have 
\begin{align} 
\label{estimate_for_the_Product_Koranyi_full_degenrate} 
\|\mathcal{M}^{\mathbf{a},K}_{full}(F,G)\|_{L^r(\mathbb{H}^n)}\lesssim\|F\|_{L^p(\mathbb{H}^n)}\|G\|_{L^q(\mathbb{H}^n)}, \quad \text{ for all } (1/p,1/q)\in \mathcal{R}_{full}^{\mathbf a,K}. 
\end{align} 
Moreover, in Case $1$, that is, when $a_1, a_2 \leq 2n+2$, estimate \eqref{estimate_for_the_Product_Koranyi_full_degenrate} holds true on the boundary $(A,B)$ and $(C,D)$ as well as the following weak-type estimates holds true at $A$ and $D$:
\begin{align}
\|\mathcal{M}^{\mathbf{a},K}_{full} (F,G)\|_{L^{1,\infty}(\mathbb{H}^n)} 
& \lesssim 
\|F\|_{L^1(\mathbb{H}^n)}\,\|G\|_{L^\infty(\mathbb{H}^n)}, 
\label{thm:main-Koranyi-degen-weak1-est} \\
\text{and} \quad \|\mathcal{M}^{\mathbf{a},K}_{full} (F,G)\|_{L^{1,\infty}(\mathbb{H}^n)} 
& \lesssim 
\|F\|_{L^\infty(\mathbb{H}^n)}\,\|G\|_{L^1(\mathbb{H}^n)}. 
\label{thm:main-Koranyi-degen-weak2-est}
\end{align}
\end{theorem}

Since $\mathcal{M}^K_{full}(F,G) \lesssim \mathcal{M}^{\mathbf{a},K}_{full} (|F|, |G|)$ for $\mathbf{a} = (4,4)$, we get the boundedness of $\mathcal{M}^K_{full}$ as an immediate corollary of Theorem \ref{Main-theorem-Koranyi-type-for-degenrate}. 
\begin{corollary} \label{Main-theorem-corollary}
Let $n \geq 2$. Let $1 \leq p, q \leq \infty$ and $0 <r \leq \infty$. Then, the  operator $\mathcal{M}^K_{full}$ satisfies the estimate \eqref{estimate_for_the_Product_Koranyi_full_degenrate}, if $r> \frac{2n+1}{4n+1}$, where $\frac{1}{p}+\frac{1}{q}=\frac{1}{r}$, except the case $(p, q, r) =(1, \infty, 1)=D \text{ or } (\infty, 1, 1)=A$, where we have the weak-type estimates as in \eqref{thm:main-Koranyi-degen-weak1-est} and \eqref{thm:main-Koranyi-degen-weak2-est}.
\end{corollary}

Next, we discuss the boundedness of $\mathcal{M}^{\mathbf{a},K}_{lac}$. As in \cite{Borges-Tainara-Foster-bounds-for-lacunary-bilinear-spherical-maximal-functions}, in our case also we shall show that if $\mathcal{R}^{\mathbf{a}, K}_{1}$ denotes the H\"older boundedness region \eqref{def:koranyi-holder-bdd-region-single-scale} for the single scale averages $\mathcal{A}^{\mathbf{a},K}_1$, then its interior is contained in the H\"older boundedness region of $\mathcal{M}^{\mathbf{a},K}_{lac}$ whenever $\mathbf{a} \in [2,2n+2]^2$. Here, 
\begin{align}
\label{def:koranyi-holder-bdd-region-single-scale} 
\mathcal{R}^{\mathbf{a}, K}_{1} = \left\{ \left( \frac{1}{p}, \frac{1}{q} \right) \in [0,1]^2 : \|\mathcal{A}^{\mathbf{a},K}_1\|_{L^p \times L^q \rightarrow L^r}<\infty  \text{ with } \frac{1}{p} + \frac{1}{q} = \frac{1}{r} \right\}. 
\end{align} 

%%%%%%%%%%%%%%%
%%%%%%%%%%%%%%%
%%%%%%%%%%%%%%%
%%%%%%%%%%%%%%%

\begin{theorem} 
\label{Main-theorem-lacunary-Koranyi-type-for-degenrate}
Let $\mathbf{a}=(\mathbf{a}_1,\mathbf{a}_2)\in [2,2n+2]^2$. Then, for any $(1/p,1/q)\in int(\mathcal{R}^\mathbf{a, K}_{1})$ such that $\frac{1}{p}+\frac{1}{q} = \frac{1}{r} $, we have
\begin{equation} \label{estimates_for_lacunary_Product_Koranyi_degenerate} \| \mathcal{M}^{\mathbf{a},K}_{lac} (F,G)\|_{L^r(\mathbb{H}^n)} \lesssim \|F\|_{L^p(\mathbb{H}^n)} \|G\|_{L^q(\mathbb{H}^n)}. 
\end{equation}
\end{theorem}

\begin{corollary} \label{Main_Theorem_Lacunary} 
Let $n\geq2$ and $\mathcal{Q} = \mathcal{R}^{\mathbf{a}, K}_{1}$ with $\mathbf{a} = (4,4)$. For any $(1/p,1/q) \in int(\mathcal{Q})$ such that $\frac{1}{p} + \frac{1}{q} = \frac{1}{r}$, the operator $\mathcal{M}^K_{lac}$ boundedly maps $L^p{(\mathbb{H}^n})\times L^q{(\mathbb{H}^n})$ to $L^r(\mathbb{H}^n)$.
\end{corollary}
\medskip

Next, in a similar way to the above, we extend the Nevo--Thangavelu case studied in \cite{ghosh2026bilinearsphericalmaximalfunction} to the following degenerate surfaces. Given $\mathbf{a}=(\mathbf{a_1,a_2})\in [1,\infty)^2$, consider 
\begin{align} 
\label{expression_of_the_non_dengenerate_surfaces_of_Nevo_Thangavelu_type}
S^{\mathbf{a}}_{NT} := \{(g_1, g_2) = (z_1,0,z_2,0) \in \mathbb{H}^n \times \mathbb{H}^n : |z_1|^{\mathbf{a_1}} +|z_2|^{\mathbf{a_2}}=1 \}, 
\end{align}
With the naturally induced surface measure $d\mu_{NT}^\mathbf{a}$ on $S^{\mathbf{a}}_{NT} $ (see Section \ref{sec:NT-maximal-degenerate} for details), consider the averaging operators $\mathcal{A}_\mathbf{s}^{\mathbf{a},NT}$ defined by 
\begin{align}\label{express-averaging-operator-nevotype}
\mathcal{A}_\mathbf{s}^{\mathbf{a},NT}(F,G)(g) = \int\limits_{S^{\mathbf{a}}_{NT}} F(g \cdot\tilde{\delta}_{s_1} g_1^{-1}) \, G(g \cdot\tilde{\delta}_{s_2}g_2^{-1}) \, d\mu_{NT}^\mathbf{a}(g_1, g_2).
\end{align}
The corresponding full maximal operator is then defined by 
\begin{align}\label{expression_of_full_maximal_operator_at_degenerate_surface_Nevo_Thangavelu}
\mathcal{M}^{\mathbf{a},NT}_{full}(F,G)(g) =\sup_{\mathbf{s}\in (0,\infty)^2} |\mathcal{A}_\mathbf{s}^{\mathbf{a},NT}(F,G)(g)|, 
\end{align}
and the lacunary one by 
\begin{align}\label{expression_of_lacunary_maximal_operator_at_degenerate_surface_Nevo_Thangavelu}
 \mathcal{M}^{\mathbf{a},NT}_{lac}(F,G)(g)=\sup_{l\in\mathbb{Z}} |\mathcal{A}_{2^l}^{\mathbf{a},NT}(F,G)(g)|.
\end{align} 

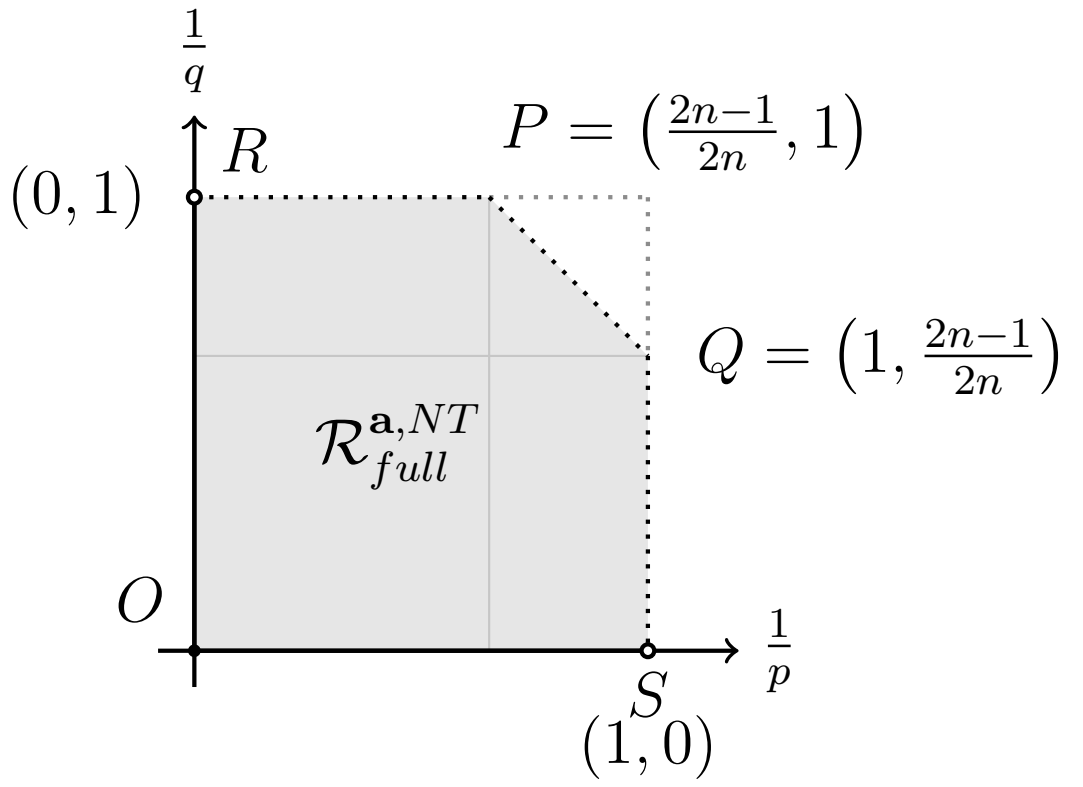
\begin{figure}[ht]
\centering
\begin{minipage}{0.25\textwidth}
\raggedright
\textbf{Case 1:} If $\mathbf{a}_1,\mathbf{a}_2\leq 2n$.
\end{minipage}
\hfill
\begin{minipage}{0.70\textwidth}
\centering

\begin{tikzpicture}[scale=
3]

    %================================================
    % Parameter
    %================================================
    \def\n{2}
     \def\a{0.65}
    %================================================
    % Coordinates
    %================================================
    \coordinate (O) at (0,0);
    \coordinate (R) at (0,1);
    \coordinate (P) at (\a,1);
    \coordinate (Q) at (1,\a);
    \coordinate (S) at (1,0);

    %================================================
    % Shaded region P
    %================================================
    \fill[gray!20]
        (O) -- (R) -- (P) -- (Q) -- (S) -- cycle;

    %================================================
    % Auxiliary horizontal and vertical lines
    %================================================
    \draw[gray!45, thin]
        (P) -- (\a,0);

    \draw[gray!45, thin]
        (0,\a) -- (Q);

    %================================================
    % Boundary of the pentagon
    %================================================
    \draw[black!45!black, thick]
        (O) -- (R);

    \draw[black!45!black, thick]
        (S) -- (O);

    %================================================
    % Dotted triangle P-(1,1)-Q
    %================================================
    \draw[black!45, dotted, thick]
        (P) -- (1,1);

         \draw[black!45!black, dotted, thick]
        (R) -- (P);

         \draw[black!45!black, dotted, thick]
        (Q) -- (S);

    \draw[black!45, dotted, thick]
        (1,1) -- (Q);

    % Dotted slanted edge PQ
    \draw[black!45!black, thick, dotted]
        (P) -- (Q);

    %================================================
    % Coordinate axes
    %================================================
    \draw[->, thick]
        (-0.08,0) -- (1.20,0)
        node[right] {$\frac{1}{p}$};

    \draw[->, thick]
        (0,-0.08) -- (0,1.18)
        node[above] {$\frac{1}{q}$};

    %================================================
    % Points
    %================================================

    % O
    \fill (O) circle (0.4pt);

    % Open circles
    \draw[fill=white, thick] (R) circle (0.4pt);
  % \draw[fill=white, thick] (P) circle (0.4pt);
   % \draw[fill=white, thick] (Q) circle (0.4pt);
    \draw[fill=white, thick] (S) circle (0.4pt);

    %================================================
    % Labels
    %================================================

    % O
    \node[above left=2pt] at (O) {$O$};

    % H
    \node[above right=1pt] at (R) {$R$};
    \node[left=5pt] at (R) {$(0,1)$};

    % P
    \node[above right=1pt and -2pt] at (P)
    {$P=\left(\frac{2n-1}{2n},1\right)$};

    % Q
    \node[right=5pt] at (Q)
        {$Q=\left(1,\frac{2n-1}{2n}\right)$};

    % G
    \node[below=0.02pt] at (S) {$S$};
    \node[below=8pt] at (S) {$(1,0)$};

    %================================================
    % Region label
    %================================================

    \node at (0.45,0.45) {$\mathcal{R}^{\mathbf{a},NT}_{full}$};

\end{tikzpicture}

\caption{Region $\mathcal{R}^{\mathbf{a},NT}_{full}$, when
$\mathbf{a}_1,\mathbf{a}_2\leq 2n$.}
\label{figure_NT-case-1}
\end{minipage}
\end{figure}

%%%%%%%%%%%%%%%%%%%%
%%%%%%%%%%%%%%%%%%%%
%%%%%%%%%%%%%%%%%%%%

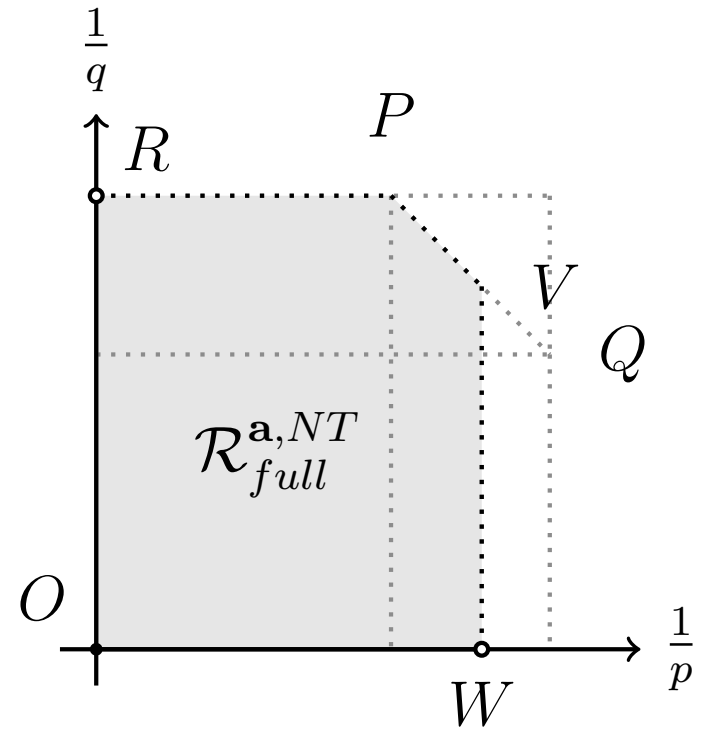
\begin{figure}[ht]
\centering
\raisebox{5cm}{%
\begin{minipage}[t]{0.38\textwidth}
\raggedright
\textbf{Case 2:} If
$\min\{\mathbf{a}_1,\mathbf{a}_2\}\leq 2n
<\max\{\mathbf{a}_1,\mathbf{a}_2\}$.
Suffices to display the case
$\mathbf{a}_1\leq 2n<\mathbf{a}_2$.
\[
V=
\left(
1-\frac{1}{2n}+\frac{1}{\mathbf{a}_2},
1-\frac{1}{\mathbf{a}_2}
\right)
\]
and
\[
W=
\left(
1-\frac{1}{2n}+\frac{1}{\mathbf{a}_2},
0
\right).
\]
\end{minipage}%
}
\hfill
\begin{minipage}[t]{0.58\textwidth}
\centering

\begin{tikzpicture}[scale=3]
 % Visual position of B and C
\def\a{0.65}
    %================================================
    % Coordinates
    %================================================
    \coordinate (O) at (0,0);
    \coordinate (R) at (0,1);
    \coordinate (P) at (\a,1);
    \coordinate (Q) at (1,\a);
    \coordinate (V) at (0.85,0.80);
    \coordinate (W) at (0.85,0);
    \coordinate (F) at (\a,0);

    %================================================
    % Shaded region
    %================================================
    \fill[gray!20]
        (O) -- (R) -- (P) -- (V) -- (W) -- (F) -- cycle;

    %================================================
    % Auxiliary dotted lines
    %================================================

    % Horizontal line through H and P
    \draw[black!45, dotted, thick]
        (R) -- (1,1);
    % (1,1) to C
    \draw[black!45, dotted, thick]
        (1,1) -- (Q);

    % Vertical line through P
    \draw[black!45, dotted, thick]
        (P) -- (P |- O);

    % Horizontal line through Z and Q
    \draw[black!45, dotted, thick]
        (0,\a) -- (Q);

    % Vertical line from Z to B
    \draw[black!45, dotted, thick]
        (V) -- (W);

    % Vertical line through Q
    \draw[black!45, dotted, thick]
        (Q) -- (Q |- O);

    %================================================
    % Boundary of shaded region
    %================================================

    % Left boundary OH
    \draw[black!45, thick]
        (O) -- (R);

    % Bottom boundary OB
    \draw[black!45!black, thick]
        (O) -- (F) -- (W);

    % Top boundary HP
    \draw[black!45!black, dotted, thick]
        (R) -- (P);

    %================================================
    % Straight line P-Z-Q
    %================================================
    \draw[black!45, dotted, thick]
        (P) -- (Q);

     \draw[black!45!black, dotted, thick]
        (P) --(V);%================================================
    % Vertical boundary Z-B
    %================================================
    \draw[black!45!black, dotted, thick]
        (V) -- (W);

    %================================================
    % Coordinate axes
    %================================================

    \draw[->, thick]
        (-0.08,0) -- (1.20,0)
        node[right] {$\frac{1}{p}$};

    \draw[->, thick]
        (0,-0.08) -- (0,1.18)
        node[above] {$\frac{1}{q}$};

    %================================================
    % Points
     % O
    \fill (O) circle (0.4pt);

    % Open circles
    \draw[fill=white, thick] (R) circle (0.4pt);
  % \draw[fill=white, thick] (P) circle (0.4pt);
   % \draw[fill=white, thick] (Q) circle (0.4pt);
    \draw[fill=white, thick] (W) circle (0.4pt);

    %================================================

    \node[above left=2pt] at (O) {$O$};

    \node[above right=1pt] at (R) {$R$};
    %\node[left=5pt] at (A) {$(0,1)$};

    \node[above=7pt] at (P) {$P$};

    \node[right=5pt] at (V) {$V$};

    \node[below=2pt] at (W) {$W$};

   % \node[below=2pt] at (F) {$F$};
     
    \node[right=5pt] at (Q) {$Q$};

    %================================================
    % Region label
    %================================================

    \node at (0.40,0.43)
        {$\mathcal{R}^{\mathbf{a},NT}_{full}$};

\end{tikzpicture}

\caption{Region $\mathcal{R}^{\mathbf{a},NT}_{full}$ for
$\mathbf a_1\leq 2n<\mathbf a_2$.}
\label{figure:NT-case-2}
\end{minipage}
\end{figure}

%%%%%%%%%%%%%%%%%%%%
%%%%%%%%%%%%%%%%%%%%
%%%%%%%%%%%%%%%%%%%%

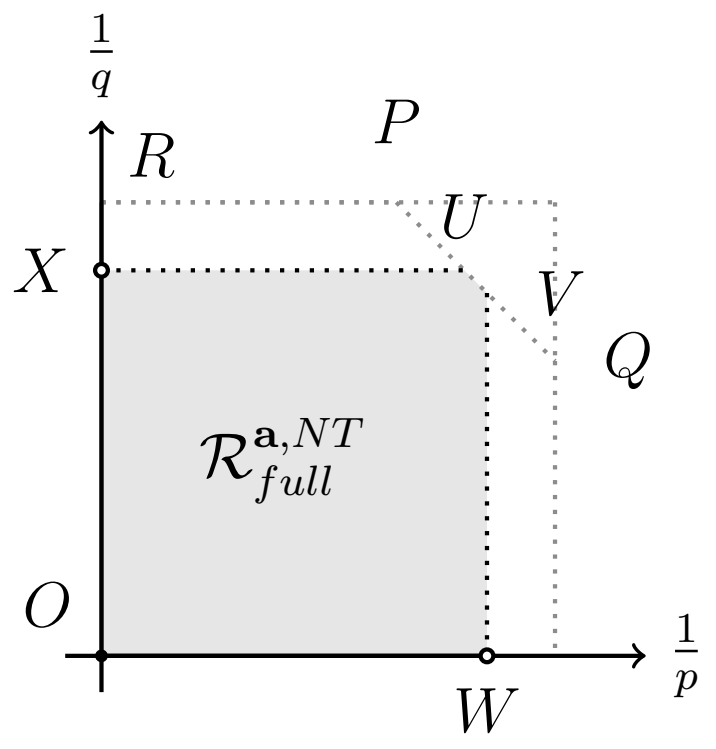
\begin{figure}[ht]
\centering
\raisebox{5cm}{%
\begin{minipage}[t]{0.38\textwidth}
\raggedright
\textbf{Case 3:} If $2n<\min{\{\mathbf{a}_1,\mathbf{a}_2\}}.$ Here also we consider the situation whenever $\frac{1}{\mathbf{a}_1}+\frac{1}{\mathbf{a}_2}>\frac{1}{2n}$. Here
 \[
U=\left(1-\frac{1}{\mathbf{a}_1},1-\frac{1}{2n}+\frac{1}{\mathbf{a}_1}\right) 
\]
and 
\[
X= \left(0,1-\frac{1}{2n}+\frac{1}{\mathbf{a}_1}\right).
\]
\end{minipage}
}
\hfill
\begin{minipage}[t]{0.58\textwidth}
\centering

\begin{tikzpicture}[scale=3]
 % Visual position of B and C
\def\a{0.65}
    %================================================
    % Coordinates
    %================================================
    \coordinate (O) at (0,0);
    \coordinate (R) at (0,1);
    \coordinate (P) at (\a,1);
    \coordinate (Q) at (1,\a);
    \coordinate (V) at (0.85,0.80);
    \coordinate (U) at (0.80,0.85);
    \coordinate (W) at (0.85,0);
    \coordinate (F) at (\a,0);
    \coordinate (X) at (0,0.85);

    %================================================
    % Shaded region
    %================================================
    \fill[gray!20]
        (O) -- (X) -- (U) -- (V) -- (W) -- cycle;

    %================================================
    % Auxiliary dotted lines
    %================================================

    % Horizontal line through H and P
    \draw[black!45, dotted, thick]
        (R) -- (1,1);
    % (1,1) to C
    \draw[black!45, dotted, thick]
        (1,1) -- (Q);

    % Horizontal line through Z and Q

    % Vertical line from Z to B
   % \draw[black!45!black, dotted, thick]
        %(U) -- (V);

    % Vertical line through Q
    \draw[black!45, dotted, thick]
        (Q) -- (Q |- O);

    %================================================
    % Boundary of shaded region
    %================================================

    % Left boundary OH
    \draw[black!45, thick]
        (O) -- (R);

    % Bottom boundary OB
    \draw[black!45!black, thick]
        (O) -- (F) -- (W);

    % Top boundary HP
    \draw[black!45, dotted, thick]
        (R) -- (P);

    %================================================
    % Straight line P-Z-Q
    %================================================
    \draw[black!45, dotted, thick]
        (P) -- (Q);

    %================================================
    % Vertical boundary Z-B
    %================================================
    \draw[black!45!black, dotted, thick]
        (V) -- (W);

       \draw[black!45!black, dotted, thick]
        (X) -- (U);
    %================================================
    % Coordinate axes
    %================================================

    \draw[->, thick]
        (-0.08,0) -- (1.20,0)
        node[right] {$\frac{1}{p}$};

    \draw[->, thick]
        (0,-0.08) -- (0,1.18)
        node[above] {$\frac{1}{q}$};

    %================================================
    % Points
    %================================================
    % O
    \fill (O) circle (0.4pt);

    % Open circles
    \draw[fill=white, thick] (X) circle (0.4pt);
  % \draw[fill=white, thick] (P) circle (0.4pt);
   % \draw[fill=white, thick] (Q) circle (0.4pt);
    \draw[fill=white, thick] (W) circle (0.4pt);

    %================================================
    % Labels
    %================================================

    \node[above left=2pt] at (O) {$O$};

    \node[above right=1pt] at (R) {$R$};
    %\node[left=5pt] at (A) {$(0,1)$};

    \node[above=7pt] at (P) {$P$};

    \node[right=5pt] at (V) {$V$};

    \node[above=2pt] at (U) {$U$};

     \node[left=2pt] at (X) {$X$};

    \node[below=2pt] at (W) {$W$};

   % \node[below=2pt] at (F) {$F$};
     
    \node[right=5pt] at (Q) {$Q$};

    %================================================
    % Region label
    %================================================

    \node at (0.40,0.43)
        {$\mathcal{R}^{\mathbf{a},NT}_{full}$};

\end{tikzpicture}

\caption{Region $\mathcal{R}^{\mathbf{a},NT}_{full}$ for
$\mathbf a_1,\mathbf a_2> 2n$ and $\frac{1}{\mathbf{a}_1}+\frac{1}{\mathbf{a}_2}>\frac{1}{2n}$. }
\label{figure:NT-case-3}
\end{minipage}
\end{figure}

In the above diagrams, for $\mathbf{a}=(\mathbf{a_1,a_2)}\in [1,\infty)^2$, we have set
\begin{equation}\label{Domain_for_the_Nevo_thangavelu_type_non_degenerate}
\mathcal{R}^{\mathbf{a},NT}_{full} := \left\{ \left( \frac{1}{p}, \frac{1}{q} \right) \in \mathcal{R}^{NT} : \frac{1}{p} < 1 - \frac{ (\mathbf{a}_2-2n)_+}{2n \, \mathbf{a_2}}, \, \frac{1}{q} < 1 - \frac{ (\mathbf{a}_1-2n)_+}{2n \, \mathbf{a_1}} \right\}, 
\end{equation}
where $\mathcal{R}^{NT}$ is given in \eqref{range_for_the_Bilinear_NT}. Also, observe that $\mathcal{R}^{\mathbf{a},NT}_{full} = \mathcal{R}^{NT}$, whenever $\mathbf{a_1,a_2} \leq 2n.$  

\medskip
Similar to the product Kor\'anyi case, here also we have three cases depending on whether $\mathbf{a}_1$, $\mathbf{a}_2$ are less than or greater than the dimension $2n$. 

%%%%%%%%%%%%%%%%%%%%

\begin{theorem}
\label{thm:full-max-NT-degenerate}
Let $n\geq 2$. For any $\mathbf{a}\in [1,\infty)^2$, and $1 < p,q\leq \infty$, if $\frac{1}{r}=\frac{1}{p}+\frac{1}{q}$, then 
\begin{equation} 
\label{estimate_for_the_Nevo_Thangavelu_full_degenrate}
\|\mathcal{M}^{\mathbf{a},NT}_{full}(F,G)\|_{L^r(\mathbb{H}^n)}\lesssim\|F\|_{L^p(\mathbb{H}^n)}\|G\|_{L^q(\mathbb{H}^n)}, \quad \text{ for all } (1/p,1/q)\in \mathcal{R}^{\mathbf{a},NT}_{full}.
\end{equation} 
\end{theorem}

One may note here that unlike Theorem \ref{Main-theorem-Koranyi-type-for-degenrate}, we haven't mentioned anything about the boundary. In particular, the boundedness situation at the extreme points $R$ and $S$ is not addressed in Theorem \ref{thm:full-max-NT-degenerate}. Actually,  in the proof of Theorem \ref{thm:full-max-NT-degenerate}, we follow the general strategy developed by Jeong--Lee \cite{Maximal-estimates-for-the-bilinear-spherical-averages} for the full bilinear maximal operator. In particular, by exploiting a suitable slicing argument, we reduce the desired estimate to an application of Hopf's ergodic maximal theorem $\Lambda$ (defined in \eqref{def:Hopfs-ergodic-maximal}) in our proof, and we do not know whether $\Lambda$ is of weak-type $(1,1)$ or not.

\medskip 
Finally, for the lacunary operator $\mathcal{M}^{\mathbf{a},NT}_{lac}$, as in the Kor\'anyi case, we shall show have that if $\mathcal{R}_1^{\mathbf{a},NT}$ denotes the H\"older boundedness region \eqref{def:NT-holder-bdd-region-single-scale} for the single scale averages $\mathcal{A}^{\mathbf{a},NT}_1$, then its interior is contained in the H\"older boundedness region of $\mathcal{M}^{\mathbf{a},NT}_{lac}$ whenever $\mathbf{a} \in [1,2n]^2$. Here, 
\begin{align}
\label{def:NT-holder-bdd-region-single-scale} 
\mathcal{R}^\mathbf{a,NT}_{1} = \left\{ \left( \frac{1}{p}, \frac{1}{q} \right) \in [0,1]^2 : \|\mathcal{A}^{\mathbf{a},NT}_1\|_{L^p \times L^q \rightarrow L^r}<\infty  \text{ with } \frac{1}{p} + \frac{1}{q} = \frac{1}{r} \right\}. 
\end{align} 

\begin{theorem}
\label{thm:lac-max-NT-degenerate}
Let $\mathbf{a}=(\mathbf{a}_1,\mathbf{a}_2)\in [1,2n]^2$. For any $(1/p,1/q) \in int(\mathcal{R}_1^{\mathbf{a},NT})$, if $\frac{1}{r}=\frac{1}{p}+\frac{1}{q}$, then 
\begin{equation} 
\label{estimates_for_lacunary_nevo_thangavelu_degenerate}    \|\mathcal{M}^{\mathbf{a},NT}_{lac}(F,G)\|_{L^r(\mathbb{H}^n)}\lesssim\|F\|_{L^p(\mathbb{H}^n)}\|G\|_{L^q(\mathbb{H}^n)}.
\end{equation}
\end{theorem}

Regarding, the H\"older boundedness region \eqref{def:NT-holder-bdd-region-single-scale} for the single scale averages $\mathcal{A}^{\mathbf{a},NT}_1$, we have the following result. 
\begin{theorem}
\label{theorem-improving-single-scales}
Let $n\geq 2$. Assume that $\mathbf{a} = (\mathbf{a}_1, \mathbf{a}_2) \in [1,\infty)^2$. Let $1 \leq p, q \leq \infty$ and $r>0$ be such that $\frac{1}{r}=\frac{1}{p}+\frac{1}{q}$. Let $\mathcal{R}^K$ denote the open pentagon with vertices $O, A, B, C$, and $D$  including line segments $[O,A)$ and $[O,D)$ as in Figure \ref{figure-1-bilinear-degenerate-Koranyi-full}. Then, 
\begin{align} 
\label{estimate:improving-single-scales} 
\| \mathcal{A}^{\mathbf{a},NT}_1 (F,G) \|_{L^r(\mathbb{H}^n)} \lesssim \|F\|_{L^p(\mathbb{H}^n)} \|G\|_{L^q(\mathbb{H}^n)}, 
\end{align} 
for all $(1/p,1/q)$ coming from $\mathcal{R}^K$ or line
segments $[O, A], [A, B), (C, D], [O, D]$. 
\end{theorem}

Proofs of Theorems \ref{thm:full-max-NT-degenerate}, \ref{thm:lac-max-NT-degenerate} and \ref{theorem-improving-single-scales} draw upon ideas and techniques developed in several works. For Theorem \ref{thm:full-max-NT-degenerate}, we follow the approach of \cite{Maximal-estimates-for-the-bilinear-spherical-averages}, particularly the use of a slicing argument. The proof of Theorem \ref{thm:lac-max-NT-degenerate} is inspired by ideas from  \cite{Borges-Tainara-Foster-bounds-for-lacunary-bilinear-spherical-maximal-functions, Aswin-Hickman, ghosh2026bilinearsphericalmaximalfunction}. Finally, the proof of Theorem \ref{theorem-improving-single-scales} draws on the techniques developed in \cite[Section 4]{Iosevich} as was also successfully used in \cite{ghosh2026bilinearsphericalmaximalfunction}.
However, working with degenerate surfaces requires additional care, particularly in setting up the coarea formula and implementing the slicing arguments appropriately.

\medskip 
In the next section, we begin with a detailed discussion on the surface measures on the Kor\'anyi sphere, the product Kor\'anyi sphere as well as the degenerate product Kor\'anyi surface, before delving into the proofs on all the main results in the subsequent sections. 

\section{Surface measures} 
\label{sec:surface-measures} 

In the first two subsections, we shall describe the Radon measure $d\sigma^{2n}_K$ appearing in the polar decomposition \eqref{polar-decompostion-withrespectto-the-Koranyi-norm} of the Haar measure on $\mathbb{H}^n$ and the Radon measure $d\mu^{4n+1}_K$ supported on $\Sigma_{K,1}$ which provides us with the polar decomposition \eqref{expres-polar-decoposition-producttt} of the Haar measure on $\mathbb{H}^n\times\mathbb{H}^n$. This helps us in defining the Radon measure $d\mu^\mathbf{a}_K$ in a natural way on the degenerate surfaces $\mathcal{S}_K^{\mathbf{a}}$ and we shall pursue the same in Subsection 
\ref{exp-surface-measure-on-the-degenerate-product-Koranyi-sphere}. Here onwards, we shall use the notation $\nabla$ to always denote the Euclidean gradient.

%%%%%%%%%%%%%%%%%%%%
%%%%%%%%%%%%%%%%%%%%
%%%%%%%%%%%%%%%%%%%%
%%%%%%%%%%%%%%%%%%%%

\subsection{Radon measure \texorpdfstring{$d\sigma^{2n}_K$}{} on the Kor\'anyi sphere \texorpdfstring{$S_{K,1}$}{}} 
\label{exp-surface-measure-on-the-linear-Koranyi-sphere}

%In this subsection, we shall explain the Radon measure $d\sigma^{2n}_K$ appearing in the polar decomposition \eqref{polar-decompostion-withrespectto-the-Koranyi-norm}. 

Take the smooth function $\phi(z,t)=|z|^4+t^2.$ Clearly, $\nabla\phi(z,t) \neq 0$ on any $r-$radius Kor\'anyi sphere $S_{K,r}$. Suppose $d\mathcal{H}^{2n}_{K,r}$ denotes the $2n$-dimensional surface measure on $S_{K,r}$. An explicit expression for $d\mathcal{H}^{2n}_{K,r}$ will be given later on (see \eqref{Expression-for-the-measure-on-Koranyi-sphere}). Since the Kor\'anyi sphere $S_{K,r}$ is a $2n$-dimensional hypersurface in $\mathbb{H}^n$, this measure is naturally the corresponding Hausdorff measure. However, $d\mathcal{H}^{2n}_{K,r}$ is not the appropriate measure for obtaining the polar decomposition of the Haar measure on $\mathbb{H}^n$. Therefore, we explicitly construct a Radon measure $d\sigma^{2n}_K$ supported on $S_{K,1}$ that yields the desired polar decomposition. 

\medskip
Applying the coarea formula and an appropriate change of variables, we obtain the following, 
\begin{align}\label{expresstion-of-the-surface-measure-we-need-for-the-polar-decomposition}
\int_{\mathbb{H}^n} f(z,t)\,dz\,dt&=\int_0^\infty\left(\int_{\phi(z,t)=s}f(z,t)\,\frac{d\mathcal{H}_{K,s^{\,\frac{1}{4}}}^{2n}(z,t)}{|\nabla(\phi(z,t)|}\right) \,ds \\ 
&=\int_0^\infty\left(\int_{S_{K,r}}f(z,t)\,d\sigma^{2n}_{K,r}(z,t)\right)\,dr\notag\\
&=\int_0^\infty\left(\int_{S_{K,1}}f(\tilde{\delta_r}w)\,d\sigma^{2n}_K(w)\right)r^{\mathcal{Q}-1} \,dr \notag,
\end{align}
which is the polar decomposition mentioned in the introduction as \eqref{polar-decompostion-withrespectto-the-Koranyi-norm} with the Radon measure $d\sigma^{2n}_K = d\sigma_{K,1}^{2n}$ supported on $S_{K,1}$. Note also that $d\sigma_{K,r}^{2n} := |\nabla \phi|^{-1} 4r^3\,d\mathcal{H}_{K,r}^{2n}$ and is supported on $S_{K,r}$. This measure is useful in developing a slicing-type argument for the product Heisenberg group (see, for example, Subsection \ref{finding-slicing-argument-form-for-product-Koranyi}), and thus we are naturally led to analyse it in more details. Let us begin with the Hausdorff measure (surface measure), $d\mathcal{H}^{2n}_K=d\mathcal{H}^{2n}_{K,1}$ on $S_{K,1}$. Let $N(z,t) = (|z|^4 + t^2)^{1/4}-1.$ Then, $\nabla N(z,t)
= \left( \frac{|z|^2 z}{(|z|^4+t^2)^{3/4}}, \; \frac{t}{2(|z|^4+t^2)^{3/4}} \right).$
In particular, $|\nabla N(z,t)| =
\sqrt{|z|^6 + \frac{t^2}{4}} = \frac{1}{4}\sqrt{16|z|^6 + 4t^2}$ on $S_{K,1}$.

%%%%%%%%%%%%%%%%%%%%
%%%%%%%%%%%%%%%%%%%%
%%%%%%%%%%%%%%%%%%%%
%%%%%%%%%%%%%%%%%%%%

\subsubsection{Surface measure $d\mathcal{H}^{2n}_K$ via graph parametrization} \label{subsection-expression-for-the-surface-measure-on-S_K,1}
Let us first decompose the surface $S_{K,1}$ as follows: 
\[
S_{K,1}
=
\{(z,t): t = \pm \sqrt{1 - |z|^4},\ |z|<1\}
\cup
\{(z,0): |z|=1\}.
\]
Note that the set $\{(z,0): |z|=1\}$ has $d\mathcal{H}_K^{2n}$-measure zero, and thus it suffices to only work with the first set in the above decomposition. 

\medskip 
Define the maps $p_\pm$ on the unit disc on $\mathbb{C}^n$ by 
\[
p_\pm(z) = \pm \sqrt{1 - |z|^4}, 
\]
%and
%\[\Phi_\pm(z) = (z, p_\pm(z)),\]
so that $S_{K,1}$ is the union of the graph of these functions $p_\pm$. For each part, we can compute the surface measure via the formula
\[
d\mathcal{H}^{2n}_K
=
\sqrt{1 + |\nabla p_\pm(z)|^2}\, dz.
\]
Since 
$
\nabla p_\pm(z)
=
\mp \frac{2|z|^2 z}{\sqrt{1 - |z|^4}},
$ 
it follows that $\sqrt{1 + |\nabla p_\pm(z)|^2}
=
\sqrt{\frac{1 - |z|^4 + 4|z|^6}{1 - |z|^4}},$

%%%%%%%%%%%%%%%%%%%
%%%%%%%%%%%%%%%%%%%
%%%%%%%%%%%%%%%%%%%
%%%%%%%%%%%%%%%%%%%
and therefore
\begin{equation} 
\label{expression-of-integral-form-of-dH-2n}
\int_{S_{K,1}} f \, d\mathcal{H}^{2n}_K 
=
\int_{|z|<1}
\left[
f(z,\sqrt{1-|z|^4})
+
f(z,-\sqrt{1-|z|^4})
\right]
\sqrt{\frac{1 - |z|^4 + 4|z|^6}{1 - |z|^4}}\, dz.
\end{equation}
%%%%%%%%%%%%%%%%%%%%%%%%%%%%%%%%%%%%%%%%%%%%%%%%%%%%%%%%%%%%%%%%%%%%%%%%%%%%%%%%%%%%%%%%%%%%%%%%%%%%%%%%%%%%%%%%%%%%%%%%%%%%%%%%%%%%%%%%%%%%%%%%%%%%%%%%%%%%%%%%%%%%%%%%%%%%%%%%%%%%%%%%%%%%%%%%%%%%%%%%%%%%%

\subsubsection{Integral with respect to $d\sigma^{2n}_K$ and the polar decomposition}
Rewriting the density in the right hand side of the integral in \eqref{expression-of-integral-form-of-dH-2n} as 
\[
\sqrt{\frac{1 - |z|^4 + 4|z|^6}{1 - |z|^4}}
=
\frac{\sqrt{4(1-|z|^4) + 16|z|^6}}{2\sqrt{1-|z|^4}}, 
\]
and using $t^2 = 1 - |z|^4$ on $S_{K,1}$, we define the measure $d\sigma^{2n}_K$ supported on $S_{K,1}$ by
\[
d\sigma^{2n}_K
=
\frac{d\mathcal{H}^{2n}_K}{|\nabla N|}=\frac{4\,d\mathcal{H}^{2n}_K}{|\nabla \phi|}.
\]

Using \eqref{expression-of-integral-form-of-dH-2n}, for any continuous function $f$ on $\mathbb{H}^n$, we have
\begin{align*}
\int_{S_{K,1}} f \, d\sigma^{2n}_K
%= \int_{S_{K,1}} \frac{f}{|\nabla N|}\, d\mathcal{H}^{2n}
= \int_{|z|<1}
\left[
f(z,\sqrt{1-|z|^4})
+
f(z,-\sqrt{1-|z|^4})
\right]
\frac{2\,dz}{\sqrt{1-|z|^4}}.
\end{align*}
Similarly, for the $r$-radius Kor\'anyi sphere, one can determine the the surface measure $d\mathcal{H}^{2n}_{K,r}$ supported on the $r$-radius Kor\'anyi sphere $S_{K,r}$, and it follows that the same takes the  form 
$$
\int_{S_{K,r}} f \, d\mathcal{H}^{2n}_{K,r}
=
\int_{|z|<r}
\left[
f(z,\sqrt{r^4-|z|^4})
+
f(z,-\sqrt{r^4-|z|^4})
\right]
\sqrt{\frac{r^4 - |z|^4 + 4|z|^6}{r^4 - |z|^4}}\, dz.
$$ 

It is obvious from the construction that the measure $d\sigma^{2n}_K$ satisfies the scaling relation
\begin{align}
\label{Koranyi-measure-scaling-relation}
d\sigma_{K,r}^{2n} = r^{\mathcal{Q}-1}\tilde{\delta}_r(d\sigma^{2n}_K), 
\end{align}
which is same as saying that 
\begin{align}
\label{Koranyi-measure-scaling-relation-integral-form}
\int_{S_{K,r}} f \, d\sigma^{2n}_{K,r} = r^{\mathcal{Q}-1}\int_{S_{K,1}} f\circ \tilde{\delta}_r(z,t)\, d\sigma^{2n}_K(z,t).
\end{align}
For the sake of clarity, let us see the details pertaining the same. For $r>0$, consider $N_r(z,t)=(|z|^4 + t^2 )^\frac{1}{4} -r$ so that the associated measure $d\sigma^{2n}_{K,r}$ supported on $S_{K,r}$ is
\begin{equation} 
\label{Expression-for-the-measure-on-Koranyi-sphere}
d \sigma^{2n}_{K,r} = \frac{d\mathcal{H}^{2n}_{K,r}}{|\nabla N_r|}=\frac{4r^3\,d\mathcal{H}^{2n}_{K,r}}{|\nabla \phi|}. 
\end{equation}
Now, 
\begin{align*}
\int_{S_{K,r}} f \, d\sigma^{2n}_{K,r}
&=
\int_{|z|<r}
\left[
f(z,\sqrt{r^4-|z|^4})
+
f(z,-\sqrt{r^4-|z|^4})
\right]
\frac{2r^3\,dz}{\sqrt{r^4-|z|^4}}\\
&=\int_{|z|<1}
\left[
f\circ \tilde{\delta}_r(z,\sqrt{1-|z|^4})
+
f(z,-\sqrt{1-|z|^4})
\right]
\frac{2r^{\mathcal{Q}-1}\,dz}{\sqrt{1-|z|^4}}\\
&=r^{\mathcal{Q}-1}\int_{S_{K,1}} f \circ \tilde{\delta}_r \, d\sigma^{2n}_K, 
\end{align*}
which is precisely \eqref{Koranyi-measure-scaling-relation-integral-form}. 
%%%%%%%%%%%%%%%%%%%%%
%%%%%%%%%%%%%%%%%%%%%
%%%%%%%%%%%%%%%%%%%%%
%%%%%%%%%%%%%%%%%%%%%
%%%%%%%%%%%%%%%%%%%%%

\subsection{Radon measure \texorpdfstring{$d\mu^{4n+1}_K$}{} on the product Kor\'anyi sphere \texorpdfstring{$\Sigma_{K,1}$}{}} \label{exp-surface-measure-on-the-product-Koranyi-sphere}
In this subsection, we will explicitly calculate the measure supported on the product Kor\'anyi sphere $\Sigma_{K,1}$ for the case $\mathbf{a}=(4,4)$, which will be useful in our context. Let $d\mathcal{H}_{K,r}^{4n+1}$ denote the Hausdorff measure (surface measure) supported on $\Sigma_{K,r}$, we shall shortly discuss it in detail in Subsection \ref{subsec:hausdorff-measure-product}, but for now let us assume it. As in the Kor\'anyi case, that is, in Subsection \ref{exp-surface-measure-on-the-linear-Koranyi-sphere}, this measure $d\mathcal{H}_{K,r}^{4n+1}$ also does not directly provide the polar decomposition \eqref{expres-polar-decoposition-producttt}, nor is it useful in decomposing the higher dimensional sphere into a lower dimensional sphere and a ball. Instead, as in Subsection \ref{exp-surface-measure-on-the-linear-Koranyi-sphere}, here also we need to modify it to define the appropriate Radon measure $d\mu_{K,r}^{4n+1}$ as follows. 

\medskip 
Consider the smooth map $\Psi : \mathbb{H}^n \times \mathbb{H}^n \to \mathbb{R}$ by $\Psi(g_1,g_2)=|g_1|_K^4+|g_2|_K^4$. Clearly, $\nabla\Psi \neq 0$ on $\Sigma_{K,r}$. Making use of the coarea formula and an appropriate change of variables, we have 
\begin{align*}
\int_{\mathbb{H}^n \times \mathbb{H}^n} f(g_1,g_2)dg_1dg_2 
&= \int_0^\infty \int_{\{\Psi=s\}} \frac{f(g_1,g_2)}{|\nabla \Psi (g_1,g_2)|} \, d\mathcal{H}_{K,s^\frac{1}{4}}^{4n+1}(g_1,g_2) \, ds \\ 
&= \int_0^\infty \int_{\{\Psi=r^4\}} \frac{f(g_1,g_2)}{|\nabla \Psi (g_1,g_2)|}4r^3 \, d\mathcal{H}_{K,r}^{4n+1}(g_1,g_2) \, dr, 
\end{align*}
and therefore if we define the Radon measure $d\mu_{K,r}^{4n+1}$ supported on $\Sigma_{K,r}$ by 
\begin{align} 
\label{surface-measure-on-r-radius-product-Koranyi-sphere}
d\mu_{K,r}^{4n+1} = 4r^3 \, |\nabla \Psi|^{-1} \, d\mathcal{H}^{4n+1}_{K,r}, 
\end{align} 
then we have the polar decomposition of the Haar measure on $\mathbb{H}^n\times \mathbb{H}^n$ as
$$\int_{\mathbb{H}^n  \times \mathbb{H}^n} f(g_1,g_2)\,dg_1\,dg_2=\int_0^\infty\int_{\Sigma_{K,r}} f(g_1,g_2)\,d\mathcal{\mu}_{K,r}^{4n+1}(g_1,g_2) \,dr.
$$
As earlier, we write $d\mathcal{H}_K^{4n+1} = d\mathcal{H}_{K,1}^{4n+1}$ and $d\mu^{4n+1}_K = d\mu_{K,1}^{4n+1}$. 

\medskip
In the next subsection we shall compute the explicit form of this measure. Even though this explicit form is not needed in the present work, we still do so as it motivated our progress on the problem in the initial phase and we also believe that this could be of independent interest in other works involving the analysis on the product Kor\'anyi spheres.
%%%%%%%%%%%%%%%%%%%%
%%%%%%%%%%%%%%%%%%%%
%%%%%%%%%%%%%%%%%%%%
%%%%%%%%%%%%%%%%%%%%

\subsubsection{Hausdorff measure \texorpdfstring{$d\mathcal{H}_K^{4n+1}$}{} on \texorpdfstring{$\Sigma_{K,1}$}{}} 
\label{subsec:hausdorff-measure-product}

We shall follow the same technique as we did for the Kor\'anyi sphere case in Subsection \ref{exp-surface-measure-on-the-linear-Koranyi-sphere}, and show that the surface measure supported on the unit product Kor\'anyi sphere $\Sigma_{K,1}$ can be given in terms of the volume element
\begin{equation}\label{exp-volume-element-form-for-unit-Koranyi sphere-1}
\frac{\sqrt{4(1-|g_1|_K^4-|z_2|^4)+|\nabla_{g_1}|g_1|_K^4|^2+16|z_2|^6}}{\sqrt{4(1-|g_1|_K^4-|z_2|^4)}}\,dg_1\,dz_2
\end{equation}
on the domain $U$ in $\mathbb{R}^{4n+1}$, where $U=\{(g_1,z_2)\in \mathbb{R}^{4n+1}: 0 \leq |g_1|_K^4+|z_2|^4<1\}.$ 

\medskip 
Recall that $\mathbb{H}^n \times \mathbb{H}^n$ is topologically nothing but $\mathbb{R}^{4n+2}$. Define $\rho_+:U\subset\mathbb{R}^{4n+1} \rightarrow \mathbb{R}$ by 
\[
\rho_+(g_1,z_2)=\sqrt{1-|g_1|_K^4-|z_2|^4}, 
\]
%and its complimentary map $\rho_-$ on $U$ as $\rho_- = -\rho_+$, 
so that 
$$
\Sigma_{K,1}=\text{Graph}(\rho_+) \cup \text{Graph}(-\rho_+) \cup \{4n\text{-dimensional surface\}.}
$$
Here also we shall ignore the final part in the above decomposition as it is of measure zero. 

\medskip 
Let us first focus on the graph $\Gamma_+$ of $\rho_+$, that is, 
\[
\Gamma_+=\{(g_1,z_2,\rho_+(g_1,z_2)):(g_1,z_2)\in U\},
\]
and consider the parametrization $(\phi_+,U)$ for $\Gamma_+$, where $\phi_+:U\rightarrow\mathbb{R}^{4n+2}$ is defined by 
\[
\phi_+(g_1,z_2)=(g_1,z_2,\sqrt{1-|g_1|_K^4-|z_2|^4}).
\]
The surface element is given by the formula $\sqrt{1+|\nabla\rho_+|^2} \, dg_1 \, dz_2.
$
If we denote $g_1=(z_1,t_1)$, then the Euclidean gradient $\nabla \rho_+ = \left( \nabla_{z_1}, \partial_{t_1}, \nabla_{z_2} \right)\rho_+$ turns out to be 
\[
\nabla \rho_+(z_1,t_1,z_2)
= \frac{-1}{\sqrt{1 - |z_1|^4 - t_1^2 - |z_2|^4}} \left(
2|z_1|^2 z_1,\;
t_1,\;
2|z_2|^2 z_2
\right), 
\]
so that 
\begin{align*}
\sqrt{1+|\nabla\rho_+|^2} \, dg_1 \, dz_2 
&=\sqrt{1+\left |\frac{(2|z_1|^2z_1,t_1,2|z_2|^2z_2)}{\sqrt{1-|z_1|^4-t_1^2-|z_2|^4}}\right |^2} dg_1 \, dz_2 \\ 
&=\sqrt{1+\frac{|\nabla_{g_1}|g_1|_K^4|^2+16|z_2|^6}{4(1-|g_1|_K^4-|z_2|^4)}} \, dg_1 \, dz_2 \\ 
&=\frac{\sqrt{4(1-|g_1|_K^4-|z_2|^4)+|\nabla_{g_1}|g_1|_K^4|^2+16|z_2|^6}}{\sqrt{4(1-|g_1|_K^4-|z_2|^4)}} \, dg_1 \, dz_2,
\end{align*}
where we made use of the relation $\nabla_{g_1}|g_1|_K^4=(4|z_1|^2z_1,2t_1)$ in the second last step. 

\medskip 
Summarising, if $f$ is a continuous function defined on the surface $\Gamma_+$, then 
\[
\int_{\Gamma_+}f \, d\mathcal{H}_+^{4n+1}=\int_Uf(g_1,z_2,\rho_+(g_1,z_2))\sqrt{1+|\nabla\rho_+|^2} \,dg_1\,dz_2,
\]
where $d\mathcal{H}_+^{4n+1}$ is the surface measure on $\Gamma_+$.

\medskip
We will have a similar expression for the surface measure $d\mathcal{H}_-^{4n+1}$ on the surface $\Gamma_-$ with respect to $-\rho_+.$ Therefore, if we define 
\[
d\mathcal{H}^{4n+1}_K(p)=
\begin{cases}
d\mathcal{H}_+^{4n+1}(p), & p \in \Gamma_+, \\
d\mathcal{H}_-^{4n+1}(p), & p \in \Gamma_-, 
\end{cases}
\]
then $d\mathcal{H}^{4n+1}_K$ is the induced surface measure or the Hausdorff measure on $\Sigma_{K,1}.$

%%%%%%%%%%%%%%%%%%%%%%%%%%%%%%%%%%%%%%%%%%%%%%%%%%%%%%%%%%%%%%%%%%%%%%%%%%%%%%%%%%%%%%%%%%%%%%%%%%%%%%%%%%%%%%%%%%%%%%%%%%%%%%%%%%%%%%%%%%%%%%%%%%%%%%%%%%%%%%%%%%%%%%%%%%%%%%%%%%%%%%%%%%%%%%%%
\subsubsection{Radon measure \texorpdfstring{$d\mu^{4n+1}_K$}{} on \texorpdfstring{$\Sigma_{K,1}$}{} suitable for polar decomposition} 
\label{dilation-structure-of-surface-measure-on-the-product-Koranyi-sphere}
Altogether, if we let $T(g_1,g_2) = \Psi(g_1, g_2)^{1/4} - 1 = (|g_1|_K^4+|g_2|_K^4)^{\frac14} - 1$, so that on $\Sigma_{K,1}$ we have 
$$|\nabla T| = \frac{1}{4} \sqrt{4(1-|g_1|_K^4-|z_2|_K^4)+|\nabla_{g_1}|g_1|_K^4|^2+16|z_2|^6},$$
and with that we define the Radon measure $d\mu^{4n+1}_K$ supported on $\Sigma^{K,1}$ by 
\[
d \mu^{4n+1}_K = |\nabla T|^{-1} \, d\mathcal{H}^{4n+1}_K = 4 |\nabla \Psi|^{-1} \, \,d\mathcal{H}^{4n+1}_K. 
\]
More precisely, for a continuous function $f$ defined on $\mathbb{H}^n\times \mathbb{H}^n$, we have
\begin{align} 
\label{expression-for-integral-form-for-product-koranyi}
\int_{\Sigma_{K,1}}f \, d{\mu}^{4n+1}_K
=& \int_U f(g_1, z_2, \rho_+(g_1,z_2)) \frac{2\,\,dg_1\,dz_2}{\sqrt{(1-|g_1|_K^4-|z_2|^4)}} \\ 
&+ \int_U f(g_1, z_2, -\rho_+(g_1,z_2)) \frac{2\,\,dg_1\,dz_2}{\sqrt{(1-|g_1|_K^4-|z_2|^4)}}. \notag
\end{align}
 
We have thus derived the measure $d\mu^{4n+1}_K$ which is supported on the unit product Kor\'anyi sphere. Let $d\mu_{K,r}^{4n+1}$ denote the analogous measure on the $r$-radius sphere $\Sigma_{K,r}$. Naturally, the following scaling relation holds true: 
\begin{align} 
\label{expression-of-dilation-of-product-surface-measure}
d\mu_{K,r}^{4n+1}=r^{4n+3}\,\delta_r(d\mu^{4n+1}_K), \quad \text{for } r>0,
\end{align}
that is, 
\begin{align}
\label{identity-with-dilation-of-product-surface-measure} 
\int_{\Sigma_{K,r}}f(g_1,g_2)\,d\mu_{K,r}^{4n+1}(g_1,g_2) = r^{4n+3}\int_{\Sigma_{K,1}}f(\delta_r(g_1,g_2))\,d\mu^{4n+1}_K(g_1,g_2).
\end{align}

Let us see the details leading to the above relation. For $r>0$, consider the map $T_r(z,t)=(|g_1|^4 + |g_2|^4 )^\frac{1}{4} -r$  and $ U_r=\{(g_1,z_2)\in \mathbb{R}^{4n+1}: 0\leq|g|_K^4+|z_2|^4<r^4\}$. Then, repeating the calculations as above, one obtains the corresponding volume element 
\begin{align*}
\frac{\sqrt{4(r^4-|g_1|_K^4-|z_2|^4)+|\nabla_{g_1}|g_1|_K^4|^2+16|z_2|^6}}{\sqrt{4(r^4-|g_1|_K^4-|z_2|^4)}}\,dg_1\,dz_2,
\end{align*}
which gives the associated measure $d\mu^{4n+1}_{K,r}$ supported on $\Sigma_{K,r}$ as
\begin{align}\label{Radon-measure-product-Koranyi-type}
d\mu^{4n+1}_{K,r} = |\nabla T_r|^{-1} d\mathcal{H}^{4n+1}_{K,r} = 4r^3 \, |\nabla \Psi|^{-1} \, d\mathcal{H}^{4n+1}_{K,r}.
\end{align}
Hence, for a continuous function $f$ on $\mathbb{H}^n$, we have
\begin{align*}
\int_{\Sigma_{K,r}} f \, d\mu^{4n+1}_{K,r} 
&= \int_{U_r}\frac{f(g_1,z_2,\sqrt{(r^4-|g_1|_K^4-|z_2|^4)})}{\sqrt{(r^4-|g_1|_K^4-|z_2|^4)}}2r^3\,dg_1\,dz_2\\
& \quad + \int_{U_r} \frac{f(g_1,z_2,-\sqrt{(r^4-|g_1|_K^4-|z_2|^4)})}{\sqrt{(r^4-|g_1|_K^4-|z_2|^4)}}2r^3\,dg_1\,dz_2 \\ 
&= \int_{U} \frac{f \circ \tilde{\delta}_r\,(g_1, z_2, \sqrt{(1-|g_1|_K^4-|z_2|^4}))}{\sqrt{(1-|g_1|_K^4-|z_2|^4)}} 2r^{4n+3}\,dg_1\,dz_2\\
& \quad + \int_{U}\frac{f \circ \tilde{\delta}_r \,(g_1,z_2,-\sqrt{(1-|g_1|_K^4-|z_2|^4}))}{\sqrt{(1-|g_1|_K^4-|z_2|^4)}} 2r^{4n+3}\,dg_1\,dz_2,
\end{align*}
which along with \eqref{expression-for-integral-form-for-product-koranyi} implies \eqref{identity-with-dilation-of-product-surface-measure}. 

%%%%%%%%%%%%%%%%%%
%%%%%%%%%%%%%%%%%%
%%%%%%%%%%%%%%%%%%
%%%%%%%%%%%%%%%%%%

\subsection{Radon measure \texorpdfstring{$d\mu^\mathbf{a}_K$}{} on the degenerate product Kor\'anyi surface
\texorpdfstring{$\mathcal{S}_K^{\mathbf{a}}$}{}} 
\label{exp-surface-measure-on-the-degenerate-product-Koranyi-sphere}

Following the technique of Subsection \ref{exp-surface-measure-on-the-product-Koranyi-sphere}, one can compute the Hausdorff surface measure $d\mathcal{H}_K^\mathbf{a}$ on $\mathcal{S}_K^\mathbf{a}.$ We do not write the details for the same as it is technically involved and we don't really need its explicit form in our proofs in the subsequent sections. 

\medskip 
Now, let us recall the following formula \cite[p.~136]{The-analysis-of-linear-partial-differential}. Let $\Omega$ be an $(n-1)$-dimensional hypersurface in $\mathbb{R}^n$ given by $\Omega=\{x\in \mathbb{R}^n: \Phi(x)=0\}$. If we have $\nabla\Phi(w) \neq 0$ for any $w \in \Omega$, then 
\begin{equation} 
\label{coarea-formula-statement}
\int_{\Omega} f \, |\nabla\Phi|^{-1} \, d\mathcal{H}^{n-1} = \int_{\mathbb{R}^n} f(x) \, \delta(\Phi(x)) \,dx,
\end{equation}
for every $f \in \mathcal{S} (\mathbb{R}^n)$, where $d\mathcal{H}^{n-1}$ denote the induced surface measure (Hausdorff measure) on $\Omega$.

\medskip 
In particular, let us consider the map $
\Phi^K: \mathbb{H}^n\times \mathbb{H}^n \longrightarrow \mathbb{R} 
$ given by $\Phi^K(g_1,g_2) = |g_1|_K^{\mathbf{a}_1} + |g_2|_K^{\mathbf{a}_2} - 1$. Note that $\mathcal{S}_K^\mathbf{a} = (\Phi^K)^{-1} (0)$. Also, for any $\mathbf{a}_1, \, \mathbf{a}_2 \geq 2$, we have 
\[
\nabla \Phi^K(g_1,g_2) = (\mathbf{a}_1 \, |g_1|_K^{\mathbf{a}_1-1} \nabla (|g_1|_K), \, \mathbf{a}_2 \, |g_2|_K^{\mathbf{a}_2-1} \nabla (|g_2|_K)) \neq 0, 
\]
for a.e. $ (g_1, g_2) \in \mathcal{S}_K^\mathbf{a}$. 

\medskip 
Thus, it follows from \eqref{coarea-formula-statement} that with the Radon measure $d\mu^\mathbf{a}_K$ on $\mathcal{S}_K^\mathbf{a}$ given by 
\begin{align}
\label{measure-degenerate-product-Koranyi}
d\mu^\mathbf{a}_K = |\nabla\Phi^K|^{-1} \, d\mathcal{H_K^{\mathbf{a}}}, 
\end{align}
we have 
\begin{align}
\label{formula-pull-back-of-the-surface-measure-Koranyi-case}
\int_{\mathcal{S}^\mathbf{a}} H \, d\mu^\mathbf{a}_K = \int_{\mathbb{R}^{4n+2}} H(g_1,g_2) \, \delta(\Phi^K(g_1,g_2)) \, dg_1 \, dg_2, 
\end{align}
for every $H\in \mathcal{S}(\mathbb{R}^{4n+2})$. 

%%%%%%%%%%%%%%%%%%%%%%%%%%%
%%%%%%%%%%%%%%%%%%%%%%%%%%%
%%%%%%%%%%%%%%%%%%%%%%%%%%%
%%%%%%%%%%%%%%%%%%%%%%%%%%%

\section{Proof of Theorem \ref{Main-theorem-Koranyi-type-for-degenrate}} 
\label{sec:proof-thm-full-max-koranyi}

One of the main ingredients in the technique developed in \cite{Maximal-estimates-for-the-bilinear-spherical-averages} was a clever use of the slicing argument. Given a continuous function $f$ on $\mathbb{R}^{2n}$, one can express
\begin{equation}\label{2.2}
\int_{S^{2n-1}} f \, d\sigma_{2n-1} = \int_{B^n(0,1)} \int_{S^{n-1}} f \left(x,\sqrt{1-|x|^2}\,y\right)
(1-|x|^2)^{\frac{n-2}{2}}
\, d\sigma_{n-1}(y)\, dx, 
\end{equation}
where $d\sigma_{n-1}$ is the standard surface measure on the unit sphere $S^{n-1}$. This identity allows one to have a pointwise bound of the bilinear maximal function by product of the well known Hardy--Littlewood maximal function and Stein's spherical maximal function on $\mathbb{R}^n$, which in turn provides boundedness estimates for the bilinear spherical maximal operator on $\mathbb{R}^n$. 
\medskip

This technique has been successfully used in the subsequent works as well. In the following subsection, we show that in our case also we can decompose the higher-dimensional sphere or more generally the degenerate surfaces $\mathcal{S}^\mathbf{a}_K$ into a lower-dimensional sphere (Kor\'anyi sphere) and a Kor\'anyi ball which can be seen as a slicing argument for this setup. 

%%%%%%%%%%%%%%%%%%%%%%%%%%%
%%%%%%%%%%%%%%%%%%%%%%%%%%%
%%%%%%%%%%%%%%%%%%%%%%%%%%%
%%%%%%%%%%%%%%%%%%%%%%%%%%%

\subsection{Slicing argument for the degenerate product Kor\'anyi surface \texorpdfstring{$\mathcal{S}_K^{\mathbf{a}}$}{}} 
\label{finding-slicing-argument-form-for-product-Koranyi} 

As seen in Subsection \ref{exp-surface-measure-on-the-degenerate-product-Koranyi-sphere}, with the map $
\Phi^K: \mathbb{H}^n\times \mathbb{H}^n \longrightarrow \mathbb{R} 
$ given by $\Phi^K(g_1,g_2) = |g_1|_K^{\mathbf{a}_1} + |g_2|_K^{\mathbf{a}_2} - 1$, for any $H\in \mathcal{S}(\mathbb{R}^{4n+2})$, we have from \eqref{formula-pull-back-of-the-surface-measure-Koranyi-case} that 
\begin{align*}
\int_{\mathcal{S}^\mathbf{a}_K} H \, d\mu_K^{\mathbf{a}}
&=\int_{\mathbb{R}^{4n+2}} H(g_1,g_2) \, \delta(\Phi^K(g_1,g_2)) \, dg_1 \, dg_2 \\ 
& = \int_{|g_1|_K<1} \int_{\mathbb{R}^{2n+1}} H(g_1,g_2) \, \delta(\Phi^K(g_1,\cdot)(g_2))\,dg_2\,dg_1 \\
&= \int_{|g_1|_K<1} \int_{S_{K,(1-|g_1|_K^{\mathbf{a}_1}})^\frac{1}{\mathbf{a}_2}} H(g_1,g_2) \, |\nabla_{g_2} \Phi^K(g_1,\cdot)(g_2)|^{-1} \, d\mathcal{H}^{2n}_{K,(1-|g_1|_K^{\mathbf{a}_1})^\frac{1}{\mathbf{a}_2}}(g_2) \, dg_1,
\end{align*}
where the last step also follows from the above formula \eqref{coarea-formula-statement} applied to the inner integral in $g_2$-variable (keeping $g_1$-variable fixed). Here, $S_{K,(1-|g_1|_K^{\mathbf{a}_1}})^\frac{1}{\mathbf{a}_2}$ is nothing but the $(1-|g_1|_K^{\mathbf{a}_1})^\frac{1}{\mathbf{a}_2}$-radius Kor\'anyi sphere on $\mathbb{H}^n.$ 

\medskip 
It follows from \eqref{Expression-for-the-measure-on-Koranyi-sphere} that the final term in the above identity equals to 
\begin{align*}
\int_{|g_1|_K<1}  \int_{S_{K,(1-|g_1|_K^{\mathbf{a}_1}})^\frac{1}{\mathbf{a}_2}} H(g_1,g_2) \, \frac{\left|\nabla_{g_2}|g_2|_K^4\right|}{4(1-|g_1|_K^{\mathbf{a}_1})^\frac{3}{\mathbf{a}_2}|\nabla_{g_2} \Phi^K(g_1,\cdot)(g_2)|}\,d{\sigma}^{2n}_{K,(1-|g_1|_K^{\mathbf{a}_1})^\frac{1}{\mathbf{a}_2}}(g_2)\, dg_1.
\end{align*}
Thus, simplifying and using the scaling relation of the measure $d\sigma^{2n}_K$ as in \eqref{Koranyi-measure-scaling-relation-integral-form}, we obtain
\begin{align} 
\label{expression-of-slicing-arguement} 
\int_{\mathcal{S}_K^\mathbf{a}} H \, d\mu_K^{\mathbf{a}} = \frac{1}{\mathbf{a}_2} \int_{B_K(0,1)} \int_{S_{K,1}}  H(g_1,\tilde{\delta}_{(1-|g_1|_K^{\mathbf{a}_1})^{\frac{1}{\mathbf{a}_2}}}, g_2^{-1}) \, d\sigma_K^{2n}(g_2) \, (1-|g_1|_K^{\mathbf{a}_1})^{\frac{2n+2-\mathbf{a}_2}{\mathbf{a}_2}} \, dg_1.
\end{align} 

\begin{remark} 
\label{remark-for-slicing argument-for-a_1=4=a_2}
For $\mathbf{a}_1 = \mathbf{a}_2=4,$ that is, when working with the product Kor\'anyi sphere $\Sigma_{K,1}$, identity \eqref{expression-of-slicing-arguement} reduces to 
\begin{align}\label{expression-of-slicing-arguement-a_1=4=a_2}
\int_{\Sigma_{K,1}} H \, d\mu^{4n+1}_K = \int_{B_K(0,1)} \int_{S_{K,1}} H(g_1,\tilde{\delta}_{(1-|g_1|_K^4)^{\frac{1}{4}}}g_2) \, d\sigma^{2n}_K(g_2) \, (1-|g_1|_K^4)^{\frac{2n-2}{4}} \, dg_1,
\end{align}
which is the slicing identity for the $1$-radius product Kor\'anyi sphere $\Sigma_{K,1}$. 
\end{remark}

%%%%%%%%%%%%%%%%%%%%
%%%%%%%%%%%%%%%%%%%%
%%%%%%%%%%%%%%%%%%%%
%%%%%%%%%%%%%%%%%%%%
%%%%%%%%%%%%%%%%%%%%
%%%%%%%%%%%%%%%%%%%%

\subsection{Proof of Theorem \ref{Main-theorem-Koranyi-type-for-degenrate}} 
\label{linearization-technique}
%%%%%%%%%%%%%%%%%%%%%%%%%%%%%%%%%%%%%%%%%%%%%%%%%%%%%%%%%%%%%%%%%%%%%%%%%%%%%%%%%%%%%%%%%%%
Using identity \eqref{expression-of-slicing-arguement}, we can see that the biparametric bilinear averages satisfies the following equality.
\begin{align}\label{equivalent expression-degenrate-Koranyi-2}
& \mathcal{A}^{\mathbf{a},K}_\mathbf{s}(F,G)(g) \\ 
&= \int_{\mathcal{S}^\mathbf{a}_K} F(g\cdot \tilde{\delta}_{s_1}{g_1}^{-1}) \, G(g\cdot \tilde{\delta}_{s_2}{g_2}^{-1}) \, d\mu^\mathbf{a}_K(g_1,g_2) \notag \\ 
&= \frac{1}{\mathbf{a}_2} \int_{B_K({0,1})} F(g\cdot \tilde{\delta}_{s_1}{g_1}^{-1}) \int_{S_{K,1}} G(g\cdot \tilde{\delta}_{s_2} \circ \tilde{\delta}_{(1-|g_1|_K^{\mathbf{a}_1})^{\frac{1}{\mathbf{a}_2}}} {g_2}^{-1}) \, d\sigma^{2n}_K(g_2) \, (1-|g_1|_K^{\mathbf{a}_1})^{\frac{2n+2-\mathbf{a}_2}{\mathbf{a}_2}} \, dg_1 \notag \\ 
%&= \frac{1}{\mathbf{a}_2} \int_{{|g_1|_K<1}}f(x\cdot\tilde{\delta}_{s_1}g_1^{-1})\left\{ \int_{S_{K,1}}g(x\cdot\tilde{\delta}_{s_2(1-|g_1|_K^{\mathbf{a}_1})^\frac{1}{\mathbf{a}_2}}g_2^{-1})\,d\mathcal{\sigma}^{2n}(g_2) \right\} (1-|g_1|_K^{\mathbf{a}_1})^\frac{(2n+2-\mathbf{a}_2)}{\mathbf{a}_2} \, dg_1\notag\\
&= \frac{1}{\mathbf{a}_2} \int_{{|g_1|_K<1}} F(g\cdot\tilde{\delta}_{s_1}g_1^{-1}) \, A^K G(g, s_2(1-|g_1|_K^{\mathbf{a}_1})^\frac{1}{\mathbf{a}_2})(1-|g_1|_K^{\mathbf{a}_1})^\frac{(2n+2-\mathbf{a}_2)}{\mathbf{a}_2}\,dg_1, \notag
\end{align}
where 
$A^K G(g, r) := \int_{S_{K,1}} G(g\cdot\tilde{\delta}_{r} g_2^{-1}) \, d\mathcal{\sigma}_K^{2n}(g_2).$

\medskip 
Interchanging the roles of $F$ and $G$, we also have the following symmetric expression: 
\begin{align}\label{equivalent expression-degenrate-Koranyi}
\mathcal{A}^{\mathbf{a},K}_\mathbf{s}(F,G)(g)
=&\frac{1}{\mathbf{a}_1} \int_{{|g_2|_K<1}} G(g\cdot\tilde{\delta}_{s_2}g_2^{-1}) \, A^K F(g,s_1(1-|g_2|_K^{\mathbf{a}_2})^\frac{1}{\mathbf{a}_1}) \, (1-|g_2|_K^{\mathbf{a}_2})^\frac{(2n+2-\mathbf{a}_1)}{\mathbf{a}_1} \, dg_2.
\end{align}
Now, whether $\mathbf{a}_1$, $\mathbf{a}_2$ are less than or greater than the homogeneous dimension $2n+2$, we will have following cases similar to \cite{Lee_Shuin} and follow their same strategy.

%%%%%%%%%%%%%%%%%%%%%%%%%%%%%%%%%%%%%%%%%%%%%%%%%%%%%%%%%%%%

\medskip
In particular, when $\max\{\mathbf{a}_1,\mathbf{a}_2\}\leq 2n+2$, the slicing identity \eqref{expression-of-slicing-arguement} allows one to have a pointwise bound of the bilinear maximal function $\mathcal{M}^{\mathbf{a},K}_{full}$ on $\mathbb{H}^n \times \mathbb{H}^n$ by the product of the spherical maximal function $M_{\mathbb{H}^n}$ (see \eqref{def:linear-Koranyi-sph-maximal}) and the centered Hardy--Littlewood maximal function $M^{K}_{\mathrm{HL}}$, as we shall see shortly in Lemma \ref{lemma-for-pointwise-control-byknown-maximal-operator}. The operator $M^{K}_{\mathrm{HL}}$ associated with the Kor\'anyi balls on $\mathbb{H}^n$ is defined by
\begin{equation}
M^{K}_{\mathrm{HL}} F(g)
= \sup_{r>0} \frac{1}{|B_K(g,r)|} \left| \int_{B_K(g,r)} F(y) \, dy \right|.
\end{equation}
It is well known (see, for example, \cite{Li-Hong-Qian-Hardy-Littlewood-maximal-function-on-the-Heisenberg-group}) that $M^{K}_{\mathrm{HL}}$ is of weak-type $(1,1)$ as well as strong-type $(p,p)$ for $1<p\leq\infty$, that is, 
\begin{align*}
|\{g\in\mathbb{H}^n : M^{K}_{\mathrm{HL}} F(g)>\lambda\}| \lesssim \frac{1}{\lambda} \|F\|_{L^1(\mathbb{H}^n)} \quad \text{and} \quad \|M^{K}_{\mathrm{HL}}F\|_{L^p(\mathbb{H}^n)} \lesssim \|F\|_{L^p(\mathbb{H}^n)}. 
\end{align*}

%%%%%%%%%%%%%%%%%%%%%%
%%%%%%%%%%%%%%%%%%%%%%
%%%%%%%%%%%%%%%%%%%%%%
%%%%%%%%%%%%%%%%%%%%%%

\begin{lemma} 
\label{lemma-for-pointwise-control-byknown-maximal-operator}
Let $n\geq1$. For $\mathbf{a}_1,\mathbf{a}_2\leq 2n+2$, we have 
\begin{align*} 
%\label{inequality-in-the-lemma}
\mathcal{M}^{\mathbf{a},K}_{full}(F,G)(g) \lesssim M^K_{HL} |F|(g) \, M_{\mathbb{H}^n}|G|(g) 
\quad \text{and} \quad 
\mathcal{M}^{\mathbf{a},K}_{full}(F,G)(g) \lesssim M_{\mathbb{H}^n}|F|(g) \, M^K_{HL} |G|(g). 
\end{align*} 
\end{lemma}

\begin{proof}
Using \eqref{equivalent expression-degenrate-Koranyi-2}, we can see that the biparametric bilinear averages $\mathcal{A}^{\mathbf{a},K}_\mathbf{s}(f,g)$ are controlled by
\begin{align*}
& \left| \mathcal{A}^{\mathbf{a},K}_\mathbf{s}(F,G)(g)\right| \\ 
& \leq \frac{1}{\mathbf{a}_2} \int_{B_K({0,1})}\left|F(g\cdot \tilde{\delta}_{s_1}{g_1}^{-1})\right| \int_{S_{K,1}}\left|G(g\cdot \tilde{\delta}_{{s_2}(1-|g_1|_K^{\mathbf{a}_1})^{\frac{1}{\mathbf{a}_2}}}{g_2}^{-1})\right| \, d\sigma^{2n}_K(g_2) \, (1-|g_1|_K^{\mathbf{a}_1})^{\frac{2n+2-\mathbf{a}_2}{\mathbf{a}_2}} \,dg_1 \\
& \lesssim M_{\mathbb{H}^n}|G|(g)\int_{B_K({0,1})}\left|F(g\cdot \tilde{\delta}_{s_1}{g_1}^{-1})\right|(1-|g_1|_K^{\mathbf{a}_1})^{\frac{2n+2-\mathbf{a}_2}{\mathbf{a}_2}}\,dg_1.
\end{align*}
Now note that as $\mathbf{a}_2\leq2n+2,$ we have $(1-|g_1|_K^4)^{\frac{2n+2-\mathbf{a}_2}{\mathbf{a}_2}}\leq 1$ on $B_K(0,1).$ Hence, we get
\begin{align*}
\left|\int_{\mathcal{S}^\mathbf{a}_K} F(g\cdot \tilde{\delta}_{s_1}{g_1}^{-1}) G(g\cdot \tilde{\delta}_{s_2}{g_2}^{-1})\, d\mu^\mathbf{a}_K(g_1,g_2)\right| & \lesssim M_{\mathbb{H}^n}|G|(g)\int_{B_K({0,1})}\left|F(g\cdot \tilde{\delta_r}{g_1}^{-1})\right| \,dg_1 \\ 
& \lesssim M_{\mathbb{H}^n}|G|(g)\, M^K_{HL}  |F|(g), 
\end{align*}
which is the first claimed inequality. The second one follows simply by interchanging the roles of $F$ and $G$.
\end{proof}

%%%%%%%%%%%%%%%%%%%%%
%%%%%%%%%%%%%%%%%%%%%
%%%%%%%%%%%%%%%%%%%%%
%%%%%%%%%%%%%%%%%%%%%

We now have all the tools to embark on the proof of Theorem \ref{Main-theorem-Koranyi-type-for-degenrate}. 

\begin{proof}[Proof of Theorem \ref{Main-theorem-Koranyi-type-for-degenrate}] 
The proof is divided into three cases, namely, for 
\begin{itemize}
\item $\max \{ \mathbf{a}_1, \mathbf{a}_2 \} \leq 2n+2$; 

\item $\min \{ \mathbf{a}_1, \mathbf{a}_2 \} \leq 2n+2 < \max \{ \mathbf{a}_1, \mathbf{a}_2 \}$; 

\item $\min \{ \mathbf{a}_1, \mathbf{a}_2\} > 2n+2.$
\end{itemize}

\medskip 
\textbf{\underline{Case $1$}:} When $\max \{ \mathbf{a}_1, \mathbf{a}_2 \} \leq 2n+2$. In this case, in order to prove the claimed estimates, we would want to use the bilinear interpolation. Unlike the linear setting, the bilinear interpolation theorem cannot be applied directly to the maximal operator. To overcome this issue, following the approach of \cite{Maximal-estimates-for-the-bilinear-spherical-averages}, we shall first linearize the operator, and on that the bilinear interpolation will be applicable. The linearisation process itself is quite simple, and we include a brief sketch of the arguments for biparametric case for the sake of completeness. But, before doing that, let us see how it helps in the proof.

\medskip
For  measurable functions $\tilde{\tau}_1,\tilde{\tau}_2:\mathbb{H}^n\rightarrow \mathbb{R}_+=\{t\in\mathbb{R}:t>0\}$, let us define 
\[
\mathcal{A}^{\mathbf{a},K}_{\tilde{\tau}_1,\tilde{\tau}_2}(F,G)(g) = \int_{\mathcal{S}_K^\mathbf{a}} F(g \cdot \tilde{\delta}_{\tilde{\tau}_1 (g)}{g_1}^{-1}) \, G(g\cdot \tilde{\delta}_{\tilde{\tau}_2(g)}{g_2}^{-1}) \, d\mu^{\mathbf{a}}_K(g_1,g_2).
\]
Obviously, $\left|\mathcal{A}^{\mathbf{a},K}_{\tilde{\tau}_1,\tilde{\tau}_2}(F,G)\right| \leq \mathcal{M}^{\mathbf{a},K}_{full}(F,G)$, 
and therefore it follows from Lemma \ref{lemma-for-pointwise-control-byknown-maximal-operator} that 
\begin{equation} 
\label{pointwise-bound}
|\mathcal{A}^{\mathbf{a},K}_{\tilde{\tau}_1,\tilde{\tau}_2}(F,G)(g)| \lesssim M^K_{HL} |F|(g) \, M_{\mathbb{H}^n}|G|(g), 
\end{equation}
with the implicit constant independent of the function $\tilde{\tau}$. 

\medskip 
Now, one can invoke the known boundedness properties of the operators $M^K_{HL}$ and $M_{\mathbb{H}^n}$ to conclude that for $n\geq 2,$ $1<p\leq \infty \text{ and } \frac{2n+1}{2n}<q\leq \infty$ with $r$ satisfying $\frac{1}{p}+\frac{1}{q}=\frac{1}{r}$, we have 
\begin{equation} 
\label{estimate-for-the-half-portion}
\|\mathcal{A}^{\mathbf{a},K}_{\tilde{\tau}_1,\tilde{\tau}_2}(F,G)\|_{L^r(\mathbb{H}^n)}
\lesssim \|M^K_{\mathrm{HL}}F\|_{L^p(\mathbb{H}^n)} \|M_{\mathbb{H}^n}G\|_{L^q(\mathbb{H}^n)} \lesssim \|F\|_{L^p(\mathbb{H}^n)} \|G\|_{L^q(\mathbb{H}^n)}, 
\end{equation}
and also for $p=1$ and $\frac{2n+1}{2n}<q\leq \infty$ with $p,q,r$ satisfying the H\"older's condition,
\begin{equation}\label{estimate-for-the-BH-part}
\|\mathcal{A}^{\mathbf{a},K}_{\tilde{\tau}_1,\tilde{\tau}_2}(F,G)\|_{L^{r,\infty}(\mathbb{H}^n)} \lesssim \|F\|_{L^1(\mathbb{H}^n)} \|G\|_{L^q(\mathbb{H}^n)}.
\end{equation}

With the above estimates at hand, one can apply the bilinear interpolation between the estimates obtained in \eqref{estimate-for-the-BH-part} to get strong-type estimate on the line $(C,D)$, that is, for $p=1$ and $\frac{2n+1}{2n}<q< \infty$ with $1+\frac{1}{q}=\frac{1}{r}$, we will have
\[
\|\mathcal{A}^{\mathbf{a},K}_{\tilde{\tau}_1,\tilde{\tau}_2}(F,G)\|_{L^r(\mathbb{H}^n)}
\lesssim
\|F\|_{L^1(\mathbb{H}^n)} \|G\|_{L^q(\mathbb{H}^n)}.
\] 
By symmetry, we can change the role of $F$ and $G$ in the estimates \eqref{estimate-for-the-half-portion} and \eqref{estimate-for-the-BH-part}. Thus, by applying bilinear interpolation we will have \eqref{estimate-for-the-half-portion} for $p,q,r$ satisfying H\"older condition in the shaded region $\mathcal{R}^{\mathbf{a},K}_{full}$ including the open line segment $(A,B)$ and $(C,D)$, as shown in Figure \ref{figure-1-bilinear-degenerate-Koranyi-full}, the closed pentagon $[O,A,B,C,D]$ except the closed line segment $[B,C]$ and the points $A$ and $D$. Note once again that the constants depend only on $p,q$ and $n$, and not on functions $\tilde{\tau}_1,\tilde{\tau}_2$. 

\medskip
Finally, following the above idea, one can also establish the claimed weak-type estimates at points $A$ and $D$ with the help of estimate \eqref{estimate-for-the-BH-part} with $q=\infty$, $r=1$ and  $p=\infty$, $r=1$

\medskip 
Hence, the proof of Theorem \ref{Main-theorem-Koranyi-type-for-degenrate} in Case $1$ will be completed once we can show that given any $F,G\in\mathcal{S}(\mathbb{H}^n)$, there exists  measurable functions  $\tilde{\tau}_1,\tilde{\tau}_2 : \mathbb{H}^n\rightarrow \mathbb{R}_+$ such that 
$$ \mathcal{M}^{\mathbf{a},K}_{full}(F,G)(g) \leq 2 \left|\mathcal{A}_{\tilde{\tau}_1(g),\tilde{\tau}_2(g)}(F,G)(g)\right|.$$
For the same, note first that  for every fixed $g\in\mathbb{H}^n$, the map
$(r_1,r_2)\mapsto\mathcal{A}^{\mathbf{a},K}_{r_1,r_2}(F,G)(g)$ is continuous on
$(0,\infty)^2$. Therefore, using the density of positive rationals $\mathbb{Q}^2_+$ in $(0,\infty)^2$, the supremum over$(0,\infty)^2$ agrees with the supremum over  $\mathbb{Q}^2_+$, that is, 
\begin{align*}
S(g) := \mathcal{M}^{\mathbf{a},K}_{full}(F,G)(g) = \sup_{(r_1,r_2)\in(0,\infty)^2} |\mathcal{A}_{(r_1,r_2)}(F,G)(g)| = \sup_{k\geq1}
\left| \mathcal{A}^{\mathbf{a},K}_{q_{1,k},q_{2,k}} (F,G)(g) \right|. 
\end{align*}
where $\{ (q_{1,k}, q_{2,k}) \}_{k\geq1}$ is a fixed a enumeration of $\mathbb{Q}_+^2$. 

\medskip 
Fix $g \in \mathbb{H}^n$. If $S(g)>0$, then there exists ${I(g)}\in \mathbb{N}$ such that $|\mathcal{A}^{\mathbf{a},K}_{q_{1, I(g)}, q_{2, I(g)}}(F,G)(g)| \geq \frac{S(g)}{2}$, and if $S(g)=0$, then $|\mathcal{A}^{\mathbf{a},K}_{q_{1,k},q_{2,k}}(F,G)(g)| \geq 0$ is true for every $k \geq 1$. This leads us to define functions $\tilde{\tau}_1,
\tilde{\tau}_2.:\mathbb{H}^n\rightarrow\mathbb{R}^+$ by 
\[
\tilde{\tau}_1(g) := q_{1,K(g)} 
\qquad \text{and} \qquad
\tilde{\tau}_2(g) := q_{2,K(g)}, 
\] 
where $K(g)= \min \left\{ I(g) \in\mathbb{N} : |\mathcal {A}^{\mathbf{a},K}_{q_{1,{I(g)}},q_{2,I(g)}}(F,G)(g)| \geq \frac{S(g)}{2} \right\}$. 

\medskip 
It is clear from the construction of $\tilde{\tau}_1$ and $\tilde{\tau}_2$ that $\displaystyle \mathcal{M}^{\mathbf{a},K}_{full}(F,G)(g)\leq 2 \left|\mathcal{A}^{\mathbf{a},K}_{\tilde{\tau}_1(g),\tilde{\tau}_2(g)}(F,G)(g)\right|$. We only need to verify that $\tilde{\tau}_1$, and $\tilde{\tau}_2$ are measurable functions. For the same, note that we can express each of the functions $\tilde{\tau}_i$ as $\tilde{\tau}_i(g)=\sum\limits_{k\geq 1} q_{i,k} \, 1_{\{K(g)=k\}}(g)$, where $i=1,2.$ Now, for each $j\in\mathbb{N}$, define
\[
E_j:=
\left\{
g\in\mathbb{H}^n:
\left|\mathcal{A}^{\mathbf{a},K}_{q_{1,j},q_{2,j}}(F,G)(g)\right|
\geq\frac{S(g)}{2}
\right\}.
\]
Since $S$ and
$\left|\mathcal{A}^{\mathbf{a},K}_{q_{1,j},q_{2,j}}(F,G)\right|$ are measurable functions, each $E_j$ is a measurable set. Moreover, $\{g:K(g)\leq m\} = \bigcup_{j=1}^{m}E_j$, and therefore $\{g: K(g)=k\}$ is a measurable set. Hence, $\tilde{\tau}_1$ and $\tilde{\tau}_2$ are non-negative measurable functions on $\mathbb{H}^n$. 

%%%%%%%%%%%%%%%%%%%%%
%%%%%%%%%%%%%%%%%%%%%
%%%%%%%%%%%%%%%%%%%%%
%%%%%%%%%%%%%%%%%%%%%

\medskip
\textbf{\underline{Case $2$}:} When $\min \{ \mathbf{a}_1, \mathbf{a}_2 \} \leq 2n+2 < \max \{ \mathbf{a}_1, \mathbf{a}_2 \}$. Let us only work for the subcase when $\mathbf{a}_1 \leq 2n+2 < \mathbf{a}_2$. Let us write 
\[
w_\mathbf{a}(g_1) = (1-|g_1|_K^{\mathbf{a}_1})^\frac{1}{\mathbf{a}_2} \quad \text{and} \quad \tilde{w}_\mathbf{a}(g_2) =(1-|g_2|_K^{\mathbf{a}_2} )^\frac{1}{\mathbf{a}_1}. 
\] 
Without loss of any generality, let us assume that the functions $F,G \geq 0.$ In the expression \eqref{equivalent expression-degenrate-Koranyi-2}, with our condition $\mathbf{a}_2>2n+2$, $w_\mathbf{a}$ appears with the negative power $2n+2-\mathbf{a}_2$, and therefore we can not bound that term. To circumvent this issue, we shall work with a dyadic  decomposition of the average away from the unit Kor\'anyi sphere. For $k\geq 1$, we localize the averaging operator to the region where the size of $1-|g_1|_K^{\mathbf{a}_1}$ is comparable to $2^{-k}$. More precisely, if we let
\[
\chi_k^{\mathbf{a}_1}(g_1)
= \chi_{[2^{-k}, \, 2^{1-k}]} \left( 1-|g_1|^{\mathbf{a}_1} \right),
\]
then 
\[
\mathcal{A}_\mathbf{s}^{\mathbf{a},K}(F,G)(g)
    \sim
    \sum_{k\geq 1}\mathcal{A}_\mathbf{s}^{\mathbf{a},k}(F,G)(g),
\]
where
\[
\mathcal{A}_\mathbf{s}^{\mathbf{a},k}(F,G)(g)
=
\int_{|g_1|<1}  F(g\cdot\tilde{\delta}_{s_1} g_1^{-1}) \, \chi_k^{\mathbf{a}_1}(g_1) \, w_\mathbf{a}^{2n+2-\mathbf{a}_2}(g_1) \,A^{K} G(g,s_2w_\mathbf{a}(g_1)) \, dg_1.
\]

We define 
\[
\mathcal{M}^k(F,G)(g) = \sup_{\mathbf{s}\in (0,\infty)^2} \left|\mathcal{A}_\mathbf{s}^{\mathbf{a},k}(F,G)(g)\right|, 
\]
so that 
\begin{align}
\label{case2-decomposition-localised}
\mathcal{M}^{\mathbf{a},K}_{full}(F,G)(g) \lesssim \sum \limits_{k\geq1}\mathcal{M}^k(F,G)(g).
\end{align} 

Now, the localization via the characteristic function for each piece gives the growth estimates in the following way. Since $\chi_k^{\mathbf{a}_1}(g_1) \, w_\mathbf{a}^{2n+2-\mathbf{a}_2}(g_1)\leq 2^{k\frac{(\mathbf{a}_2-2n-2)}{\mathbf{a}_2}}\chi_k^{\mathbf{a}_1}(g_1),$ we obtain 
\begin{equation} 
\label{growth-estimate-K} \|\mathcal{M}^k(F,G)\|_{L^r} \lesssim 2^{k(\frac{\mathbf{a}_2-(2n+2)}{\mathbf{a}_2})} \|F\|_{L^p}\|G\|_{L^q},
\end{equation} 
for $(\frac{1}{p},\frac{1}{q})\in [0,1)\times[0,\frac{2n}{2n+1}).$

\medskip
On the other hand, $(g_1, g_2) \in \mathcal{S}_K^\mathbf{a}$ implies that $1-|g_1|_K^{\mathbf{a}_1} \in [2^{-k}, 2^{1-k}]$ if and only if $|g_2|_K^{\mathbf{a}_2}\in [2^{-k},2^{1-k}].$ Therefore, we get via \eqref{equivalent expression-degenrate-Koranyi} that $\mathcal{A}_\mathbf{s}^{\mathbf{a},k}(F,G)(g)$ is similar to
\[
\int\limits_{{2^{-k} < |g_2|^{\mathbf{a}_2}_K < 2^{1-k}}} G(g \cdot \tilde{\delta}_{s_2} g_2^{-1}) \, A^K F(g,s_1\tilde{w}_\mathbf{a}(g_2)) \, \tilde{{w}}_\mathbf{a}^{(2n+2-\mathbf{a}_1)}(g_2) \, dg_2, 
\]
which in turn imply the decay estimate 
\begin{equation} 
\label{decay-estimate-K}
\|{\mathcal{M}^k}(F,G)\|_{L^r}\lesssim2^{-k(\frac{2n+2}{\mathbf{a}_2})} \|F\|_{L^p}\|G\|_{L^q}, 
\end{equation}
where $(\frac{1}{p},\frac{1}{q})\in [0,\frac{2n}{2n+1})\times[0,1).$ 

\medskip 
With estimates \eqref{growth-estimate-K} and \eqref{decay-estimate-K} at hand, we can make use of the bilinear interpolation (after doing the linearisation as in Case $1$) to conclude that for any $\epsilon > 0$, 
\[
\|\mathcal{M}^k(F,G)\|_{L^r}\lesssim 2^{-\epsilon k}\|F\|_{L^p}\|G\|_{L^q},
\] 
for $(\frac{1}{p},\frac{1}{q})\in\mathcal{R}^\mathbf{a}$, the pentagon with vertices $O, A, B, Z$, and $E$ (see Figure \ref{figure-2-bilinear-degenerate-Koranyi-full}), except the boundary $[A,B],$ $[B,Z]$, and $[Z,E]$, where 
\[
Z = \left( 1 - \frac{1}{2n+1} + \frac{2n+2}{(2n+1)\mathbf{a}_2}, \, 1 - \frac{2n+2}{(2n+1)\mathbf{a}_2} \right), \quad E = \left( 1 - \frac{1}{2n+1} + \frac{2n+2}{(2n+1) \mathbf{a}_2}, \, 0 \right),
\]
and hence, thanks to \eqref{case2-decomposition-localised}, we have the claimed boundedness of $\mathcal{M}^{\mathbf{a},K}_{full}$ in the same region as above. 

%%%%%%%%%%%%%%%%%
%%%%%%%%%%%%%%%%%
%%%%%%%%%%%%%%%%%
%%%%%%%%%%%%%%%%%

\medskip
\textbf{\underline{Case $3$}:} When $\min \{ \mathbf{a}_1, \mathbf{a}_2\} > 2n+2.$ 
With $\mathbf{a}_2 > 2n+2$, we can follow the decomposition as in the analysis of Case $2$ along with the notation from there. We thus obtain the growth estimate \eqref{growth-estimate-K} for all $k\geq1$ here also. As earlier, $(g_1, g_2) \in \mathcal{S}_K^\mathbf{a}$ implies that $1-|g_1|_K^{\mathbf{a}_1} \in [2^{-k}, 2^{1-k}]$ iff $|g_2|_K^{\mathbf{a}_2}\in [2^{-k},2^{1-k}].$ Thus, in this range, we get $\tilde{w}_\mathbf{a}^{(2n+2-\mathbf{a}_1)}(g_2) \leq C$, for some constant $C>0$, uniformly for $k\geq 2.$ And then, we can once again apply bilinear interpolation to conclude that $\sum\limits_{k\geq 2}\mathcal{M}^k$ is bounded in the range same as in Case $2$. So, we are left to analyse $\mathcal{M}^1$. Now, note that for $\mathcal{M}^1$, we already have region using the growth estimate \eqref{growth-estimate-K}, but we can include more via the same type of decomposition on the average as follows. Since 
\[
1-|g_1|_K^{\mathbf{a}_1} \in [1/2,1] \iff 1-|g_2|_K^{\mathbf{a}_2} \in [0,1/2], 
\]
we can split the operator according to the dyadic regions of
\(1-|g_2|^{a_2}\), namely, 
\[
\mathcal{A}_\mathbf{s}^{\mathbf{a},1}(F,G)(g)
=
\sum_{j\geq 2}\tilde{\mathcal{A}}_\mathbf{s}^j(F,G)(g),
\]
where
\[
\tilde{\mathcal{A}}_\mathbf{s}^j (F,G)(g)) = \int
G(g \cdot \tilde{\delta}_{s_2} g_2^{-1}) \, \tilde{w}_{\mathbf{a}}^{(2n+2-\mathbf{a}_1)}(g_2) \, \chi_{[2^{-j}, 2^{-j+1}]} \, (1-|g_2|^{\mathbf{a}_2}) \, A^{K}F(g,s_1\tilde{w}_\mathbf{a}(g_2)) \, dg_2.
\]
Thus, if we define
\[
\widetilde{\mathcal{M}}^j(F,G)(g)
:=
\sup_{\mathbf{s}\in(0,\infty)^2}
\left|
\tilde{\mathcal{A}}_\mathbf{s}^j(F,G)(g), 
\right|
\]
then 
\[
\mathcal{M}^1(F,G)(g)\lesssim\sum\limits_{j\geq2}\widetilde{\mathcal{M}}^j(F,G)(g).
\]
With this, we can repeat the same computation as we did in Case $2$, just that this time we first make use of the expression \eqref{equivalent expression-degenrate-Koranyi} of the averaging operator to get the growth estimate and after that only we invoke the expression \eqref{equivalent expression-degenrate-Koranyi-2}. Now, $w_\mathbf{a}^{(2n+2-\mathbf{a}_2)}(g_1) \leq\tilde{C},$ for some constant $\tilde{C}>0$, uniformly for $j \geq 2$. 
Hence, we obtain the following two estimates for $\widetilde{\mathcal{M}}^j$, for every $j\geq2$: 
\begin{align*}
\| \widetilde{ \mathcal{M}}^j(F,G) \|_{L^r} \lesssim 
\begin{cases} 2^\frac{{j(\mathbf{a}_1-(2n+2))}}{\mathbf{a}_1} \|F\|_{L^p} \, \|G\|_{L^q}, & \text{for } \left( \frac{1}{p}, \frac{1}{q} \right) \in \left[ 0, \frac{2n}{2n+1} \right) \times [0,1), \\[2mm] 
2^{-j\frac{(2n+2)}{\mathbf{a}_1}} \|F\|_{L^p} \, \|G\|_{L^q}, & \text{for } \left( \frac{1}{p},\frac{1}{q} \right) \in [0,1) \times \left[0,\frac{2n}{2n+1} \right), 
\end{cases}
\end{align*}
where $\frac{1}{p}+\frac{1}{q}=\frac{1}{r}.$ 

\medskip 
With the above estimate at hand, we can yet again make use of the bilinear interpolation to conclude that $\sum \limits_{j\geq2} \widetilde{\mathcal{M}}^j (f,g)$ boundedly maps from $L^p \times L^q \rightarrow L^r$ in the pentagon with vertices $O, K, Y, C$, and $D$,  except the boundary $[K,Y]$, $[Y,C]$, and $[C,D]$, where 
\[
Y = \left( 1-\frac{2n+2}{(2n+1) \mathbf{a}_1}, \, 1 - \frac{1}{2n+1} + \frac{2n+2}{(2n+1) \mathbf{a}_1} \right), \quad G = \left( 0, \, 1-\frac{1}{2n+1} + \frac{2n+2}{(2n+1) \mathbf{a}_1} \right).
\]
Hence,   $\mathcal{M}^{\mathbf{a},K}_{full}$ is bounded in the intersection of the two regions where $\mathcal{M}^1$ and $\sum\limits_{k\geq 2}\mathcal{M}^k$ are bounded, and that is exactly the shaded portion in Figure \ref{figure-3-bilinear-degenerate-Koranyi-full}, drawn for the situation when $\frac{1}{\mathbf{a}_1}+\frac{1}{\mathbf{a}_2}>\frac{1}{2n+2}$. Just note that when $\frac{1}{\mathbf{a}_1}+\frac{1}{\mathbf{a}_2}\leq\frac{1}{2n+2}$, the boundedness region of $\mathcal{M}^{\mathbf{a}, K}_{full}$ will not intersect the line segment $BC$ and it will be a rectangle. 

%%%%%%%%%%%%%%%%%%%
%%%%%%%%%%%%%%%%%%%
%%%%%%%%%%%%%%%%%%%
%%%%%%%%%%%%%%%%%%%

\medskip 
Finally, let us see that the estimate \eqref{estimate_for_the_Product_Koranyi_full_degenrate} cannot be true at the points $A$ and $D$ in Figure \ref{figure-1-bilinear-degenerate-Koranyi-full}, that is, for $p=1=r$ and $q=\infty$. Assume on contrary that $\mathcal{M}^{\mathbf{a},K}_{full} : L^1(\mathbb{H}^n) \times L^\infty (\mathbb{H}^n) \rightarrow L^1(\mathbb{H}^n)$ is bounded. Take $G$ to be the constant function $1$. Now, for any $F\in L^1(\mathbb{H}^n)$, identity \eqref{equivalent expression-degenrate-Koranyi-2} gives 
\begin{align*}
\mathcal{M}^{\mathbf{a},K}_{full}(|F|,G)(g) \sim |\Sigma_{K,1}| \sup_{r>0} \int_{B_K(0,1)} |F(g \cdot \tilde{\delta}_r
g_1)| (1-|g_1|_K^{\mathbf{a}_1})^{\frac{2n+2-\mathbf{a}_2}{\mathbf{a}_2}} \, dg_1 \gtrsim M^K_{HL}(|F|)(g).
\end{align*}
Thus, the boundedness assumption for $\mathcal{M}^{\mathbf{a},K}_{full}$ implies a strong-type $(1,1)$ estimate for $M^K_{HL}$, which is not possible. Similarly, it follows that at the point $A$, that is, for $q=1=r$ and $p=\infty$ too, estimate \eqref{estimate_for_the_Product_Koranyi_full_degenrate} cannot hold true. 
\end{proof}

%%%%%%%%%%%%%%%%%%%%%%%
%%%%%%%%%%%%%%%%%%%%%%%
%%%%%%%%%%%%%%%%%%%%%%%
%%%%%%%%%%%%%%%%%%%%%%%

\section{Proof of Theorem \ref{Main-theorem-lacunary-Koranyi-type-for-degenrate}} 
\label{sec:proof-thm-lacunary-max-koranyi}

Recall (from introduction) the H\"older boundedness region \eqref{def:koranyi-holder-bdd-region-single-scale} for the single scale averages $\mathcal{A}^{\mathbf{a},K}_1$, given by 
\begin{align*}
\mathcal{R}^{\mathbf a,K}_1 = \left\{ \left( \frac{1}{p}, \frac{1}{q} \right) \in [0,1]^2 : \|\mathcal{A}^{\mathbf{a},K}_1\|_{L^p \times L^q \rightarrow L^r}<\infty  \text{ with } \frac{1}{p}+\frac{1}{q} = \frac{1}{r} \right\}.
\end{align*}

Note that for any $t>0$, \begin{equation}\label{relation_between_A_1_and_A_t}
\mathcal{A}^{\mathbf{a},K}_t(F,G)=\mathcal{A}^{\mathbf{a},K}_1(F\circ \tilde{\delta}_t,G\circ \tilde{\delta}_t)\circ\tilde{\delta}_\frac{1}{t},
\end{equation}
which implies that $\mathcal{A}^{\mathbf{a},K}_t$ has the same H\"older boundedness region as that of $\mathcal{A}^{\mathbf{a},K}_1$. 

\medskip
We shall closely follow the proof of \cite{ghosh2026bilinearsphericalmaximalfunction} which itself builds upon the strategy used in \cite{Aswin-Hickman}. A fundamental ingredient in the argument is an $L^2(\mathbb{H}^n)\times L^2(\mathbb{H}^n) \to L^1(\mathbb{H}^n)$ decay estimate for suitable localized components of the averaging operator $\mathcal{A}^{\mathbf{a},K}_1$. Then using bilinear interpolation we will have decay estimates for the localized pieces in the $int({R}^{\mathbf a,K}_1)$. 

\medskip 
Recall now the slicing argument as in \eqref{equivalent expression-degenrate-Koranyi-2} for the degenerate product Kor\'anyi case to rewrite the bilinear averages $\mathcal{A}^{\mathbf{a},K}_1 (F,G)$ as follows: 
\begin{align*}
& \mathcal{A}^{\mathbf{a},K}_1
(F,G)(g) \\ 
&= \frac{1}{\mathbf{a}_2} \int_{{|g_1|_K<1}} F(g \cdot g_1^{-1}) \, A^K G(g, (1-|g_1|_K^{\mathbf{a}_1})^\frac{1}{\mathbf{a}_2})(1-|g_1|_K^{\mathbf{a}_1})^\frac{(2n+2-\mathbf{a}_2)}{\mathbf{a}_2}\,dg_1 \\ 
&=\frac{1}{\mathbf{a}_2}\int_0^1 s^{2n+1}(1-s^{\mathbf{a}_1})^\frac{2n+2-\mathbf{a}_2}{\mathbf{a}_2} A^K_s(F)(g)\,{A}^K_{(1-s^{\mathbf{a}_1})^\frac{1}{\mathbf{a}_2}} (G)(g) \, ds, 
\end{align*}
where the last step follows from the polar decomposition \eqref{expresstion-of-the-surface-measure-we-need-for-the-polar-decomposition}, and here we are using the notation $A^K(F)(g,s) := A^K_s(F)(g)$. 
Now, it follows from the Minkowski's and H\"older's inequalities that 
\begin{align}\label{expression-of-sliced-averaging-operator}
   \| \mathcal{A}^{\mathbf{a},K}_1(F,G)\|_{L^1(\mathbb{H}^n)} \leq  \int_0^1 s^{2n+1}(1-s^{\mathbf{a}_1})^\frac{2n+2-\mathbf{a}_2}{\mathbf{a}_2}\, \|A^K_s(F)\|_{L^2(\mathbb{H}^n)} \|A^K_{(1-s^{\mathbf{a}_1})^\frac{1}{\mathbf{a}_2}}G\|_{L^2(\mathbb{H}^n)}  \, ds.
\end{align}

\medskip 
We shall apply the above estimate on localised pieces of $F$ and $G$. For the same, take and fix a function $\phi\in C_c^\infty(\mathbb{H}^n)$ such that $\int_{\mathbb{H}^n} \phi(g) \, dg = 1$. For each $r>0$, let $\phi_r$ denote the corresponding dilate of
$\phi$, that is, $\phi_r=r^{-\mathcal{Q}}(\phi\circ \tilde{\delta}_{r^{-1}})$ where $\mathcal{Q}=2n+2$ is the Homogeneous dimension of $\mathbb{H}^n$ so that $\{\phi_r\}_{r>0}$ forms an approximate identity. For every $k\in\mathbb{Z}$, if we define $\psi_r(g) = \phi_r(g) - \phi_{2r}(g)$, then each $\psi_r$ satisfies the cancellation property, namely, $\int_{\mathbb{H}^n} \psi_r = 0$, and the collection $\{\psi_{2^k}\}_{k\in\mathbb{Z}}$ gives rise to a
Littlewood--Paley decomposition, that is, for any $F\in C_c^\infty(\mathbb{H}^n)$, one gets 
\begin{equation} 
\label{littlewood_paley_1}
F=\sum_{k\in\mathbb{Z}} F*\psi_{2^k}. 
\end{equation} 
Using \eqref{littlewood_paley_1}, we decompose $F$ and $G$ as $F=\sum\limits_{j\in \mathbb{Z}} F*\psi_{2^j}$ and $G=\sum\limits_{k\in \mathbb{Z}} G*\psi_{2^k}$, which gives
\[
\mathcal{A}^{\mathbf{a},K}_1(F,G)(g)=\sum\limits_{j,k\in \mathbb{Z}} \mathcal{A}_1^{j,k}(F,G)(g),
\]
where 
\[
\mathcal{A}_1^{j,k}(F,G)(g)=\mathcal{A}^{\mathbf{a},K}_1(F*\psi_{2^j},G*\psi_{2^k})(g).
\]

As mentioned earlier, following is the main $L^2(\mathbb{H}^n)\times L^2(\mathbb{H}^n) \to L^1(\mathbb{H}^n)$ estimate that we shall shortly prove. There exists $0 < \delta < 2n+2$ such that 
\begin{align} 
\label{2-2-1-decay_estimate_for_frequency_pieces} 
\| \mathcal{A}^{j,k}_1(F,G)\|_{L^1(\mathbb{H}^n)} \lesssim 2^{(k+j)\delta} \|F\|_{L^2(\mathbb{H}^n)}\|G\|_{L^2(\mathbb{H}^n)}, \end{align} 
for every $j, k \in \mathbb{Z}$. 

\medskip In fact, once we have \eqref{2-2-1-decay_estimate_for_frequency_pieces}, we can make use of the fact that $\mathcal{A}_1^{j,k}$ is bounded in ${R}^{\mathbf a,K}_1$ region, and then the bilinear interpolation implies that for any $(1/p, 1/q) \in int ({R}^{\mathbf a,K}_1)$, there exist $\delta =: \delta (p,q,n) >0$ so that
\begin{align} 
\label{decay_estimate_for_frequency_pieces}
\| \mathcal{A}^{j,k}_1(F,G)\|_{L^r(\mathbb{H}^n)} \lesssim 2^{(k+j)\delta} \|F\|_{L^p(\mathbb{H}^n)}\|G\|_{L^q(\mathbb{H}^n)}, 
\end{align}
for every $j,k \in \mathbb{Z}$ and with $\frac{1}{p}+\frac{1}{q}=\frac{1}{r}$. 

\medskip 
So, let us embark on the proof of \eqref{2-2-1-decay_estimate_for_frequency_pieces}. For the same, note first that we can make use of the estimate \eqref{expression-of-sliced-averaging-operator} for the localised pieces $\mathcal{A}_1^{j,k}$
to obtain
\begin{align} \label{estimate-of-A^j,k(f,g)-for-2,2,1-case}
& \| \mathcal{A}^{j,k}_1(F,G) \|_{L^1(\mathbb{H}^n)} \\ 
& \leq \int_0^1 s^{2n+1}(1-s^{\mathbf{a}_1})^\frac{2n+2-\mathbf{a}_2}{\mathbf{a}_2} \, \|A^K_s(F*\psi_{2^j})\|_{L^2(\mathbb{H}^n)} \|A^K_{(1-s^{\mathbf{a}_1})^\frac{1}{\mathbf{a}_2}}(G*\psi_{2^k})\|_{L^2(\mathbb{H}^n)} \, ds, \notag
\end{align}
for every $j,k \in \mathbb{Z}$. We shall now prove the following $L^2-L^2$ estimate: there exists $\gamma >0$ such that for large $N \in \mathbb{N}$: 
\begin{align} 
\label{L^2-L^2-estimate}
\|A^K_s(F*\psi_{2^j})\|_{L^2(\mathbb{H}^n)} \lesssim_N C \left( 2^j/s \right)^{\gamma/2^N} \|F\|_{L^2(\mathbb{H}^n)}, 
\end{align}
but before proving it, let us see that using it in \eqref{estimate-of-A^j,k(f,g)-for-2,2,1-case} in fact implies \eqref{2-2-1-decay_estimate_for_frequency_pieces}. 

\medskip 
Take $\delta=\frac{\gamma}{2^N}$ for a large enough $N\geq1,$ so that $0<\delta<2n+2$. Now, using \eqref{L^2-L^2-estimate} in \eqref{estimate-of-A^j,k(f,g)-for-2,2,1-case}, we apply Young's inequality to obtain 
\begin{align*}
\| \mathcal{A}^{j,k}_1(F,G)\|_{L^1(\mathbb{H}^n)} &\lesssim  \|F\|_{L^2(\mathbb{H}^n)}\|G\|_{L^2(\mathbb{H}^n)}\int_0^1 s^{2n+1}(1-s^{\mathbf{a}_1})^\frac{2n+2-\mathbf{a}_2}{\mathbf{a}_2}\,  \left(\frac{2^j}{s}\right)^\delta \left(\frac{2^j}{(1-s^{\mathbf{a}_1})^\frac{1}{\mathbf{a}_2}}\right)^\delta \, ds\\
&\lesssim 2^{(k+j)\delta} \|F\|_{L^2(\mathbb{H}^n)}\|G\|_{L^2(\mathbb{H}^n)}, 
\end{align*}
which is exactly the claimed estimate \eqref{2-2-1-decay_estimate_for_frequency_pieces}. 

\medskip 
Let us now prove estimate \eqref{L^2-L^2-estimate}. For the same, observe first that 
\[
A^K_s(F*\psi_{2^j})=A_1^K(F\circ \tilde{\delta}_s*\psi_\frac{2^j}{s})\circ \tilde{\delta}_\frac{1}{s}, 
\]
and therefore, if we define $TF=F*\psi_\frac{2^j}{s}*\sigma^{2n}_K$, then the estimate \eqref{L^2-L^2-estimate} is equivalent to 
proving that for large $N\in \mathbb{N}$, 
\begin{align} 
\label{L2-L2-estimate}
\|TF\|_{L^2(\mathbb{H}^n)} \lesssim_N \left( 2^j/s \right)^{\gamma/2^N} \|F\|_{L^2(\mathbb{H}^n)},
\end{align}
uniformly in $j \in \mathbb{Z}$. 

\medskip 
In order to prove \eqref{L2-L2-estimate}, let us begin with recalling the following from \cite{Aswin-Hickman}. Let $\mu$ be a finite Borel measure on $\mathbb{H}^n$. The iterated convolutions of $\mu$ are defined as follows. Set $\mu^{(0)}=\mu$, and for $N\geq 1$, define 
\[
\mu^{(N)} = \begin{cases} \mu^{(N-1)} * \widetilde{\mu}, & N \text{ odd},\\[2mm] \mu^{(N-1)} * \mu, & N\geq 2 \text{ even}, \end{cases} 
\] where $\widetilde{\mu}$ denotes the reflected measure. We say that $\nu$ satisfies the \emph{curvature assumption} (CA) if there exists an $N\in\mathbb{N}$ such that the $N$-fold convolution $\nu^{(N)}$ is absolutely continuous with respect to the Haar measure $dg$ on $\mathbb{H}^n$, that is, 
\[ 
d\nu^{(N)}(g) = h(g) \, dg, 
\]
and moreover there exist constants $\gamma>0$ and $C_\nu>0$ such that 
\begin{equation} 
\label{expression-of-Holder-L^!-inequality}
\int_{\mathbb{H}^n} \left( \left|h(gg'^{-1})-h(g)\right| + \left|h(g'^{-1}g)-h(g)\right| \right) dg \leq C_\nu |g'|^\gamma 
\end{equation}
for every $g'\in\mathbb{H}^n$.

\medskip 
In particular, if we take
$\mu = \sigma^{2n}_K$, the surface measure on the Kor\'anyi sphere, it is already known from \cite[Lemma 3.2]{Aswin-Hickman} that $\mu$ satisfies the curvature assumption (CA), so there exists $M \geq 1$ such that for every $N \geq M$,  we have that $\mu^{(N-1)} \ll dx$ and it satisfies the H\"older $L^1$-type inequality \eqref{expression-of-Holder-L^!-inequality}. 

\medskip 
In order to prove \eqref{L2-L2-estimate}, it suffices to work with the operator $T^*T$. In general, given any $N \geq 1$, let us consider the operator $A^{(N)}$ defined by $A^{(N)}F := F*\psi_\frac{2^j}{s}*\mu^{(N-1)}$. Clearly, $T = A^{(1)}$. Now, note that in view of the fact that $\widetilde{\sigma_K^{2n}} = \sigma_K^{2n}$, we get ${A^{(N)}}^*F = F*\mu^{(N-1)}*\widetilde{\psi_\frac{2^j}{s}}$, where $\widetilde{\psi_\frac{2^j}{s}}(g) = \overline{\psi_\frac{2^j}{s}(g^{-1})}$, and hence
\begin{align}
\label{identity-T*T-for-N}
{A^{(N)}}^*{A^{(N)}} F = F * \psi_\frac{2^j}{s} * \mu^{(N-1)} * \mu^{(N-1)} * \widetilde{\psi_\frac{2^j}{s}}.
\end{align}
Now, if we apply Young's inequality on \eqref{identity-T*T-for-N}, for each $N \geq 2$, we obtain  
\begin{align*}
\|(A^{(N)})^*A^{(N)}F\|_{L^2(\mathbb{H}^n)}&\leq \|A^{(N+1)}F\|_{L^2(\mathbb{H}^n)} \, \|\psi_\frac{2^j}{s}*\mu^{(N-2)}\|_{L^1(\mathbb{H}^n)}, 
\end{align*}
whereas for $N=1$, we will have 
\begin{align*}
\|(A^{(1)})^*A^{(1)}F\|_{L^2(\mathbb{H}^n)} \lesssim \|A^{(2)} F \|_{L^2(\mathbb{H}^n)} \, \|\psi_\frac{2^j}{s} \|_{L^1(\mathbb{H}^n)}.
\end{align*}
Now, using the above recurrence relation, we get that for every $N \geq 1$, we have 
\begin{align}\label{recurrence-relation}
\|A^{(1)}F\|_{L^2(\mathbb{H}^n)}\lesssim_N  \|\psi_\frac{2^j}{s}*\mu^{(N-1)}\|_{L^1(\mathbb{H}^n)}^{1/2^N}\|F\|_{L^2(\mathbb{H}^n)}, 
\end{align}
where we have made use of the basic fact that $\|A^{(N)}\|_{L^2 \rightarrow L^2} \leq C_N,$ where the constant $C_N$ is independent of $j$ and $s$. 

\medskip 
With the above estimate at hand, we can complete the proof of estimate \eqref{L2-L2-estimate} with the help of the H\"older $L^1$-type inequality \eqref{expression-of-Holder-L^!-inequality} as follows. Let $h^{(N)}$ denote the Radon-Nikodym derivative of the measure $\mu^{(N)}$. Then, using the cancellation of $\psi$, we have 
\begin{align*} 
\|\psi_\frac{2^j}{s}*\mu^{(N-1)}\|_{L^1(\mathbb{H}^n)} & = \int_{\mathbb{H}^n}\left|\int_{\mathbb{H}^n}\left(h^{(N-1)}(g'^{-1}\cdot g)-h^{(N-1)}(g) \right) \psi_{\frac{2^j}{s}}(g')\,dg' \right|\,dg \\ 
&\lesssim_N \int_{\mathbb{H}^n} |g'|^\gamma \,|\psi_{\frac{2^j}{s}}|(g') \, dg' \lesssim_N \left( 2^j/s \right)^{\gamma}. 
\end{align*}

Thanks to estimate \eqref{decay_estimate_for_frequency_pieces}, we are now ready to complete the proof of Theorem \ref{Main-theorem-lacunary-Koranyi-type-for-degenrate}. 

\begin{proof}[Proof of Theorem \ref{Main-theorem-lacunary-Koranyi-type-for-degenrate}] Take and fix $(\frac{1}{p},\frac{1}{q}) \in int ({R}^{\mathbf a,K}_1)$. Note that by shifting the indices by any $l \in \mathbb{Z}$, we deduce from \eqref{littlewood_paley_1} that 
\begin{equation}\label{littlewood_paley_2}
F=F * \phi_{2^\ell}+\sum_{j<0} F * \psi_{2^{j+\ell}}.
\end{equation}
Identity \eqref{littlewood_paley_2} helps in decomposing the single lacunary average $\mathcal{A}^{\mathbf{a},K}_{2^l}$ into low-high, high-low, low-low, and high-high frequency pieces as follows:
\begin{align} 
\label{expression-of-frequency-decomposition}
\mathcal{A}^{\mathbf{a},K}_{2^l}(F,G)(g)=&\mathcal{A}^{\mathbf{a},K}_{2^l}(F*\phi_{2^l},G)(g)+\mathcal{A}^{\mathbf{a},K}_{2^l}(F,G*\phi_{2^l})(g)-\mathcal{A}^{\mathbf{a},K}_{2^l}(F*\phi_{2^l},G*\phi_{2^l})(g)\\&+\sum\limits_{j,k<0}\mathcal{A}^{\mathbf{a},K}_{2^l}(F*\psi_{2^{j+l}},G*\psi_{2^{k+l}})(g).\notag
\end{align}

We shall work with each term in \eqref{expression-of-frequency-decomposition} separately. To begin with, since $\mathcal{A}^{\mathbf{a},K}_{2^l}(F * \phi_{2^l},G)=\mathcal{A}^{\mathbf{a},K}_{1}(F \circ  \tilde{\delta}_{2^l}* \phi,G \circ  \tilde{\delta}_{2^l})\circ \tilde{\delta}_{\frac{1}{2^l}}$, it suffices to work with $\mathcal{A}^{\mathbf{a},K}_1(F*\phi,G)$. Note that since 
$\phi\in C_c^\infty(\mathbb H^n),$ it is straightforward to verify that 
\begin{align} 
\label{low-high-maximal}
|(F*\phi)(g\cdot \tilde{g}^{-1})|\lesssim_\phi M^K_{HL} |F|(g), 
\end{align}
for all $\tilde{g}\in \mathbb{H}^n$ with $|\tilde{g}|_K<1$. 

\medskip 
Now, unlike the proof of \cite{ghosh2026bilinearsphericalmaximalfunction}, we estimate $\mathcal{A}^{\mathbf{a},K}_1(F*\phi,G)$ with the help of \eqref{low-high-maximal} and the slicing argument \eqref{equivalent expression-degenrate-Koranyi} as follows: 
\begin{align*}
|\mathcal{A}^{\mathbf{a},K}_1(F*\phi,G)(g)|&\leq \int_{\mathcal{S}_K^\mathbf{a}} |F*\phi(g\cdot g_1^{-1})|\,|G(g\cdot g_2^{-1})| \, d\mu^{\mathbf{a}}_K(g_1,g_2)\\
&\lesssim  M^K_{HL} |F|(g) \int_{B_K(0,1)}|G(g\cdot g_2^{-1})| \, (1-|g_2|_K^{\mathbf{a}_2})^\frac{(2n+2-\mathbf{a}_1)}{\mathbf{a}_1}\,dg_2 \\
&\lesssim M^K_{HL} |F|(g) \, M^K_{HL} |G|(g),
\end{align*}
and it is this step where we require the assumption $\mathbf{a}_1\leq2n+2.$ 

\medskip 
Altogether, the above estimate implies that 
\begin{equation} 
\label{estimate_for_low_high_frequency_piece}
\mathcal{M}^{LH}_{lac}(F,G)(g) := \sup_{l \in \mathbb{Z}} |\mathcal{A}^{\mathbf{a},K}_{2^l} (F*\phi,G)(g)| \lesssim M^K_{HL} |F|(g) \, M^K_{HL} |G|(g).
\end{equation}
Arguing similarly, we can also control the high-low, and low-low frequency part by the product $M^K_{HL} |F|(g) \, M^K_{HL} |G|(g)$, provided  $\mathbf{a}_1, \mathbf{a}_2 \leq 2n+2.$  Summarising, for $n\geq2$, 
\[
\|\mathcal{M}^{LH}_{lac}+\mathcal{M}^{HL}_{lac}+\mathcal{M}^{LL}_{lac}\|_{L^p\times L^q \rightarrow L^r}<\infty,
\] for any $1< p,q \leq \infty$ with $\frac{1}{p}+\frac{1}{q}=\frac{1}{r}$.

\medskip
Finally, for the high-high frequency pieces in \eqref{expression-of-frequency-decomposition}, one can essentially repeat the proof of \cite{ghosh2026bilinearsphericalmaximalfunction}. In doing so, whenever there is a further decomposition involving low-high or high-low frequency terms, we can make use of the idea as above to dominate by the Hardy--Littlewood maximal function $M^K_{HL}$ instead of the Hopf's maximal function associated with Birkhoff averages. This completes the proof of Theorem \ref{Main-theorem-lacunary-Koranyi-type-for-degenrate}. 
\end{proof}

%%%%%%%%%%%%%%
%%%%%%%%%%%%%%
%%%%%%%%%%%%%%
%%%%%%%%%%%%%%
%%%%%%%%%%%%%%

\section{Bilinear maximal functions for Nevo--Thangavelu type degenerate surfaces} 
\label{sec:NT-maximal-degenerate}

This section concerns with the proofs of Theorems \ref{thm:full-max-NT-degenerate}, \ref{thm:lac-max-NT-degenerate} and \ref{theorem-improving-single-scales}, which deals with the boundedness of the full and lacunary maximal functions corresponding to  the following degenerate surfaces associated with the Nevo--Thangavelu type averaging structure. Given $\mathbf{a}=(\mathbf{a_1,a_2})$ with $\mathbf{a_1,a_2} > 0$, let 
\begin{align*} 
S_{NT}^{\mathbf{a}} := \{(g_1, g_2) = (z_1,0,z_2,0) \in \mathbb{H}^n \times \mathbb{H}^n : |z_1|^{\mathbf{a_1}} +|z_2|^{\mathbf{a_2}}=1 \}. 
\end{align*}
Clearly, $S_{NT}^\mathbf{a}$ is a $4n-1$ dimensional surface in $\mathbb{R}^{4n+2}$. 

\medskip 
Consider the natural embedding
\[
\iota:\mathbb{R}^{4n}\longrightarrow\mathbb{R}^{4n+2}
\]
defined by $\iota(z_1,z_2)=(z_1,0,z_2,0),$
which preserves the Euclidean metric. Clearly, \(S
_{NT}^{\mathbf a} = \iota(\Sigma^{\mathbf a})\), where $\Sigma^{\mathbf a}$ is the $4n-1$ dimensional hypersurface in \(\mathbb{R}^{4n}\) given by 
\[
\Sigma^{\mathbf a} = \left\{ (z_1,z_2) \in \mathbb{R}^{4n}: |z_1|^{\mathbf a_1} + |z_2|^{\mathbf a_2} = 1
\right\}, 
\]
equipped with its natural \(4n-1\) dimensional Hausdorff measure \(d\mu^{\mathbf a}\).

\medskip
Let us consider the measure $d\mu_{NT}^{\mathbf a}$ defined on \(S_{NT}^{\mathbf a}\) by the pushforward of \(d\mu^{\mathbf a}\) under \(\iota\), that is, $d\mu_{NT}^{\mathbf a}:=\iota_*d\mu^{\mathbf a}.$ Since $\iota$ is an isometry, this pushforward measure is naturally the surface measure on $S_{NT}^{\mathbf a}$. With this, the biparametric bilinear averaging operator  $\mathcal{A}_\mathbf{s}^{\mathbf{a},NT}$ over the surface $S_{NT}^{\mathbf a}$ can be defined as in \eqref{express-averaging-operator-nevotype}, and we write it again for the sake of convenience: 
\begin{align} \label{expression_of_averaging_operator_at_degenerate_surface_Nevo_Thangavelu}
\mathcal{A}_\mathbf{s}^{\mathbf{a},NT}(F,G)(g) &= \int\limits_{S_{NT}^{\mathbf a}} F(g \cdot \tilde{\delta}_{s_1} g_1^{-1}) \, G(g \cdot \tilde{\delta}_{s_2} g_2^{-1}) \, d\mu_{NT}^\mathbf{a} (g_1, g_2) \notag\\
&=\int\limits_{\Sigma^\mathbf{a}} F(g\cdot\tilde{\delta}_{s_1}(z_1,0)^{-1}) \, G(g\cdot\tilde{\delta}_{s_2}(z_2,0)^{-1}) \, d\mu^\mathbf{a}(z_1,z_2).
\end{align}

In the next subsection, we discuss the proofs of Theorems \ref{thm:full-max-NT-degenerate} and \ref{thm:lac-max-NT-degenerate}, which follow from more or less the same ideas as those in the proofs of Theorems \ref{Main-theorem-Koranyi-type-for-degenrate} and \ref{Main-theorem-lacunary-Koranyi-type-for-degenrate} respectively. We shall mainly explain how to deploy an appropriate slicing argument. We shall discuss the idea behind the proof of Theorem \ref{theorem-improving-single-scales} in Subsection \ref{subsec:thm-single-scale-NT-degenrate}. 

%%%%%%%%%%%%%%%%%%%%%%%%
%%%%%%%%%%%%%%%%%%%%%%%%
%%%%%%%%%%%%%%%%%%%%%%%%
%%%%%%%%%%%%%%%%%%%%%%%%

\subsection{Proofs of Theorems \ref{thm:full-max-NT-degenerate} and \ref{thm:lac-max-NT-degenerate}} 
\label{subsec:thms-full-and-lacunary-NT-degenrate} 

It suffices to work with functions $F, G \geq 0.$ Similar to Section \ref{sec:proof-thm-full-max-koranyi}, here also we shall resort to an appropriate slicing argument. We shall work with $\mathbf{a}_1, \mathbf{a}_2 \geq 1$. As in Subsection \ref{exp-surface-measure-on-the-degenerate-product-Koranyi-sphere}, here also we exploit the formula \eqref{coarea-formula-statement}. For the same, consider the map  $\Phi_{NT}^\mathbf{a}(z_1,z_2)=|z_1|^{\mathbf{a}_1}+|z_2|^{\mathbf{a}_2}-1.$ Clearly, ${\Phi_{NT}^\mathbf{a}}^{-1}\{0\} = \Sigma^\mathbf{a}.$ Moreover, $|\nabla \Phi^\mathbf{a}_{NT}(z_1,z_2)| = \sqrt{ \mathbf{a}_1^2 |z_1|^{2(\mathbf{a}_1-1)} + \mathbf{a}_2^2 |z_2|^{2(\mathbf{a}_2-1)}}$ is non vanishing a.e. on the surface $\Sigma^\mathbf{a}$, and it is easy to verify that there exist positive constants $C$ and $D$ such that 
\[
C \leq |\nabla \Phi^\mathbf{a}_{NT}(z_1,z_2)| \leq D
\]
for a.e. $(z_1,z_2) \in \Sigma^{\mathbf a}.$ Thus, it follows from \eqref{expression_of_averaging_operator_at_degenerate_surface_Nevo_Thangavelu} that  
\begin{align*}
\mathcal{A}_\mathbf{s}^{\mathbf{a}, NT} (F,G)(g) 
& \sim \int \limits_{\Sigma^\mathbf{a}} F(g\cdot\tilde{\delta}_{s_1}(z_1,0)^{-1}) \, G(g\cdot\tilde{\delta}_{s_2}(z_2,0)^{-1}) \, \frac{d\mu^\mathbf{a}(z_1,z_2)}{|\nabla \Phi^\mathbf{a}_{NT}(z_1,z_2)|} \\
&= \int_{\mathbb{R}^{4n}} F(g\cdot\tilde{\delta}_{s_1}(z_1,0)^{-1}) \, G(g\cdot\tilde{\delta}_{s_2}(z_2,0)^{-1}) \, \delta((\Phi_{NT}^\mathbf{a}(z_1,z_2)) \, dz_1 \, dz_2. 
\end{align*}

\medskip 
In the following, let us write $\Omega_{z_1} = \Phi_{NT}^\mathbf{a}(z_1,\cdot)^{-1}(0)$, the standard sphere of radius $(1-|z_1|^{\mathbf{a}_1})^\frac{1}{{\mathbf{a}_2}}$ in $\mathbb{R}^{2n}$. For fixed $z_1$, if we apply the formula \eqref{coarea-formula-statement} again for $\Phi_{NT}^\mathbf{a}(z_1,\cdot)(z_2)=|z_1|^{\mathbf{a}_1}+|z_2|^{\mathbf{a}_2}-1$, then we obtain 
\begin{align} \label{equivalent-expresion-1-nevo-type}
\mathcal{A}_\mathbf{s}^{\mathbf{a},NT}(F,G)(g) 
&\sim \int_{|z_1|<1} F(g\cdot\tilde{\delta}_{s_1}(z_1,0)^{-1})\left\{\int_{\Omega_{z_1}} G(g\cdot\tilde{\delta}_{s_2}(z_2,0)^{-1})\, \frac{d\sigma_{(1-|z_1|^{\mathbf{a}_1})^\frac{1}{\mathbf{a}_2}}(z_2)}{|\nabla_{z_2}(\Phi_{NT}^\mathbf{a}(z_1,\cdot) (z_2)|} \right\} \, dz_1 \notag \\ 
&\sim \int_{|z_1|<1} F(g\cdot\tilde{\delta}_{s_1}(z_1,0)^{-1}) \, A^{NT} G(g, s_2 w^{NT}_\mathbf{a}(z_1)) \, \left( w^{NT}_\mathbf{a}(z_1) \right)^{(2n-\mathbf{a}_2)} \,dz_1,
\end{align}
\medskip
where $w^{NT}_\mathbf{a}(z_1) = (1-|z_1|^{\mathbf{a}_1})^\frac{1}{{\mathbf{a}_2}}$ and $A^{NT}G(g, r)$ are nothing but the spherical averages define in \eqref{def:linear-NT-averages}, namely, 
$$ A^{NT} G(g, r) := A^{NT}_r G(g) = \int_{S^{2n-1}} G\left(g \cdot \tilde{\delta}_r(z,0)^{-1} \right)\, d\sigma(z).$$

\medskip 
Now, similar to the product Kor\'anyi case, depending on whether $\mathbf{a}_1$ or $\mathbf{a}_2$ is less than or greater than $2n$, the following three cases arise: 
\begin{itemize}
\item $\max \{\mathbf{a}_1,\mathbf{a}_2\} \leq 2n$; 

\item $\min \{ \mathbf{a}_1, \mathbf{a}_2\} \leq 2n < \max \{ \mathbf{a}_1, \mathbf{a}_2\}$; 

\item $2n < \min \{ \mathbf{a}_1, \mathbf{a}_2\}.$
\end{itemize}

\medskip 
In Case $1$, that is, when $\max \{\mathbf{a}_1,\mathbf{a}_2\} \leq 2n$, following the idea of \cite{ghosh2026bilinearsphericalmaximalfunction}, we can dominate the full maximal operator $ \mathcal{M}^{\mathbf{a},NT}_{full}$ by a product of the Nevo--Thangavelu spherical maximal function $M^{NT}_{\mathbb{H}^n}$ and Hopf's ergodic maximal function $\Lambda$ on $\mathbb{H}^n$ associated with the uniform averages $\Lambda_r$ of Nevo--Thangavelu spherical averages $A_s^{NT}$. More precisely, let 
\[
\Lambda_r F(g)
:=
\frac{1}{r}\int_0^r A_s^{NT} F(g) \, ds,
\]
and define the corresponding Hopf's ergodic maximal function $\Lambda$ by 
\begin{align}
\label{def:Hopfs-ergodic-maximal}
\Lambda F(g)
:=
\sup_{r>0}\left|\Lambda_r F(g)\right|
=
\sup_{r>0}
\left|
\frac{1}{r}\int_0^r A_s^{NT}F(g)\,ds
\right|. 
\end{align}
It is well known from Hopf's maximal ergodic theorem that $\Lambda$ is of strong-type $(p,p)$ for every $1<p\leq\infty$, (see \cite[Proposition 3.6.1]{Harmonic-analysis-on-the-Heisenberg-Group-Thangavelu}).

\medskip
One can essentially repeat the proof of Lemma \ref{lemma-for-pointwise-control-byknown-maximal-operator} to show that 
\begin{align*} 
\mathcal{M}^{\mathbf{a},NT}_{full}(F,G)(g)\lesssim \Lambda |F|(g) \, M^{NT}_{\mathbb{H}^n} |G|(g) 
\quad \text{and} \quad 
\mathcal{M}^{\mathbf{a},NT}_{full}(F,G)(g) \lesssim M^{NT}_{\mathbb{H}^n} |F|(g) \, \Lambda |G|(g), 
\end{align*}
and with that one can complete the proof of the Theorem \ref{thm:full-max-NT-degenerate} in all three cases by following the proof of Theorem \ref{Main-theorem-Koranyi-type-for-degenrate}, with the only care that one should use Hopf's ergodic maximal function $\Lambda$ instead of the Hardy--Littlewood maximal function $M^K_{HL}$ and the Nevo--Thangavelu spherical maximal function $M^{NT}_{\mathbb{H}^n}$ in place of the Kor\'anyi spherical maximal function $M_{\mathbb{H}^n}$. Unlike Theorem \ref{Main-theorem-Koranyi-type-for-degenrate}, we don't know anything about the boundary, namely, the weak-type estimates at the extreme points. This is because we do not know whether the Hopf's ergodic maximal function $\Lambda$ is weak-type $(1,1)$ bounded or not. 

\medskip
Finally, let us discuss the lacunary case of Theorem \ref{thm:lac-max-NT-degenerate}. Recall from \eqref{def:NT-holder-bdd-region-single-scale} the H\"older boundedness region of $\mathcal{A}^{\mathbf{a},NT}_1$, that is, 
\[
\mathcal{R}_1^{\mathbf{a},NT}=\left\{\left(\frac{1}{p},\frac{1}{q} \right) :1\leq p,q\leq\infty,\frac{1}{p}+\frac{1}{q}=\frac{1}{r} \text{ and } \|\mathcal{A}^{\mathbf{a},NT}_1\|_{L^p \times L^q \rightarrow L^r}<\infty\right\}.
\]
Following the proof of Theorem \ref {Main-theorem-lacunary-Koranyi-type-for-degenrate} and using the slicing argument as in \eqref{equivalent-expresion-1-nevo-type}, we already have 
\[
 \mathcal{A}^{\mathbf{a},NT}_1
(F,G)(g)\sim \int_{|z_1|<1} F(g\cdot(z_1,0)^{-1}) \, A^{NT} G(g, w_\mathbf{a}(z_1)) \, \left( w^{NT}_\mathbf{a}(z_1)\right)^{2n-\mathbf{a}_2} \, dz_1,
\] 
where $w_\mathbf{a}(z_1)=(1-|z_1|^{\mathbf{a}_1})^\frac{1}{{\mathbf{a}_2}}$. Applying the polar decomposition on $z_1$-variable, we obtain
\[
\mathcal{A}^{\mathbf{a},NT}_1
(F,G)(g)\sim \int_0^1 s^{2n-1} \, (1-s^{\mathbf{a}_1})^\frac{2n-\mathbf{a}_2}{\mathbf{a}_2} \, A^{NT} F (g,s) \, A^{NT} G(g,{(1-s^{\mathbf{a}_1})^\frac{1}{\mathbf{a}_2}}) \, ds.
\]

Now, we can make use of the Littlewood-Paley decomposition of $F$ and $G$ as in \eqref{littlewood_paley_1} and repeat the proof of Theorem \ref{Main_Theorem_Lacunary}. Thus, replacing $F$ by $F*\psi_{2^j}$ and $G$ by $G*\psi_{2^k}$ and then taking $L^1$-norm, right hand side of the above expression will be dominated by
\begin{align*}
\int_0^1 s^{2n-1} \, (1-s^{\mathbf{a}_1})^\frac{2n-\mathbf{a}_2}{\mathbf{a}_2} \, \|A^{NT}(F*\psi_{2^j})(\cdot,s)\|_{L^2} \, \|A^{NT}(G*\psi_{2^k})(\cdot,(1-s^{\mathbf{a}_1})^\frac{1}{\mathbf{a}_2})\|_{L^2} \, ds.
\end{align*}
Analogous to \eqref{2-2-1-decay_estimate_for_frequency_pieces}, the main $L^2 \times L^2 \to L^1$ estimate in this case reads as 
\begin{align}
\label{NT-2-2-1-decay_estimate_for_frequency_pieces} 
\|\mathcal{A}^{\mathbf{a},NT}_1(F*\psi_{2^j},G*\psi_{2^k})\|_{L^1} \lesssim 2^{\delta(j+k)}\|F\|_{L^2}\|G\|_{L^2}, 
\end{align}
which follows once we know that $\|A^{NT}(F*\psi_{2^j})\|_{L^2}\lesssim (\frac{2^j}{s})^\delta \|F\|_{L^2},$ for some $0<\delta<2n,$ and the same is already known from \cite[Proposition $5.1$]{ghosh2026bilinearsphericalmaximalfunction}. 

\medskip 
With the estimate \eqref{NT-2-2-1-decay_estimate_for_frequency_pieces} at hand, one can essentially repeat the rest of the proof of Theorem \ref{Main-theorem-lacunary-Koranyi-type-for-degenrate} to complete the proof of Theorem \ref{thm:lac-max-NT-degenerate}. 

%%%%%%%%%%%%%%%%%
%%%%%%%%%%%%%%%%%
%%%%%%%%%%%%%%%%%
%%%%%%%%%%%%%%%%%

\subsection{Proof of Theorem \ref{theorem-improving-single-scales}} 
\label{subsec:thm-single-scale-NT-degenrate} 

It suffices to prove the following $L^1(\mathbb{H}^n) \times L^q(\mathbb{H}^n) \to L^1(\mathbb{H}^n)$ bounds which is analogous to Proposition 3.1 of \cite{ghosh2026bilinearsphericalmaximalfunction}, and with this one can complete the proof of Theorem \ref{theorem-improving-single-scales} by repeating the rest of the proof of Theorem 1.1 of \cite{ghosh2026bilinearsphericalmaximalfunction}, and we shall omit that. 

\begin{proposition}
\label{Proposition-improving-L^1-times-L^q-to-L^1}
Let $\frac{2n+1}{2n}< q \leq \infty$ and  $\mathbf{a}_1,\mathbf{a}_2 \geq 1$. Then the single scale averaging operator $\mathcal{A}_1^{\mathbf{a},NT}$ has the following improving property:
\begin{align}
\label{expression-L^1-times-L^q-to-L^1}
\|\mathcal{A}_1^{\mathbf{a},NT}(F,G)\|_{L^1(\mathbb{H}^n)} \lesssim 
\begin{cases}
\|F\|_{L^1(\mathbb{H}^n)} \, \|G\|_{L^q(\mathbb{H}^n)}, & \textup{if } \, \frac{2n+1}{2n}< q \leq \infty, \\ \\ \|F\|_{L^p(\mathbb{H}^n)} \, \|G\|_{L^1(\mathbb{H}^n)}, & \text{if } \, \frac{2n+1}{2n}< p \leq \infty.
\end{cases}
\end{align}
\end{proposition}

\begin{proof}
It is sufficient to prove the first estimate of \eqref{expression-L^1-times-L^q-to-L^1}. Once again, let us take $F, \, G \geq 0$. Then, 
\begin{align*}
\|\mathcal{A}_1^{\mathbf{a},NT}(F,G)\|_{L^1(\mathbb{H}^n)} 
& = \int_{\mathbb{H}^n} \int_{\Sigma^\mathbf{a}} F(g\cdot(z_1,0)^{-1}) \, G(g\cdot (z_2,0)^{-1}) \, d\mu^\mathbf{a}(z_1,z_2) \, dg \\ 
& = \int_{\mathbb{H}^n} F(g) \, \mathcal{I}_{\mathbf a}G(g)\, dg,
\end{align*}
where 
\begin{align}
\label{step1-expression-L^1-times-L^q-to-L^1} \mathcal{I}_{\mathbf a}G(g) = \int_{\Sigma^\mathbf{a}} \tau_g G \left( z_1-z_2,-\frac{1}{2}\Im(z_1 \Bar{z_2}) \right) d\mu^\mathbf{a}(z_1,z_2) .
\end{align}
Clearly, the claimed estimate in \eqref{expression-L^1-times-L^q-to-L^1} shall follow immediately once we prove that 
\begin{align}
\label{step2-expression-L^1-times-L^q-to-L^1} 
\sup_{g \in \mathbb{H}^n} \, \mathcal{I}_{\mathbf a} G(g) \lesssim \|G\|_{L^q(\mathbb{H}^n)},
\end{align}
for every \(\frac{2n+1}{2n}<q\leq\infty\) and \(a_1,a_2\geq 1\), and this is what we shall now show. 

\medskip
Performing the change of variables $(z_1, \, z_2) \to (\frac{u+b}{2}, \, \frac{u-b}{2})$ in \ref{step1-expression-L^1-times-L^q-to-L^1}, we get 
\begin{align*}
\mathcal{I}_{\mathbf a} G(g) & = \int_{S^\mathbf{a}} \tau_g G \left(b,-\frac{1}{4} \Im (b \cdot \Bar{u})\right)\, d\sigma^\mathbf{a}(u,b), 
\end{align*}
where $S^\mathbf{a} = \left\{ (u,b) \in \mathbb{C}^n \times \mathbb{C}^n : H(u,b) := |\frac{u+b}{2}|^{\mathbf{a}_1} + |\frac{u-b}{2}|^{\mathbf{a}_2} - 1 = 0 \right\}$ and $d\sigma^\mathbf{a}$ is the induced measure on the surface $S^\mathbf{a}$. Now, note that 
\[
|\nabla H(u,b)|
=\frac{1}{\sqrt{2}}
\left(
a_1^2\left|\frac{u+b}{2}\right|^{2a_1-2}
+a_2^2\left|\frac{u-b}{2}\right|^{2a_2-2}
\right)^{1/2}, \quad \text{for a.e. } (u,b),
\]
and in view of the formula \eqref{coarea-formula-statement}, we have 
\begin{align*}
\mathcal{I}_{\mathbf a}G(g) &= \int_{\mathbb{R}^{2n}} \int_{\mathbb{R}^{2n}} \tau_g G \left(b,-\frac{1}{4} \Im (b \cdot \Bar{u})\right)\, |\nabla H(u,b)|\, \delta(H(u,b)) \, du \, db \\ 
&= \int_{S^{2n-1}} \int_0^\infty r^{2n-1} \int_{\mathbb{R}^{2n}} \tau_g G \left(rw,-\frac{1}{4} \Im (rw \cdot \Bar{u})\right) |\nabla H(u,rw)|\, \delta(H(u,rw)) \, du \, dr \, dw, 
\end{align*}
where the last step is nothing but the polar-coordinate decomposition in $b$-variables. 

\medskip
Fix $w\in S^{2n-1}$, and choose $U_\omega\in U(n)$ such that $U_\omega (\omega) = e_1$, where  $e_1=(1,0,\ldots,0)\in\mathbb C^n$. Since $U_\omega$ is complex-linear, its real representation commutes with the standard symplectic $2n \times 2n$ matrix
\[
J=\begin{pmatrix}
0&-I_n\\
I_n&0
\end{pmatrix},
\]
and therefore, via the orthogonal transformation $U_\omega$, we get that $\mathcal{I}_{\mathbf a}G(g)$ equals to 
\[
\int_{S^{2n-1}} \int_0^\infty r^{2n-1} \int_{\mathbb{R}^{2n}} \tau_g G \left(rw,-\frac{1}{4} \Im (re_1 \cdot \Bar{u'})\right) |\nabla H(u',re_1)|\, \delta(H(u',re_1)) \, du' \, dr \, dw.
\]

\medskip 
Let us write
\[
u'=(x+iy)e_1+v,
\]
where $x,y\in\mathbb{R}$, and
$v\in\operatorname{span}_{\mathbb{C}}\{e_2,\ldots,e_n\}$ in $\mathbb{C}^n$. Clearly, $\frac{1}{4} \Im(re_1\cdot \Bar{u'}) = - \frac{1}{4}ry$. Now, writing $q = |v|^2$, in terms of the above co-ordinates, we define 
\begin{align*}
\widetilde{H}(x,y,r;q) :&= H(u',re_1) \\ 
&= \left(\frac{\sqrt{(x + r)^2 + y^2 + q}}{2} \right)^{\mathbf{a}_1} + \left(\frac{\sqrt{(x-r)^2 + y^2 +q}}{2} \right)^{\mathbf{a}_2} -1,
\end{align*} 
so that 
\begin{align}
\label{step3-expression-L^1-times-L^q-to-L^1} 
\mathcal{I}_{\mathbf a}G(g)
&= C \int_{S^{2n-1}} \int_0^\infty r^{2n-1}
\int_{\mathbb R^2}
\tau_g G\left(rw,\frac14 ry\right) \int_{S^{2n-3}} \int_0^\infty
|\nabla H((x+iy,\sqrt q\,\theta),re_1)| \\
& \qquad \delta ( \widetilde{H}(x,y,r;q) ) \, q^{n-2}\,dq \,d\theta \,dx\,dy\,dr\,dw, \notag
\end{align}
where we have applied the polar decomposition on $\mathbb{C}^{n-1} = \mathbb{R}^{2n-2}$, namely, with $v = \sqrt{q} \, \theta$, for $q > 0$ and $\theta\in S^{2n-3}$. 

\medskip 
It is clear from the definition that the function $\widetilde{H}$ is strictly increasing in $q>0$. Moreover, $\widetilde{H}(q) \to +\infty$ as $q \to \infty$. Hence, if we consider the partition
of $\mathbb R^2\times(0,\infty) = \Omega_- \cup \Omega_-$, where 
\[
\Omega_- = \{ (x,y,r) \in \mathbb R^2 \times (0,\infty) : \widetilde{H} (x,y,r;0) < 0 \},
\]
and
\[
\Omega_+=\{(x,y,r)\in\mathbb R^2\times(0,\infty):
\widetilde{H} (x,y,r;0) \geq 0\},
\]
we get that as a function of $q$-variable, $\widetilde{H}$ has no root in $(0,\infty)$ for any $(x,y,r) \in \Omega_+$, and therefore this part does not contribute in \eqref{step3-expression-L^1-times-L^q-to-L^1}. On the other hand, for every $(x,y,r) \in \Omega_-$, there exists unique simple root $q_*>0$ of $\widetilde{H}$, which altogether implies that 
\begin{align}
\label{exp-of-I_aG(g)}
\mathcal{I}_{\mathbf a}G(g)&= C \int_{S^{2n-1}} \int_0^\infty r^{2n-1}
\int_{\{(x,y):(x,y,r)\in \Omega_-\}}
\tau_g G\left(rw,\frac{1}{4} ry\right)\\
&\qquad\times
\int_{S^{2n-3}} q_*^{n-2}\,
\frac{|\nabla H((x+iy, \sqrt{q_*} \, \theta),re_1)|}{| \widetilde{H}_q(x,y,q_*,r)|} \, d\theta \,dx\,dy\,dr\,dw. \notag
\end{align}
By simple calculation, one can verify that  for any arbitrary $q>0$, 
\[
|\nabla H((x+iy, \sqrt{q} \, \theta),re_1)| = \frac{1}{\sqrt{2}} \left(a_1^2 A^{2a_1-2} + a_2^2 B^{2a_2- 2}\right)^{1/2}, 
\] 
and 
\[
| \widetilde{H}_q (x,y,r;q) | = \frac{a_1}{8}A^{a_1-2} + \frac{a_2}{8}B^{a_2-2}, 
\]
where 
\[
A=\frac{\sqrt{(x+r)^2+y^2+q}}{2} \quad \text{and} \quad
B=\frac{\sqrt{(x-r)^2+y^2+q}}{2}. 
\] 
Observe also that within the domain of the integral in \eqref{exp-of-I_aG(g)}, $A$ and $B$ are uniformly bounded quantities, and thus the above relations lead to 
\[
\sup_{x,y,r,q,\theta} \frac{|\nabla H((x+iy, \sqrt{q} \, \theta),re_1)|}{| \widetilde{H}_q (x,y,r;q) |} \lesssim 
\sup_{x,y,r,q} \frac{\left(a_1^2 A^{2a_1-2} + a_2^2 B^{2a_2- 2}\right)^{1/2}}{a_1 A^{a_1-2} + a_2 B^{a_2-2}} 
\leq C_{\mathbf{a}_1,\mathbf{a_2}} < \infty.
\]
Moreover, $H(x,y,q_*,r)=0$ implies that $q_*$ is also uniformly bounded over the domain of the integral in \eqref{exp-of-I_aG(g)}. Using these observations in \eqref{exp-of-I_aG(g)}, we get that 
\[
\mathcal{I}_{\mathbf a} G(g) \lesssim  \int_{S^{2n-1}} \int_0^\infty r^{2n-1}
\int_{\{(x,y):(x,y,r)\in \Omega_-\}}
\tau_g G \left(rw,\frac{1}{4} ry \right) \,dx\,dy\,dr\,dw.
\]
Now, note that for any $(x,y)\in \Omega_-$, we have $r \leq 2$, and $|x|, \, |y| \leq 4$, and therefore 
\begin{align*}
\mathcal{I}_{\mathbf a} G(g) 
& \lesssim \int_{S^{2n-1}} \int_0^2 r^{2n-1}
\int_{|y| \leq 4}
\tau_g G \left(rw, \frac{1}{4} ry \right) \,dy \,dr\,dw \\ 
& = 4 \int_{S^{2n-1}} \int_0^2 r^{2n-1}
\int_{|t| \leq r}
\frac{\tau_g G \left(rw, t \right)}{r} \, dt \,dr\,dw \\ 
& = 4 \int_{|z|\leq 2} \int_{|t|\leq |z|}\frac{\tau_g G \left(z,t \right)}{|z|} \,dt\,dz.
\end{align*}
With the above pointwise inequality at hand, we can make use of the H\"older's inequality to conclude that 
\[
\mathcal{I}_{\mathbf a} G(g) \lesssim \| \tau_g G\|_{L^q(\mathbb{H}^n)} \left(\int_{|z|\leq 2} \int_{|t|\leq |z|}\frac{1}{|z|^{q'}} \,dt\,dz. \right)^\frac{1}{q'} \lesssim \| \tau_g G\|_{L^q(\mathbb{H}^n)},
\] 
provided that $q'-1-2n+1<1$ which is same as $q>\frac{2n+1}{2n}.$ This completes the claimed estimate in \eqref{step2-expression-L^1-times-L^q-to-L^1} and with that the proof of Proposition \ref{Proposition-improving-L^1-times-L^q-to-L^1}.
\end{proof}

%%%%%%%%%%%%%%%%%
%%%%%%%%%%%%%%%%%
%%%%%%%%%%%%%%%%%
%%%%%%%%%%%%%%%%%

\section{On sharpness of Theorems \ref{Main-theorem-Koranyi-type-for-degenrate} and \ref{thm:full-max-NT-degenerate}}
\label{sec:sharpness-examples}
In this section, we shall show that the regions for the strong-type boundedness in estimates \eqref{estimate_for_the_Product_Koranyi_full_degenrate} in Theorem \ref{Main-theorem-Koranyi-type-for-degenrate} and \eqref{estimate_for_the_Nevo_Thangavelu_full_degenrate} in Theorem \ref{thm:full-max-NT-degenerate} are sharp.

\subsection{Sharpness of estimate \texorpdfstring{
\eqref{estimate_for_the_Product_Koranyi_full_degenrate}}{} in Theorem
\ref{Main-theorem-Koranyi-type-for-degenrate}
}
\label{subsec:sharpness-examples-Koranyi}

Here we shall show that estimate \eqref{estimate_for_the_Product_Koranyi_full_degenrate} implies 
\begin{align}
\frac{1}{r} &= \frac{1}{p}+\frac{1}{q}\leq \frac{4n+1}{2n+1}, 
\label{sharpness-1st-product-Koranyi} \\ 
\frac{1}{p} &\leq 1- \frac{1}{2n+1} +\frac{2n+2}{(2n+1)\mathbf{a}_2}, \quad  \mathbf{a}_2>2n+2, 
\label{sharpness-2nd-product-Koranyi} \\ 
\frac{1}{q} &\leq 1- \frac{1}{2n+1} +\frac{2n+2}{(2n+1)\mathbf{a}_1}, \quad  \mathbf{a}_1>2n+2. 
\label{sharpness-3rd-product-Koranyi}
\end{align}
We prove the above claims by constructing Knapp-type examples. 
The same is motivated by the construction in \cite{Rajula-Srivastava-Koranyi_sphere, Lee_Shuin}. 

\medskip 
\textbf{\underline{Proof of \eqref{sharpness-1st-product-Koranyi}}:} Consider 
$h: [0, \infty
) \to \mathbb{R}$ defined by $h(s)=s^{\mathbf{a}_1}+s^{\mathbf{a}_2}-1$. It is easy to see that there exists unique $s_*\in (0,1)$ so that $h(s_*)=0.$ Now, choose and fix $\varepsilon_0>0$ small enough such that $s_*+\varepsilon_0<1.$ 

\medskip
Let $0<\delta \ll \varepsilon_0$ be a small quantity and let $B_\delta = \{(z,t)\in \mathbb{H}^n: |z|\leq \delta , |t|\leq\delta\}$. Clearly, the Lebesgue measure satisfies $|B_\delta| \sim \delta^{2n+1}$. For $c_1, c_2>0$ that we shall choose shortly (see the paragraphs right after estimates \eqref{slicing-argument-lower-bound-first-example}), let $F_\delta$ and $G_\delta$ be the characteristic functions of dilates of $B_{\delta}$, To be more precise, let $F_\delta=\chi_{B_{c_1\delta}}$ and $G_\delta=\chi_{B_{c_2\delta}}$. Then, 
\begin{align*}
& \|F_\delta\|_p^p \sim \|G_\delta\|_q^q \sim \delta^{2n+1}.
\end{align*}

Now, consider the evaluation region
\begin{align} 
\label{region-R}
R:=\{g \in \mathbb{H}^n : s_* < |g|_K < s_*+\varepsilon_0
\}.
\end{align}
Clearly, $|R| \sim 1$. 

\medskip 
Given any $g=(z,t)\in R$, consider the following caps on the Kor\'anyi ball and Kor\'anyi sphere: 
\begin{align*}
E^1_{g} &= \left\{ g_1=(z_1,t_1)\in B_K\left(0,1\right):
\left| z_1-\frac{z\,s_*}{|g|_K} \right|<\delta, \,
\left|t_1- \frac{t\,s_*^2}{|g|_K^2} \right|<\delta
\right\}, \\ 
E^2_{g} &= \left\{ g_2=(z_2,t_2)\in S_{K,1}:
\left| z_2-\frac{z}{|g|_K} \right|<\delta,\;
\left| t_2-\frac{t}{|g|_K^2} \right|<\delta
\right\}.
\end{align*}
It is not difficult to verify that the Lebesgue measure of $E^1_{g}$ is equivalent to $\delta^{2n+1}$ and the surface measure of $E^2_{g}$ is equivalent to $\delta^{2n}$. 

\medskip
Now, applying the corresponding slicing argument we developed earlier in \eqref{equivalent expression-degenrate-Koranyi-2}, we observe that for any $g\in R$,
\begin{align}
\label{slicing-argument-lower-bound-first-example} 
& \mathcal{M}^{\mathbf{a},K}_{full}(F_\delta,G_\delta)(g) \geq \mathcal{A}^{\mathbf a, K}_{\frac{|g|_K}{s_*}}(F_\delta, G_\delta)(g)\\
&= \int_{ \mathcal{S}^\mathbf{a}_K} F_\delta (g \cdot \tilde{\delta}_{s_*^{-1}|g|_K}{g_1}^{-1}) \, G_\delta (g \cdot \tilde{\delta}_{s_*^{-1}|g|_K}{g_2}^{-1})\, d\mu^{\mathbf{a}}(g_1,g_2) \notag \\
& \gtrsim \int_{E^1_g} F_\delta (g\cdot \tilde{\delta}_{s_*^{-1}|g|_K} {g_1}^{-1}) \int_{E^2_g} G_\delta (g\cdot \tilde{\delta}_{s_*^{-1}|g|_K (1-|g_1|_K^{\mathbf{a}_1})^{\frac{1}{\mathbf{a}_2}}}{g_2}^{-1}) \, d\sigma^{2n}_K(g_2) \, (1-|g_1|_K^{\mathbf{a}_1})^{\frac{2n+2-\mathbf{a}_2}{\mathbf{a}_2}} \,dg_1. \notag 
\end{align}

\medskip
We shall show that by choosing $\delta>0$ small enough, one can ensure that $w_\mathbf{a}(g_1) = \left( 1-|g_1|_K^{\mathbf{a}_1} \right)^{\frac{1}{\mathbf{a}_2}} \gtrsim 1$ for all $g \in R$. Moreover, in this case, we can also choose $c_1, c_2 > 0$ in such a way that the functions $F_\delta$ and $G_\delta$ in the above integral evaluate to exactly $1$. 

\medskip 
Let us assume the above claims for now and see how they are useful in completing the proof of claim 
\eqref{sharpness-1st-product-Koranyi}. It follows from \eqref{slicing-argument-lower-bound-first-example} that on the set $R$, we have 
\[
\mathcal{M}^K_{full}(F_\delta, G_\delta)(g) \gtrsim \delta^{2n+1} \delta^{2n} = \delta^{4n+1},
\]
which implies that 
\begin{align}
\label{lower-bound-first-example}
\|\mathcal{M}^K_{full}\|_{L^r(\mathbb{H}^n)}\gtrsim\delta^{4n+1}.
\end{align}

Suppose now we have boundedness as in \eqref{estimate_for_the_Product_Koranyi_full_degenrate}, then with $F_\delta$ and $G_\delta$ as above, we can make use of the estimate \eqref{lower-bound-first-example} to conclude that 
\[
\delta^{4n+1}\lesssim \delta^{\frac{2n+1}{p}} \delta^{\frac{2n+1}{q}}, 
\]
which upon letting $\delta \rightarrow 0^+$ forces 
\[
4n+1-\left(\frac{2n+1}{p} + \frac{2n+1}{q}\right)\geq 0,
\]
which is same as saying that $\frac{1}{r} = \frac{1}{p}+\frac{1}{q}\leq \frac{4n+1}{2n+1},$ and this completes the proof of claim 
\eqref{sharpness-1st-product-Koranyi}. 

\medskip 
We are thus left with proving the claims written in the paragraph right after \eqref{slicing-argument-lower-bound-first-example} which we take up now. Note that with $g=(z,t)\in R$ and $g_1=(z_1,t_1)\in E^1_{g}$, it is not difficult to verify that 
$$ \left| z-\frac{|g|_K}{s_*}z_1 \right| = \frac{|g|_K}{s_*} \left| z_1-\frac{z\,s_*}{|g|_K} \right| \leq d_1 \delta,$$ 
for some $d_1>0,$ and 
$$\left|t-\frac{|g|_K}{2s_*}\Im(z\overline{z_1})-\frac{|g|_K^2}{s_*^2}t_1\right|\leq\frac{|g|_K^2}{s_*^2}\left|t_1-\frac{t\,s_*^2}{|g|_K^2}\right|+\frac{|g|_K}{2s_*}|\Im(z\overline{z_1})| =  d_1^2 \delta + \frac{1}{2} |\Im ((z-\frac{|g|_K}{s_*}z_1)\overline{z})|\leq d_2\delta,$$ 
for some $d_2>0. $ So, if we choose $c_1 = \max \{d_1, d_2\}$, then it follows that $F_\delta (g \cdot \tilde{\delta}_{ s_*^{-1} |g|_K}{g_1}^{-1}) = 1$.

\medskip
On the other hand, let us write $\lambda=\frac{s_*}{|g|_K}$, then $|(\lambda z,\lambda^2t)|_K=s_*,$ and
\begin{align*}
    ||g_1|_K^4-s_*^4|&\leq||z_1|^4-|\lambda z|^4|+|t_1^2-\lambda^4t^2|\\
    &\lesssim  |z_1-\lambda z|+ |t_1-\lambda^2t|\leq d_3\delta, \quad \text{ for some } d_3>0,
\end{align*}
which, in particular, implies that by keeping $\delta>0$ small enough, we have $w_\mathbf{a}(g_1) = \left( 1-|g_1|_K^{\mathbf{a}_1} \right)^{\frac{1}{\mathbf{a}_2}} \gtrsim 1$ and $w_\mathbf{a}(g_1) \simeq s_*.$ It is then easy to verify that for any $g_2=(z_2,t_2)\in E^2_g,$ we have 
$$ \left| z-\frac{|g|_Kw_\mathbf{a}(g_1)}{s_*}z_2 \right| \lesssim |z-|g|_Kz_2|\leq d_4 \delta,$$ 
for some $d_4>0,$ and 
$$ \left| t-\frac{|g|^2_Kw_\mathbf{a}^2(g_1)}{s_*^2}t_2-\frac{|g|_Kw_\mathbf{a}(g_1)}{2s_*}\Im(z\overline{z_2}) \right| \lesssim |t-|g|_K^2t_2-\frac{|g|_K}{2}\Im(z\overline{z_2})|\leq d_5\delta,$$ for some $d_5>0.$ Hence, if we take $c_2=\max\{d_4,d_5\}$, then $ G_\delta (x \cdot \tilde{\delta}_{s_*^{-1} |g|_K\,(1-|g_1|_K^{\mathbf{a}_1})^{\frac{1}{\mathbf{a}_2}}}{g_2}^{-1}) = 1$. 

\medskip 
\textbf{\underline{Proof of \eqref{sharpness-2nd-product-Koranyi}}:} For $0<\delta<\frac{1}{2}$, consider $F_\delta=\chi_{B_{c \delta}}$, for some constant $c>0$ which can be chosen appropriately (following the idea similar to that in the proof of \eqref{sharpness-1st-product-Koranyi}). Also, let $G=\chi_{B_{10}}.$ Then $\|F_\delta\|_{L^p(\mathbb{H}^n)}\sim\delta^\frac{2n+1}{p}$ and  $ \|G\|_{L^q(\mathbb{H}^n)} \gtrsim 1,$ for any $1\leq p,q \leq\infty.$ 

\medskip 
Given $g\in R = \{1 \leq |g|_K \leq2\}$, it follows from \eqref{equivalent expression-degenrate-Koranyi} that 
\begin{align*}
& \mathcal{M}^{\mathbf{a},K}_{full}(F_\delta,G)(g) \\ 
& \gtrsim \int \limits_{|g_2|^{\mathbf{a}_2}_K\leq\delta} G (g\cdot \tilde{\delta}_{|g|_K}{g_2}^{-1}) \int_{E^2_g} F_\delta (g \cdot \tilde{\delta}_{|g|_K\,(1-|g_2|_K^{\mathbf{a}_2})^{\frac{1}{\mathbf{a}_1}}}{g_1}^{-1}) \, d\sigma^{2n}_K(g_1) \, (1-|g_2|_K^{\mathbf{a}_2})^{\frac{2n+2-\mathbf{a}_1}{\mathbf{a}_1}} \, dg_2.
\end{align*}
Since $|g_2|_K^{\mathbf{a}_2} \leq \delta$, by choosing $\delta>0$ to be small enough, we can ensure that $\tilde{w}_\mathbf{a}(g_2)\sim1$ and $|\tilde{w}_\mathbf{a}(g_2)-1|\lesssim \delta,$ and therefore, for any $g\in R$, we have
\[
\mathcal{M}^{\mathbf{a},K}_{full}(F_\delta,G)(g)\gtrsim\delta^{\frac{2n+2}{\mathbf{a}_2}+2n}. 
\] 
As an immediate consequence of the above estimate, we get that if $\mathcal{M}^{\mathbf{a},K}_{full}$ satisfies the inequality \eqref{Main-theorem-Koranyi-type-for-degenrate}, then we must have $\delta^{\frac{2n+2}{\mathbf{a}_2}+2n+1-1}\lesssim \delta^\frac{2n+1}{p}$ for small $\delta > 0$. With that, we obtain $\frac{2n+2}{\mathbf{a}_2}+2n+1-1\ge \frac{2n+1}{p}$ by letting $\delta\to 0^+$, and hence we obtain \eqref{sharpness-2nd-product-Koranyi}. In particular, when $\mathbf{a}_2>2n+2$, this gives a nontrivial necessary condition on $p$.

\medskip 
Finally, one can observe that Claim \eqref{sharpness-3rd-product-Koranyi} is symmetric to Claim \eqref{sharpness-2nd-product-Koranyi}. 

\subsection{Sharpness of estimate 
\texorpdfstring{
\eqref{estimate_for_the_Nevo_Thangavelu_full_degenrate}
}{} in Theorem \ref{thm:full-max-NT-degenerate}} 
\label{subsec:sharpness-examples-NT}

Here we shall show that estimate \eqref{estimate_for_the_Nevo_Thangavelu_full_degenrate} implies 
\begin{align}
\frac{1}{r} &= \frac{1}{p}+\frac{1}{q}\leq \frac{4n-1}{2n}, 
\label{sharpness-1st-product-NT} \\ 
\frac{1}{p} &\leq 1- \frac{1}{2n} +\frac{1}{\mathbf{a}_2}, \quad  \mathbf{a}_2>2n, 
\label{sharpness-2nd-product-NT} \\ 
\frac{1}{q} &\leq 1- \frac{1}{2n} +\frac{1}{\mathbf{a}_1}, \quad  \mathbf{a}_1>2n. 
\label{sharpness-3rd-product-NT}
\end{align}

Once again, we shall construct some Knapp-type examples for the above claims. These examples are motivated from \cite{ghosh2026bilinearsphericalmaximalfunction,Lee_Shuin}. Mostly, the arguments here are on the lines of that from Subsection \ref{subsec:sharpness-examples-Koranyi}, but we write the necessary details for the sake of completion. 

\medskip 
\textbf{\underline{Proof of \eqref{sharpness-1st-product-NT}}:} We shall choose the evaluation region $R_\delta$ in such a way that within it $w^{NT}_\mathbf{a}(z_1)=(1-|z_1|^{\mathbf{a}_1})^\frac{1}{\mathbf{a}_2}$ remains away from $0$ and $1$. Taking the function $h$ same as in the proof of \eqref{sharpness-1st-product-Koranyi}, we have unique $s_*\in (0,1)$ such that $h(s_*)=0.$ 
Fix $\epsilon_0 > 0$ so small such that $s_* + \epsilon_0 < 1.$

\medskip
Let $0<\delta\ll \epsilon_0$ be very small. With that, consider two functions $F_\delta$ and $G_\delta$ be the characteristic functions of dilates of $B_{\delta}$ (Euclidean ball in $\mathbb{H}^n$), that is, let $F_\delta=\chi_{B_{c_1\delta}}$ and $G_\delta=\chi_{B_{c_2\delta}}$, for appropriate $c_1, \, c_2>0$. Then, 
\begin{align*}
& \|F_\delta\|_p^p \sim \|G_\delta\|_q^q \sim \delta^{2n+1}.
\end{align*}

Now, define the evaluation region $R_\delta$ as follows: 
\[
R_\delta:=\{g=(z,t) \in \mathbb{H}^n : s_* < |z| < s_*+\varepsilon_0 , \, |t|<\delta
\}. 
\] 
Clearly, $|R| \sim \delta$. 
Given any $g=(z,t)\in R_\delta$, consider the following caps on the Euclidean ball and sphere: 
\begin{align*}
D^1_{g} = \left\{ z_1\in B^{2n}\left(0,1\right): \left| \hat{z}-\frac{z_1}{s_*} \right|<\delta \right\} \quad \text{and} \quad  
D^2_{g} = \left\{ z_2\in S^{2n-1}: \left| \hat{z}-z_2 \right|<\delta
\right\}.
\end{align*}
It is easy to verify that the Lebesgue measure of $D^1_{g}$ is equivalent to $\delta^{2n}$ and the surface measure of $D^2_{g}$ is equivalent to $\delta^{2n-1}$. 

\medskip
Now, for any $g\in R_\delta $ and with $s_1=s_2=\frac{|z|}{s_*}$, using the slicing argument from \eqref{equivalent-expresion-1-nevo-type}, we have 
\begin{align*}
& \mathcal{M}^{\mathbf{a},NT}_{full}(F_\delta,G_\delta)(g) \\ 
& \gtrsim \int_{D^1_g} F_\delta (g\cdot \tilde{\delta}_{s_*^{-1}|z|} {(z_1,0)}^{-1}) \int_{D^2_g} G_\delta (g\cdot \tilde{\delta}_{s_*^{-1}|z| w^{NT}_\mathbf{a}(z_1)}{(z_2,0)}^{-1}) \, d\sigma(z_2) \,  \left(w^{NT}_\mathbf{a}(z_1)\right)^{2n-\mathbf{a}_2} \,dz_1. \notag 
\end{align*}

With that, we can repeat the last part of the proof of Claim \eqref{sharpness-1st-product-Koranyi} to argue that one can choose $c_1, \, c_2 > 0$ such that for small enough $\delta > 0$, functions $F_\delta$ and $G_\delta$ evaluate to exactly $1$ and $w^{NT}_\mathbf{a}(z_1) \simeq s_*$ in the right hand side of the above integral. This would imply that 
\[
\mathcal{M}^{NT}_{full}(F_\delta, G_\delta)(g) \gtrsim \delta^{2n} \delta^{2n-1} = \delta^{4n-1},
\]
for small enough $\delta > 0$ and for all $g\in R_\delta
$. Consequently, if we have the boundedness \eqref{estimate_for_the_Nevo_Thangavelu_full_degenrate}, then we must have 
\[
\delta^{4n-1}\delta^\frac{1}{r}\lesssim \delta^{\frac{2n+1}{p}} \delta^{\frac{2n+1}{q}}. 
\]
for small enough $\delta > 0$. Hence, letting $\delta \rightarrow 0^+$ forces 
\[
\frac{1}{r} =\frac{1}{p}+\frac{1}{q} \leq \frac{4n-1}{2n}
\]
and this completes the proof of claim 
 \eqref{sharpness-1st-product-NT}. 

 \medskip
\textbf{\underline{Proof of \eqref{sharpness-2nd-product-NT}}:} The proof follows by repeating the proof of \eqref{sharpness-2nd-product-Koranyi} with only the following in mind. As earlier, we only consider small enough $\delta > 0$, and with that, owing to the definition of the bilinear maximal function, it is natural to explore the estimates with the characteristic function $F_\delta$ corresponding to the set with small region for $z$-variable, and fixed bounded range for $t$-variable. More precisely, we consider the set $U_\delta = \{(z,t)\in \mathbb{H}^n: |z| < c \, \delta, \, |t|<10\}$ for an appropriate constant $c>0$ and define $F_\delta = \chi_{U_\delta}$. Also, take $G=\chi_{B_{10}}.$ With these functions, one can essentially repeat the proof of \eqref{sharpness-2nd-product-Koranyi} to complete the proof of \eqref{sharpness-2nd-product-NT}. 

%%%%%%%%%%%%%%%%%%
%%%%%%%%%%%%%%%%%%
%%%%%%%%%%%%%%%%%%
%%%%%%%%%%%%%%%%%% 

\section*{Acknowledgments}
SB and RG were partially supported by the Anusandhan National Research Foundation (ANRF), India, under the research project ANRF/ARG/2025/003732/MS. SKR gratefully acknowledges the Council of Scientific and Industrial Research (CSIR) for the award of a Senior Research Fellowship (SRF) under file no. 09/1020(22850)/2025-EMR-I at the Indian Institute of Science Education and Research (IISER) Bhopal. 

%%%%%%%%%%%%%%%%%
%%%%%%%%%%%%%%%%%
%%%%%%%%%%%%%%%%%
%%%%%%%%%%%%%%%%%

%\bibliographystyle{amsalpha}
%\bibliography{mybibliogra}

\begin{thebibliography}{GSHW24}

\bibitem[BF24]{Borges-Tainara-Foster-bounds-for-lacunary-bilinear-spherical-maximal-functions}
Tainara Borges and Benjamin Foster, \emph{Bounds for lacunary bilinear spherical and triangle maximal functions}, J. Fourier Anal. Appl. \textbf{30} (2024), no.~5, Paper No. 55, 29. \MR{4800832}

\bibitem[BGH{\etalchar{+}}18]{Loukas-Grafakos-bilinear-spherical-maximal-function}
J.~A. Barrionuevo, Loukas Grafakos, Danqing He, Petr Honz\'ik, and Lucas Oliveira, \emph{Bilinear spherical maximal function}, Math. Res. Lett. \textbf{25} (2018), no.~5, 1369--1388. \MR{3917731}

\bibitem[BHRT21]{Sayan-Sourav-Luz-Thangavelu-on-the-maximal-function-associated-to-the-spherical-means-on-the-Heisenberg-groups}
Sayan Bagchi, Sourav Hait, Luz Roncal, and Sundaram Thangavelu, \emph{On the maximal function associated to the spherical means on the {H}eisenberg group}, New York J. Math. \textbf{27} (2021), 631--675. \MR{4250270}

\bibitem[Bou86]{Bourgain-Averages-in-the-plane-over-convex-curves-and-maximal-operators}
J.~Bourgain, \emph{Averages in the plane over convex curves and maximal operators}, J. Analyse Math. \textbf{47} (1986), 69--85. \MR{874045}

\bibitem[Cal79]{Calderon-lacunary-spherical-means}
Calixto~P. Calder\'on, \emph{Lacunary spherical means}, Illinois J. Math. \textbf{23} (1979), no.~3, 476--484. \MR{537803}

\bibitem[Cow79]{Cowling-Littlewood-paley-stein-theory}
Michael~G. Cowling, \emph{On {L}ittlewood--{P}aley--{S}tein theory}, Rendiconti del Seminario Matematico e Fisico di Milano \textbf{49} (1979), no.~1, 79--87.

\bibitem[GGI{\etalchar{+}}13]{Geba-Greenleaf-multilinear-operators}
Dan-Andrei Geba, Allan Greenleaf, Alex Iosevich, Eyvindur Palsson, and Eric Sawyer, \emph{Restricted convolution inequalities, multilinear operators and applications}, Math. Res. Lett. \textbf{20} (2013), no.~4, 675--694. \MR{3188026}

\bibitem[GS26]{ghosh2026bilinearsphericalmaximalfunction}
Abhishek Ghosh and Rajesh~K. Singh, \emph{Bilinear spherical maximal function on the heisenberg group}, 2026.

\bibitem[GSHW24]{Aswin-Hickman}
Aswin Govindan~Sheri, Jonathan Hickman, and James Wright, \emph{Lacunary maximal functions on homogeneous groups}, J. Funct. Anal. \textbf{286} (2024), no.~3, Paper No. 110250, 25. \MR{4669588}

\bibitem[GT21]{On-the-lacubary-spherical-maximal-function-on-the-Heisenberg-group}
Pritam Ganguly and Sundaram Thangavelu, \emph{On the lacunary spherical maximal function on the {H}eisenberg group}, J. Funct. Anal. \textbf{280} (2021), no.~3, Paper No. 108832, 32. \MR{4170795}

\bibitem[H\"83]{The-analysis-of-linear-partial-differential}
Lars H\"ormander, \emph{The analysis of linear partial differential operators. {I}}, Grundlehren der mathematischen Wissenschaften [Fundamental Principles of Mathematical Sciences], vol. 256, Springer-Verlag, Berlin, 1983, Distribution theory and Fourier analysis. \MR{717035}

\bibitem[HHY20]{Heo-Hong-Improved-bounds-for-the-bilinear-spherical-maximal-function}
Yaryong Heo, Sunggeum Hong, and Chan~Woo Yang, \emph{Improved bounds for the bilinear spherical maximal operators}, Math. Res. Lett. \textbf{27} (2020), no.~2, 397--434. \MR{4117081}

\bibitem[IPS22]{Iosevich}
Alex Iosevich, Eyvindur~Ari Palsson, and Sean~R. Sovine, \emph{Simplex averaging operators: quasi-{B}anach and {$L^p$}-improving bounds in lower dimensions}, J. Geom. Anal. \textbf{32} (2022), no.~3, Paper No. 87, 16. \MR{4363760}

\bibitem[JL20]{Maximal-estimates-for-the-bilinear-spherical-averages}
Eunhee Jeong and Sanghyuk Lee, \emph{Maximal estimates for the bilinear spherical averages and the bilinear {B}ochner-{R}iesz operators}, J. Funct. Anal. \textbf{279} (2020), no.~7, 108629, 29. \MR{4103874}

\bibitem[Lac19]{Lacey-Sparse-bounds-for-spherical-maximal-functions}
Michael~T. Lacey, \emph{Sparse bounds for spherical maximal functions}, J. Anal. Math. \textbf{139} (2019), no.~2, 613--635. \MR{4041115}

\bibitem[LQ14]{Li-Hong-Qian-Hardy-Littlewood-maximal-function-on-the-Heisenberg-group}
Hong-Quan Li and Bin Qian, \emph{Centered {H}ardy-{L}ittlewood maximal functions on {H}eisenberg type groups}, Trans. Amer. Math. Soc. \textbf{366} (2014), no.~3, 1497--1524. \MR{3145740}

\bibitem[LS23]{Lee_Shuin}
Sanghyuk Lee and Kalachand Shuin, \emph{Bilinear maximal functions associated with degenerate surfaces}, J. Funct. Anal. \textbf{285} (2023), no.~8, Paper No. 110070, 26. \MR{4609073}

\bibitem[MS04]{Muller-Seeger-spherical-maximal-function-on-twostep-nilpotent-group}
Detlef M\"uller and Andreas Seeger, \emph{Singular spherical maximal operators on a class of two step nilpotent {L}ie groups}, Israel J. Math. \textbf{141} (2004), 315--340. \MR{2063040}

\bibitem[NT97]{Nevo-Thangavelu-pointwise-ergodic-theorems-for-the-radial-averages-on-the-Heisenberg-groups}
Amos Nevo and Sundaram Thangavelu, \emph{Pointwise ergodic theorems for radial averages on the {H}eisenberg group}, Adv. Math. \textbf{127} (1997), no.~2, 307--334. \MR{1448717}

\bibitem[NT04]{Naryanan-Thangavelu-optimal-theorem-for-the-spherical-maximal-operator-on-the-Heisenberg-group}
E.~K. Narayanan and S.~Thangavelu, \emph{An optimal theorem for the spherical maximal operator on the {H}eisenberg group}, Israel J. Math. \textbf{144} (2004), 211--219. \MR{2121541}

\bibitem[RSS22]{Roos-seeger-Rajula-lebesgue-space-estimates-for-spherical-maximal-functions}
Joris Roos, Andreas Seeger, and Rajula Srivastava, \emph{Lebesgue space estimates for spherical maximal functions on {H}eisenberg groups}, Int. Math. Res. Not. IMRN (2022), no.~24, 19222--19257. \MR{4523247}

\bibitem[RSS23]{Roos-seeger-Rajula-Restricted-dilation-sets}
\bysame, \emph{Spherical maximal operators on {H}eisenberg groups: restricted dilation sets}, Studia Math. \textbf{273} (2023), no.~1, 1--28. \MR{4660419}

\bibitem[Sri24]{Rajula-Srivastava-Koranyi_sphere}
Rajula Srivastava, \emph{On the {K}or\'anyi spherical maximal function on {H}eisenberg groups}, Math. Ann. \textbf{388} (2024), no.~1, 191--247. \MR{4693933}

\bibitem[Ste76]{Stein-Elias-Maximal-functions-Spherical-means}
Elias~M. Stein, \emph{Maximal functions. {I}. {S}pherical means}, Proc. Nat. Acad. Sci. U.S.A. \textbf{73} (1976), no.~7, 2174--2175. \MR{420116}

\bibitem[Tha98]{Harmonic-analysis-on-the-Heisenberg-Group-Thangavelu}
Sundaram Thangavelu, \emph{Harmonic analysis on the {H}eisenberg group}, Progress in Mathematics, vol. 159, Birkh\"auser Boston, Inc., Boston, MA, 1998. \MR{1633042}

\end{thebibliography}
\newcommand{\etalchar}[1]{$^{#1}$}
\providecommand{\bysame}{\leavevmode\hbox to3em{\hrulefill}\thinspace}
\providecommand{\MR}{\relax\ifhmode\unskip\space\fi MR }
\providecommand{\MRhref}[2]{%
  \href{http://www.ams.org/mathscinet-getitem?mr=#1}{#2}
}
\providecommand{\href}[2]{#2}

\end{document}